\documentclass{amsart}

\usepackage{palatino, epsfig, amsmath, amsthm, amssymb, amscd, mathrsfs}
\input epsf
\def\relabelbox{%
  \hbox\bgroup%
}%
\def\endrelabelbox{%
\egroup%
}%
\def\relabel #1#2 {%
  \special{ps:/a {} def}%
  \smash{\rlap{#2}}%
}%
\def\adjustrelabel <#1,#2> #3#4 {%
  \special{ps:/a {} def}%
  \smash{\rlap{\kern #1 \raise #2\hbox{#4}}}%
}%
\def\extralabel <#1,#2> #3 {\smash{\rlap{\kern #1 \raise #2\hbox{#3}}}}%

\begin{document}

\newtheorem{thm}{Theorem}[subsection]
\newtheorem{lem}[thm]{Lemma}
\newtheorem{cor}[thm]{Corollary}
\newtheorem{conj}[thm]{Conjecture}
\newtheorem{qn}[thm]{Question}
\newtheorem*{claim}{Claim}
\newtheorem*{crit}{Criterion}
\newtheorem*{lam_exist_thm}{Theorem A}

\theoremstyle{definition}
\newtheorem{defn}[thm]{Definition}
\newtheorem{construct}[thm]{Construction}
\newtheorem{note}[thm]{Notation}

\theoremstyle{remark}
\newtheorem{rmk}[thm]{Remark}
\newtheorem{exa}[thm]{Example}

\def\R{\mathbb R}
\def\Z{\mathbb Z}
\def\CP{\mathbb {CP}}
\def\H{\mathbb H}
\def\F{\mathscr F}
\def\E{\mathscr E}
\def\HH{\mathscr H}
\def\M{\mathscr M}
\def\C{\mathbb C}
\def\Q{\mathscr Q}
\def\P{\mathscr P}
\def\SS{\mathscr S}
\def\T{\mathscr T}
\def\G{\mathscr G}
\def\L{\mathscr L}
\def\fc{\mathfrak c}
\def\fG{\mathfrak G}
\def\PSL{\text{PSL}}
\def\homeo{\text{Homeo}}
\def\isom{\text{Isom}}
\def\core{\text{core}}
\def\u{{\text{univ}}}
\def\geo{{\text{geo}}}
\def\split{{\text{split}}}
\def\pos{{\text{pos}}}
\def\neg{{\text{neg}}}
\def\stb{{\text{stab}}}
\def\sli{{\text{sl}}}
\def\til{\widetilde}
\def\I{\mathscr I}
\def\N{\mathbb N}

\title{Promoting essential laminations}

\author{Danny Calegari}

\address{Department of Mathematics \\ California Institute of Technology \\ 
Pasadena CA, 91125}
\email{dannyc@its.caltech.edu}
\date{4/11/2006. Version 0.30}

\begin{abstract}
We show that a co--orientable taut foliation of a closed, orientable, algebraically
atoroidal $3$--manifold is either the weak stable foliation of an Anosov flow,
or else there are a pair of very full genuine laminations transverse to the foliation.
\end{abstract}


\maketitle

\section{Introduction}\label{introduction}

A topological manifold is a very flabby object. It has no local
internal structure, and except in very special cases,
the group of automorphisms is
transitive on the set of subsets of a fixed finite cardinality.
The same manifold can appear in a myriad of different forms, and the
question of recognizing or distinguishing manifolds, or of certifying
a useful property, is in the best case very hard, and in the worst
(typical) case algorithmically unsolvable.

It is therefore desirable to {\em stiffen} or {\em rigidify} the
structure of a manifold, by introducing geometry in some form or other, in
order to reduce this ambiguity of form to a manageable amount. But
exactly what sort of geometric constraints are neither too much
(so that there are no examples) or too little (so that the geometric
structure does not help with the problem of understanding or recognizing
the underlying manifold) is very dimension dependent. As a general
rule, smaller dimensional objects are easier to understand. Important
principles become easier to apply and yield
more and richer structure. In this paper, the notion of {\em monotonicity} is very
important, especially as it relates to natural order or partial-order structures on certain sets.
More generally, such order structures provide a bridge from geometric
problems to algebraic language, and permit one to perform experiments
and construct certificates with the use of computers.

As an organizational tool, monotonicity loses effectiveness
as dimension goes up;
consequently it is most powerful when used in the context of certain
dynamical systems, which can effectively
reduce the study of a manifold to two complementary problems of
{\em strictly smaller dimension}: the
study of the {\em orbits} of the system, and the study of the
{\em parameter space} or {\em leaf space} 
of the orbits. Geometric or analytic qualities
of the dynamical system are reflected in the properties of the dimensionally
reduced systems. 

\subsection{Foliations and arboreal group theory}

In the case of the study of $3$--manifolds, a very effective
tool for dimensional reduction is the structure of a $2$--dimensional foliation,
especially a {\em taut foliation} $\F$ of $M$ which, at least
when $M$ is atoroidal, may be defined as
a $2$--dimensional foliation without spherical or torus leaves. 
A basic structure theorem of
Novikov implies that the leaf space $L$ of the universal cover $\til{\F}$ of
such a foliation is a {\em typically
non--Hausdorff simply--connected $1$--manifold}, on which $\pi_1(M)$ acts
naturally.

Such $1$--manifolds are not unlike $\R$--trees in some ways, and many
of the tools of arboreal group theory (e.g. \cite{Serre_trees}, \cite{Alperin_arboreal})
can be used to study the action of $\pi_1(M)$. If $\F$ is co--orientable,
the leaf space $L$ is an {\em oriented} $1$--manifold, and this orientation
defines a {\em partial order} on the elements of $L$. This global partial 
order structure adds extra nuances to the arboreal theory, and is the
source of many important constructions. For example, in a remarkable {\it tour de force},
Roberts, Shareshian and Stein (\cite{RSS}) recently 
managed to give examples of an infinite family of
hyperbolic $3$--manifolds which do not admit taut foliations, simply by studying the
action of their fundamental groups on (non--Hausdorff) simply--connected $1$--manifolds.

\subsection{The classification of surface homeomorphisms}

If $L$ and the action of $\pi_1(M)$ are understood, it remains to understand
the leaves of $\F$ themselves, and the way they fill out the manifold $M$;
the relevant subject is the theory of surface homeomorphisms.

In this subsection we discuss the simplest case of the theory of
surface homeomorphisms. We derive Thurston's famous theorem on the 
classification of surface homeomorphisms 
by a route which is nearly the opposite of
the historical and, for that matter, the logical direction. The reason
is mainly pedagogical: this order of exposition more clearly reveals the
order of development of some analogous ideas for more general taut foliations.
We will necessarily cover a lot of material very briefly. Most of
the details can be found in the papers \cite{wT88}, \cite{wT98b} and
\cite{Cannon_Thurston}.

First we recall the statement of the theorem, in its most basic form.

\begin{thm}[Thurston, \cite{wT88} Classification of 
surface homeomorphisms]\label{classification_of_homeomorphisms}
Let $\Sigma$ be a closed, orientable surface of genus at least $2$, and
let $\phi:\Sigma \to \Sigma$ be a homeomorphism. Then one of the following
three alternatives holds:
\begin{enumerate}
\item{$\phi$ is {\em periodic}; that is, some finite power of $\phi$ is isotopic
to the identity.}
\item{$\phi$ is {\em reducible}; that is, there is some finite collection of disjoint
essential simple closed curves in $\Sigma$ which are permuted up to isotopy by
$\phi$.}
\item{$\phi$ is {\em pseudo--Anosov}.}
\end{enumerate}
\end{thm}

For the moment, we postpone the definition of a pseudo--Anosov diffeomorphism of
a surface, since this will be the punchline of our revisionist story.

Given the pair $(\Sigma,\phi)$ one forms the mapping torus $M_\phi$ which is
the quotient of the product $\Sigma \times I$ by the equivalence relation
$(s,1) \sim (\phi(s),0)$. $M_\phi$ is a fibration over $S^1$, with fiber
$\Sigma$ and monodromy $\phi$. By analogy with the notation for a short exact
sequence, we denote this
$$\Sigma \to M_\phi \to S^1$$
and there is a corresponding short exact sequence of groups
$$\pi_1(\Sigma) \to \pi_1(M_\phi) \to \Z$$
which represents $\pi_1(M_\phi)$ as an HNN extension. The automorphism $\phi$
of $\Sigma$ induces an automorphism $\phi_*$ of $\pi_1(\Sigma)$, well--defined
up to inner automorphisms. A presentation for $\pi_1(M_\phi)$ is then
given by
$$\pi_1(M_\phi) = \langle \pi_1(\Sigma), t \; 
| \; t^{-1} \alpha t = \phi_*(\alpha) \text{ for each } \alpha \text{ in } \pi_1(\Sigma) 
\rangle$$
The homeomorphism type of this $3$--manifold only
depends on the isotopy class of $\phi$. Then the classification of $\phi$
neatly reflects the geometry of the mapping torus.

Recall that for $M$ a closed, topological $3$--manifold, and $X$ a simply--connected
locally symmetric Riemannian $3$--manifold, an {\em $X$ geometry} on $M$ is a
homeomorphism $$\varphi:M \to X/\Gamma$$ 
where $\Gamma$ is a free, discrete, cocompact, properly discontinuous subgroup
of $\isom(X)$. See \cite{wT_book} for more details.

\begin{thm}[Thurston, \cite{wT98b} Geometrization of 
surface bundles]\label{geometrization_of_bundles}
Let $\Sigma$ be a surface of genus at least $2$, and let
$\phi:\Sigma \to \Sigma$ be a homeomorphism. Then the mapping torus
$M_\phi$ satisfies the following:
\begin{enumerate}
\item{If $\phi$ is periodic, $M_\phi$ admits an $\H^2 \times \R$ geometry.}
\item{If $\phi$ is reducible, $M_\phi$ has a non--trivial JSJ decomposition.}
\item{If $\phi$ is pseudo--Anosov, $M_\phi$ admits an $\H^3$ geometry.}
\end{enumerate}
\end{thm}

From now on we consider the case where $\phi$ is pseudo--Anosov, and
therefore $M_\phi$ is hyperbolic, and we can identify its universal
cover with hyperbolic $3$--space $$\til{M_\phi} = \H^3$$
The action of
$\pi_1(M_\phi)$ on $\til{M_\phi}$ extends continuously to an action
on the {\em ideal boundary} of $\H^3$, which is a topological sphere
which we denote by $S^2_\infty$, and the action of $\pi_1(M_\phi)$ on
this sphere is by M\"obius transformations. We denote the representation 
inducing this action by
$$\rho_{\text{geo}}:\pi_1(M_\phi) \to \homeo(S^2_\infty)$$

There is another view of $\til{M_\phi}$ which comes from the foliated
structure of $M_\phi$. To describe this point of view, we make use
of some ideas of coarse geometry from Gromov as developed in \cite{mG87}.

The foliation of $\Sigma \times I$ descends to a (taut) foliation of
$M_\phi$ by surfaces which are the fiber of the fibration over $S^1$.
This gives $\til{M_\phi}$ the structure of an open solid cylinder
$$\til{M_\phi} = \til{\Sigma} \times \R$$
The universal cover of each fiber $\Sigma_\theta$ is quasi--isometric
with its pulled back intrinsic metric to the hyperbolic plane $\H^2$,
and can therefore be compactified by its ideal boundary, which is
a topological circle $S^1_\infty$. 

This circle $S^1_\infty$ can just as well be thought of as the Gromov
boundary of the {\em group} $\pi_1(\Sigma)$. The group $\pi_1(M_\phi)$
acts on $\pi_1(\Sigma)$ in the obvious way: the subgroup $\pi_1(\Sigma)$
acts on the left by multiplication, and the element $t$ acts by the
automorphism $\phi_*$. This action on $\pi_1(\Sigma)$ induces an action
of $\pi_1(M_\phi)$ on $S^1_\infty(\pi_1(\Sigma))$, and together with
the action on $\R$ given by the homomorphism to $\Z$, this gives a
(product respecting) action of $\pi_1(M_\phi)$ on $S^1 \times \R$ which
partially compactifies the action on the open cylinder $\til{\Sigma} \times \R$.

The action of $\pi_1(M_\phi)$ on $\R$ is boring; all the information
is already contained in the action on $S^1_\infty$. We denote the
representation inducing this action by
$$\rho_{\text{fol}}:\pi_1(M_\phi) \to \homeo(S^1_\infty)$$

\begin{thm}[Cannon--Thurston \cite{Cannon_Thurston} Continuity 
of Peano map]\label{Cannon_Thurston_map}
Suppose $M_\phi$ is a hyperbolic surface bundle over $S^1$ with fiber
$\Sigma$ and monodromy $\phi$. Then there is a continuous, surjective map
$$P:S^1_\infty \to S^2_\infty$$
which is a semiconjugacy between the two actions of $\pi_1(M_\phi)$. 
That is, for each $\alpha \in \pi_1(M_\phi)$,
$$P \circ \rho_{\text{fol}}(\alpha) = \rho_{\text{geo}}(\alpha) \circ P$$
\end{thm}

Since the image of $S^1_\infty$ under $P$ is closed and
invariant under the action of $\pi_1(M_\phi)$, it is equal to
the entire sphere $S^2_\infty$; that is, it is a {\em Peano curve}, or
{\em sphere--filling map}.

The fact that $P$ is sphere--filling is disconcerting and beautiful, but it is
not the whole story. More interesting is the fact that
$P$ can be {\em approximated by embeddings} in a natural way.

If $\til{\Sigma_\theta}$ denotes the universal cover of a fiber, then
$\til{\Sigma_\theta}$ is a properly embedded plane in $\til{M_\phi} = \H^3$.
By theorem~\ref{Cannon_Thurston_map}, the embedding of $\til{\Sigma_\theta}$
extends continuously to the Peano map on the boundary, by the canonical
identification of $S^1_\infty(\til{\Sigma_\theta})$ with $S^1_\infty(\pi_1(\Sigma))$.
In the unit ball model of $\H^3$, let $p \in \til{\Sigma_\theta}$ be a base point
at the origin. Let $T_i \subset \til{\Sigma_\theta}$ be (a component of)
the intersection of $\til{\Sigma_\theta}$ with 
a family of concentric spheres about $p$. Radial projection from $p$
in $\Sigma_\theta$ identifies each $T_i$ with $S^1_\infty(\til{\Sigma})$, and
radial projection from $p$ in $\H^3$ identifies each $T_i$ with an 
{\em embedded} circle in $S^2_\infty$. We denote the composition of these
identifications by
$$P_i:S^1_\infty \to S^2_\infty$$
which gives a family of maps
which converge in the compact--open topology to $P$. Each $P_i$ decomposes
the complement of its image into two sides, which we can consistently label as
the positive and negative sides, compatibly with an orientation on $S^1_\infty$
and $S^2_\infty$.

Define a {\em positive pair} to be a pair of elements $p,q \in S^1$, and a
choice, for each $P_i$, of an arc $\gamma_i \subset S^2_\infty$ from $P_i(p)$ to
$P_i(q)$ whose interior is disjoint from $P_i(S^1)$ and contained on the
positive side, and which satisfies
$$\lim_{i \to \infty} \text{diameter}(\gamma_i) \to 0$$
We denote a positive pair by $(p,q,\lbrace \gamma_i \rbrace)$.

Now, if $(p_1,p_2,\lbrace \gamma_i \rbrace)$ is one
positive pair and $(q_1,q_2,\lbrace \delta_i \rbrace)$ is another, then either
$\lbrace p_1,p_2\rbrace$ and $\lbrace q_1,q_2\rbrace$ 
are {\em unlinked} as copies of $S^0$ in $S^1_\infty$,
or else all four points are mapped to the same point by $P$. The reason is
that if $\lbrace p_1,p_2\rbrace$ and $\lbrace q_1,q_2 \rbrace$ are linked in
$S^1_\infty$, then $\gamma_i$ and $\delta_i$ lying on the same side of the image
of $P_i$ must {\em intersect}. Since their
lengths converge to $0$ as $i\to \infty$, the claim follows.

The positive pairs define a subset of $S^1 \times S^1$ which generates
a closed equivalence relation, which we denote by $\sim^+$. Similarly, we
can define $\sim^-$ in terms of negative pairs. Note that
distinct equivalence classes of $\sim^+$ say have the property that
they are {\em unlinked} as subsets of $S^1$, in the sense that if
$S^0_1,S^0_2$ are two embedded copies of $S^0$ in $S^1$
which are each contained in distinct equivalence classes of $\sim^+$, 
then the homological linking number of the $S^0_i$'s is $0$. 

Applying this fact to the map $P:S^1_\infty \to S^2_\infty$ lets us 
construct a pair of geodesic laminations $\til{\Lambda}^\pm$ 
of $\til{\Sigma}$ as follows (see section \S~\ref{topology_of_circle} for
a definition and a discussion of geodesic laminations).
The lamination $\til{\Lambda}^+$ is the
union, over equivalence classes $[p]$ of $\sim^+$, of the boundary of the
convex hull of $[p]$, thought of as a subset of $S^1_\infty(\til{\Sigma})$.
It is crucial here that distinct equivalence classes are unlinked, so that
the result is a lamination, and not merely a collection of geodesics.
The action of $\pi_1(\Sigma)$ on $\til{\Sigma}$ preserves these laminations, and
they descend to geodesic laminations $\Lambda^\pm$ on $\Sigma$ which are
preserved by the action of (a homeomorphism isotopic to)
$\phi$. It is not hard to show that these laminations
are transverse, and {\em bind} $\Sigma$, in the sense that complementary
regions are (compact) finite sided polygons. The usual Perron--Frobenius
theory shows that $\Lambda^\pm$ admit transverse measures $\mu^\pm$
which are multiplied by $\lambda,\lambda^{-1}$ respectively by $\phi$, for
some $\lambda > 1$. This is one of the definitions of a
pseudo--Anosov map, and Thurston's theorem on the classification 
of surface homeomorphisms is recovered via a very non--standard route.

Note that {\em a posteriori}, it can be seen that the laminations
$\til{\Lambda}^\pm$ are determined uniquely by the action of $\pi_1(M_\phi)$
on the circle $S^1_\infty$, and can be recovered from the fixed point
data of $\phi_*$ and its conjugates.

Of course this is not a {\em logical} deduction, since the usual
proofs of both theorem~\ref{geometrization_of_bundles} and 
theorem~\ref{Cannon_Thurston_map} depend essentially on 
theorem~\ref{classification_of_homeomorphisms}.

More useful information can be derived from this picture. Each of the invariant
geodesic laminations $\Lambda^\pm$ on $\Sigma$ suspend in the mapping torus to 
two dimensional laminations. Notice that such laminations have some
useful properties. The leaves are covered in $\til{M}_\phi$ by planes. The
finitely many complementary regions are topologically open solid tori, 
which have the extra structure of
finite sided ideal polygon bundles over $S^1$. They are the prototypical
example of {\em very full genuine laminations}, a particularly well behaved
subclass of the class of genuine laminations, introduced by Gabai and Oertel
in \cite{Gabai_Oertel}. Such laminations certify important properties of their
ambient manifold. In \cite{dGwK97},\cite{dGwK98} and \cite{dGwK98b} 
Gabai and Kazez show that an atoroidal $3$--manifold $M$ with a genuine 
lamination has word--hyperbolic fundamental group, has a finite mapping 
class group, and that every self--homeomorphism homotopic to the identity 
is isotopic to the identity.

In another direction, the laminations $\Lambda^\pm$ can be used to produce
a particularly nice flow $X$ transverse to the fibration. The projectively
measured laminations $\til{\Lambda}^\pm$ are dual to a pair of topological
$\R$--trees $T^\pm$. Then the tautological quotient maps of $\til{\Sigma}$ to
$T^+$ and $T^-$ define a map to the product $T^+ \times T^-$ whose image is
topologically a plane. Moreover, the image inherits a pair of singular
foliations $\F^\pm$ by the intersection with factors $T^+ \times \text{point}$
and $\text{point} \times T^-$ of the product structure. This structure is
equivariant, and defines a pair of transversely measured singular foliations
on $\Sigma$ which are transverse to each other and invariant by a suitable
element in the isotopy class of $\phi$.

The suspension flow $X$ of this homeomorphism is {\em pseudo--Anosov}. That is,
away from finitely many orbits, there is a decomposition of the tangent
space $TM$ into a sum $TX \oplus TE^s \oplus TE^u$ which is preserved by the
flow, and where the time $t$ flow multiplies the vectors in the
sub--bundles $TE^s$ and $TE^u$ by factors $O(e^{t\lambda})$ and $O(e^{-t\lambda})$
respectively, for some $\lambda > 0$. Moreover, the singular orbits look like
branched covers of the ordinary orbits, with branch index $n/2$ for some integer
$n\ge 3$. A pseudo--Anosov without such singular orbits is {\em Anosov}.
This pseudo--Anosov flow has the property of being the {\em minimal
entropy flow} transverse to the foliation, and it also has the
property of being {\em quasigeodesic}. That is, flowlines of the lift 
$\til{X}$ in the universal cover are a bounded distance from hyperbolic geodesics
in $\til{M_\phi} = \H^3$.

\subsection{Circle of ideas}

This pencil sketch of the theory of surface diffeomorphisms outlines the
application of this dimensional reduction idea to $3$--manifold theory.
A certain kind of foliation --- namely a fibration --- reduces a $3$--manifold 
$M_\phi$ to a $2$--manifold $\Sigma$ together with some dynamics $\phi$.
Ideal geometry reduces the $2$--manifold $\Sigma$ to a $1$--manifold $S^1_\infty$
together with some further dynamics, namely the action of $\pi_1(\Sigma)$.
The relationship between $S^1_\infty$ and $S^2_\infty$ can be encoded in
another pair of $1$--dimensional objects, namely the laminations $\til{\Lambda}^\pm$,
which are actually encoded in data living only on $S^1_\infty$, and which can
be recovered in principle purely from the dynamics of $\pi_1(M_\phi)$ on 
this $1$--dimensional object.

The goal of this paper is to reproduce as much of this structure as possible
in the context of a more general kind of foliation, namely a {\em taut foliation}.
We follow the principle that smaller is better when it comes to dimension.
Accordingly, we aim to reduce our $3$--manifold, via the use of
some auxiliary dynamical data, to a canonical circle $S^1_\u$ called a {\em universal
circle}, together with a natural representation 
$$\rho_\u:\pi_1(M) \to S^1_\u$$
This circle and representation encodes the original dynamical data, 
or as much of it as is important. In particular, from $S^1_\u$ and $\rho_\u$ we can 
reconstruct the original $3$--manifold $M$ and certify important
topological, geometric and dynamical properties of it.

Loosely speaking, the sources of universal circles are threefold: they arise from
the following three objects, which are all present in the example of surface 
bundle over a circle.
\begin{enumerate}
\item{Taut foliations}
\item{Very full genuine laminations}
\item{Quasigeodesic pseudo--Anosov flows}
\end{enumerate}
Precise definitions of these structures will be deferred until \S~\ref{laminations_theory}.

In the best situation, all three structures give rise to and can be recovered
from the universal circle, and their interactions are encoded in a uniform way. 
For details, consult \cite{dCnD02}. In this paper we aim to show how, 
under suitable circumstances, one of
the structures --- a taut foliation --- gives rise to another: a (pair of) very
full genuine laminations.

\subsection{Atoroidal versus algebraically atoroidal}

Throughout this paper, we use the term {\em atoroidal} in a slightly nonstandard
way as shorthand for
{\em algebraically atoroidal}. A $3$--manifold $M$ is algebraically
atoroidal if there is no $\Z \oplus \Z$ in $\pi_1(M)$. A $3$--manifold $M$
is {\em geometrically atoroidal} if every essential embedded torus is boundary parallel.
For closed $3$--manifolds, the two terms are interchangeable except when $M$ is
a small Seifert fibered space; i.e. a Seifert fibered space over a triangle orbifold.

In statements of important theorems, in order to minimize confusion, we try to use
the longer term {\em algebraically atoroidal}.

\subsection{Statement of results}

In this section we state our results precisely. 

\S~\ref{topology_of_circle} gathers basic results and constructions in
the point set topology of $S^1$ which are used again and again throughout
the rest of the paper. We define laminations of $S^1$, laminar relations
on $S^1$, and geodesic laminations of the hyperbolic plane $\H^2$, and we
show how to move back and forward between these three kinds of objects.
We also define {\em monotone maps} between circles, which are degree one
maps whose point preimages are connected. The most important theorem is
theorem~\ref{unlinked_extension}, which concerns
continuous families of monotone maps. This is a technical theorem which is
used in later sections, especially \S~\ref{constructing_laminations}.

\S~\ref{laminations_theory} is mostly expository, being a brief introduction
to the theory of taut foliations and their cousins, essential and genuine
laminations in $3$--manifolds. It is hardly a comprehensive survey, but
it gives the definitions of the most important objects and constructions, and
gives statements of and references to all the basic foundational results
that we make use of in this paper. The subsection \S~\ref{Candel_uniformization}
describes Candel's uniformization theorem~\ref{Candel_uniformize} 
for hyperbolic laminations, and describes how to use this theorem to construct
the circle bundle $E_\infty$ over the leaf space $L$ of $\til{\F}$ for
a taut foliation $\F$. The fiber of $E_\infty$ over a leaf $\lambda$ of
$\til{\F}$ is just the ideal boundary of $\lambda$ in the sense of Gromov
(see \cite{mG87}). The more precise details of Candel's theorem are necessary
to define the correct topology on $E_\infty$. The circle bundle $E_\infty$
is used repeatedly throughout the rest of the paper.

\S~\ref{constructing_circle_section} is also expository. We present the
outlines of proofs of the Leaf Pocket theorem and 
the Universal Circle theorem from \cite{dCnD02} (theorems 5.2 and 6.2 respectively
in \cite{dCnD02}). 
For most of \S~\ref{constructing_laminations}, we do not need the
details of the proofs of these theorems, and we proceed as far as possible
from the axiomatic statements of these theorems. However, later in the paper
we need to make use of some of the properties of the
universal circles constructed in \cite{dCnD02}, and therefore it is necessary
to explain the constructions in some detail.

\S~\ref{constructing_laminations} contains the really new results in
this paper. We construct a pair of laminations $\Lambda^\pm_\u$ of the
universal circle $S^1_\u$ constructed in \S~\ref{constructing_circle_section}, and
use these laminations to construct a pair of $2$--dimensional 
laminations $\Lambda^\pm_\split$ of $M$ which are transverse to $\F$. We then
go on to establish basic properties of these laminations. We say a foliation $\F$
has {\em $2$--sided branching} if the leaf space $L$ of the pullback foliation $\til{\F}$
on the universal cover branches in both the positive and the negative directions.
Our main result is the following:

\begin{lam_exist_thm}
Let $\F$ be a co--orientable taut foliation of a closed, orientable
algebraically atoroidal $3$--manifold $M$. Then either $\F$ has $2$--sided branching
and is the weak stable foliation of an Anosov flow, 
or else there are a pair of very full
genuine laminations $\Lambda^\pm_\split$ transverse to $\F$.
\end{lam_exist_thm}

It follows by work of Gabai and Kazez, that a closed $3$--manifold with
a taut foliation either contains a $\Z \oplus \Z$ in its fundamental group,
or contains an Anosov flow whose stable and unstable foliation have $2$--sided branching,
or else it has word--hyperbolic fundamental group, the mapping class group is
finite, and every self--homeomorphism homotopic to the identity is isotopic
to the identity.

Finally, in \S~\ref{lamination_dynamics} we discuss the dynamics of the
laminations $\Lambda^\pm_\split$. This section is mainly descriptive, and serves
to illustrate some of the structure developed in earlier sections.

\subsection{Notation}

We make use of certain conventions for notation throughout this
paper, and try to be consistent throughout. For an object or structure
$X$ in a $3$--manifold $M$, $\til{X}$ will denote the pull back
of $X$ to the universal cover $\til{M}$, where this makes sense. Surfaces and
manifolds will be denoted by upper case Roman letters $S,M,N$ etc.
and points by lower case Roman letters $p,q,r$ etc. Foliations
will be denoted by script letters $\F,\G$ etc. and
laminations by upper case Greek letters $\Lambda$ etc. Leaves will
be denoted by lower case Greek letters $\lambda,\mu,\nu$ etc. Guts of
genuine laminations will be denoted by Gothic $\fG$ and core circles
of interstitial annuli by lower case Gothic letters $\fc$.

\subsection{Acknowledgements}
I am indebted to Andrew Casson and S\'ergio Fenley for conversations and
discussions concerning this material over a number of years. More recently,
I am greatly indebted to Nathan Dunfield for careful comments about the
content and readability of this manuscript. I would also like to thank the
anonymous referee for his very careful reading of this manuscript.
Finally, it should be stressed
that I learned much of the overarching philosophy and organizing principles
of the theory of foliations and geometrization which permeate this paper and
its siblings \cite{dC99b}, \cite{dC00} and \cite{dCnD02} from very long,
extensive conversations with Bill Thurston, mainly during 1996--1999. Many of
the theorems and constructions in this paper are at least inspired by
ideas and comments from Bill.

While writing this paper, I was partially supported by a Sloan Research Fellowship,
and NSF grant DMS 0405491.

\section{The topology of $S^1$}\label{topology_of_circle}

In this section we establish basic properties of the point set topology of
$S^1$ which will be used in the rest of the paper. Good general references for
point set topology in low dimensions are
\cite{Hocking_Young},\cite{Bing} and \cite{Moise}.

\subsection{Laminations of $S^1$}

\begin{defn}
We let $S^0$ denote the $0$ sphere; i.e. the discrete, two element set.
Two disjoint copies of $S^0$ in $S^1$ are {\em homologically linked},
or just {\em linked} if the points in one of the $S^0$'s are contained
in different components of the complement of the other. Otherwise
we say they are {\em unlinked}.
\end{defn}

Note that the definition of linking is symmetric.

\begin{defn}
A {\em lamination} $\Lambda$ 
of $S^1$ is a closed subset of the space of unordered
pairs of distinct points in $S^1$ with the property that no two elements
of the lamination are linked as $S^0$'s in $S^1$. The elements of
$\Lambda$ are called the {\em leaves} of the lamination.
\end{defn}

The space of unordered pairs of distinct points in $S^1$ may be thought
of as a quotient of $S^1 \times S^1\backslash \text{diagonal}$ by
the $\Z/2\Z$ action which interchanges the two components. Topologically,
this space is homeomorphic to a M\"obius band.

Most readers will be familiar with the concept of a {\em geodesic lamination}
on a hyperbolic surface.

\begin{defn}
A {\em geodesic lamination} $\Lambda$ on a complete hyperbolic surface
$\Sigma$ is a closed union of disjoint embedded complete geodesics.
\end{defn}

For a thorough development of the elementary theory of geodesic laminations,
see \cite{Casson_Bleiler}. A geodesic lamination of $\Sigma$ pulls back
to define a geodesic lamination of $\H^2$. Geodesic laminations of $\H^2$
and laminations of $S^1$ are essentially equivalent objects, as the following
construction shows:

\begin{construct}\label{circle_geodesic}
Let $\Lambda$ be a lamination of $S^1$. We think of $S^1$ as the boundary
of $\H^2$ in the unit disk model. Then we construct a geodesic lamination of
$\H^2$ whose leaves are just the geodesics whose endpoints are leaves of
$\Lambda$. We will sometimes denote this geodesic lamination by $\Lambda_\geo$.
Conversely, given a geodesic lamination $\Lambda$ of $\H^2$, we
get a lamination of the ideal boundary $S^1_\infty$ whose leaves are
just the pairs of endpoints of the leaves of $\Lambda$. 
\end{construct}

There is another perspective on circle laminations, coming from
equivalence relations. The correct class of equivalence relations for
our purposes are {\em upper semi--continuous decompositions}.

\begin{defn}\label{monotone_definition}
A {\em decomposition} of a topological 
space $X$ is a partition into compact subsets.
A decomposition $\G$ 
is {\em upper semi--continuous} if for every decomposition element
$\zeta \in \G$ and every open set $U$ with $\zeta \subset U$, there
exists an open set $V \subset U$ with $\zeta \subset V$ such that every
$\zeta' \in \G$ with $\zeta' \cap V \neq \emptyset$ has
$\zeta' \subset U$. The decomposition is {\em monotone} if its elements are connected.
\end{defn}

A proper map from a Hausdorff space $X$ to a Hausdorff space $Y$
induces a decomposition of $X$ by its point preimages which is
upper semi--continuous. Conversely, the quotient of a Hausdorff space by
an upper semi--continuous decomposition is Hausdorff, and the tautological
map to the quotient space is continuous and proper. See e.g. \cite{Hocking_Young}.

\begin{defn}
An equivalence relation $\sim$ on $S^1$ is {\em laminar} if the equivalence
classes are {\em closed}, if the resulting decomposition
is upper semicontinuous, and if distinct
equivalence classes are unlinked as subsets of $S^1$. That is, if
$S^0_1,S^0_2 \subset S^1$ are two $S^0$'s which are contained in distinct
equivalence classes, then they are not homologically linked in $S^1$.
\end{defn}

We now show how to move back and forth between circle laminations and
laminar relations.

\begin{construct}\label{equivalence_circle}
Given a laminar equivalence relation $\sim$ of $S^1$, we think of $S^1$ as
the ideal boundary of $\H^2$. Then for every equivalence class $[p]$
of $\sim$ we form the convex hull 
$$H([p]) \subset \H^2$$ 
and the boundary of the
convex hull 
$$\Lambda([p]) = \partial H([p]) \subset \H^2$$
We let $\Lambda$ denote the union over all equivalence classes $[p]$:
$$\Lambda = \bigcup_{[p]} \Lambda([p])$$
Then the fact that the equivalence classes are unlinked implies that the
geodesics making up $\Lambda$ are disjoint. Moreover, the fact that
$\sim$ is upper semicontinuous implies that $\Lambda$ is {\em closed}
as a subset of $\H^2$. That is, it is a geodesic lamination, and
determines a lamination of $S^1$ by construction~\ref{circle_geodesic}.

Conversely, given a lamination $\Lambda$ of $S^1$, we may form the quotient
$Q$ of $S^1$ by the smallest equivalence relation which collapses every leaf to
a point. This is not necessarily Hausdorff; we let $Q'$ denote the
Hausdorffification. Then the map from $S^1$ to $Q'$ induces an
upper semi--continuous decomposition of $S^1$. Moreover, this equivalence
relation is obviously unlinked; in particular, it is laminar.
\end{construct}

We abstract part of construction~\ref{equivalence_circle} to show that
every subset $K \subset S^1$ gives rise to a lamination, as follows

\begin{construct}\label{boundary_hull}
Let $K \subset S^1$ be arbitrary. Think of $S^1$ as $\partial \H^2$, and let
$H(\overline{K}) \subset \H^2$ be the convex hull of the {\em closure} of $K$ in
$S^1$. Then the boundary $\partial H(\overline{K})$ is a geodesic lamination of
$\H^2$, which determines a lamination of $S^1$ by construction~\ref{circle_geodesic}.
We denote this lamination of $S^1$ by $\Lambda(K)$.
\end{construct}

\subsection{Monotone maps}

\begin{defn}
Let $S^1_X,S^1_Y$ be homeomorphic to $S^1$.
A continuous map $\phi:S^1_X \to S^1_Y$ is {\em monotone} if it is
degree one, and if it induces a monotone decomposition of $S^1_X$, in the
sense of definition~\ref{monotone_definition}.
\end{defn}

Note that the target and image circle should not be thought of as the same
circle.

Equivalently, a map between circles
is monotone if the point preimages are connected
and contractible. Said yet another way, a map is monotone if it does
not reverse the cyclic order on triples of points for some choice of
orientations on the target and image circle.

\begin{defn}
Let $\phi:S^1_X \to S^1_Y$ be monotone. The {\em gaps} of $\phi$ are the
maximal open connected intervals in $S^1_X$ in the preimage of single points
of $S^1_Y$. The {\em core} of $\phi$ is the complement of the union of
the gaps.
\end{defn}

Note that the gaps of $\phi$ are the connected
components of the set where $\phi$ is locally constant. 

Recall that a set is {\em perfect} if no element is isolated.

\begin{lem}\label{perfect_core}
Let $\phi:S^1_X \to S^1_Y$ be monotone. Then the core of $\phi$ is
perfect.
\end{lem}
\begin{proof}
The core of $\phi$ is closed. If it is not perfect, there is some
point $p \in \core(\phi)$ which is isolated in $\core(\phi)$. Let
$p^\pm$ be the nearest points in $\core(\phi)$ to $p$ on either side,
so that the open oriented intervals $p^-p$ and $pp^+$ are gaps of $\phi$. But
then by definition, 
$$\phi(p^-) = \phi(p) = \phi(p^+) = \phi(r)$$
for any $r$ in the oriented interval $p^-p^+$. So by definition, the interior
of this interval is contained in a single gap of $\phi$. In particular, $p$
is contained in a gap of $\phi$, contrary to hypothesis.
\end{proof}

It follows that the set of points in $\core(\phi)$ which are nontrivial
limits from both directions is dense in $\core(\phi)$.

\begin{figure}[ht]
\centerline{\relabelbox\small \epsfxsize 2.0truein
\epsfbox{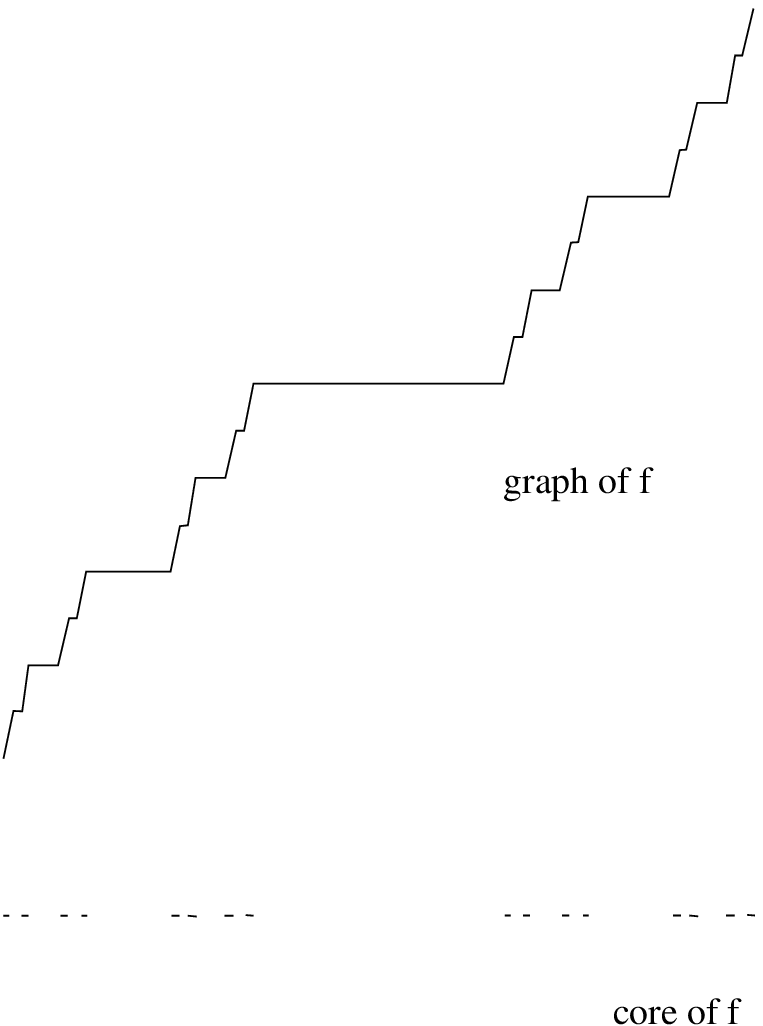}
\endrelabelbox}
\caption{The graph of $f$ is called the {\em Devil's staircase}.}
\label{devil_stair}
\end{figure}

\begin{exa}[The Devil's staircase]
Let $f:[0,1] \to [0,1]$ be the function defined as follows. If $t \in [0,1]$,
let
$$0 \cdot t_1t_2t_3 \cdots$$
denote the base $3$ expansion of $t$. Let $i$ be the smallest index for which $t_i = 1$.
Then $f(t) = s$ is the number whose base $2$ expansion is
$$0 \cdot s_1s_2s_3 \cdots s_i 0 0 \cdots$$
where each $s_j = 1$ iff $t_j = 1 \text{ or } 2$ and $j \le i$, and $s_j = 0$
otherwise. The graph of this function is illustrated in figure~\ref{devil_stair}.

The core of this map is the usual middle third Cantor set.
\end{exa}

\begin{defn}
Let $B$ be a topological space, and $E$ a circle bundle over $B$. A {\em monotone
family} of maps is a continuous map $$\phi:S^1 \times B \to E$$
which covers the identity map on $B$, and which restricts for each
$b \in B$ to a monotone map of circles
$$\phi_b = \phi|_{S^1 \times b}:S^1 \times b \to E_b$$
We denote a monotone family by the triple $(E,B,\phi)$.
\end{defn}

\begin{lem}\label{cores_semicontinuous}
Let $(E,B,\phi)$ be a monotone family. Then the 
family of subsets $\lbrace \core(\phi_b)\rbrace$
vary lower semicontinuously as a function of $b \in B$, in the Hausdorff
topology. That is, if $x \in \core(\phi_b)$, then if $b_i \to b$ in $B$,
there are points $x_i \in \core(\phi_{b_i})$ such that $x_i \to x$.
\end{lem}
\begin{proof}
Let $x \in \core(\phi_b)$. By lemma~\ref{perfect_core} it follows there is some
sequence of distinct points $x_i \to x$ such that $\phi_b(x_i) \ne \phi_b(x_j)$
for each $i,j$. It follows that for each $i$ there is a $k$ such that
$\phi_{b_K}(x_i) \ne \phi_{b_K}(x_{i+1})$ for all $K \ge k$. In particular, the
core of $\phi_{b_K}$ contains some point between $x_i$ and $x_{i+1}$. The lemma
follows.
\end{proof}

It follows that the closure of the union of gaps of $\phi_b$ varies
upper semicontinuously as a function of $b$. An alternate proof of 
lemma~\ref{cores_semicontinuous} uses the fact that the closures of gaps
are exactly the nontrivial elements in the decomposition of $S^1 \times B$ induced by
$\phi$.

\begin{defn}\label{define_core_of_union}
Let $(E,B,\phi)$ be a monotone family. Let $X \subset B$ be a subspace. Define
$$\core(X) = \overline{\bigcup_{b \in X} \core(\phi_b)}$$
\end{defn}

Notice that we define $\core(X)$ to be the {\em closure} of the union of
the cores of $\phi_b$ over all $b \in X$, and not simply the ordinary
union. This is important to keep in mind in the sequel; we will refer
to this construction on a number of occasions in section~\ref{constructing_laminations}.

\begin{thm}\label{unlinked_extension}
Let $(E,B,\phi)$ be a monotone family, and suppose $X,Y$ are {\em path connected}
subsets of $B$.
Suppose for each $x \in X$ and $y \in Y$ that $\core(\phi_x)$ and $\core(\phi_y)$
are unlinked. Then $\core(X)$ and $\core(Y)$ are unlinked.
\end{thm}
\begin{proof}
Since $\core(\phi_x)$ and $\core(\phi_y)$ are
unlinked for each pair $x \in X,y \in Y$,
it follows that $\core(\phi_x)$ is contained in the closure of
a single gap of $\core(\phi_y)$, and {\it vice versa}.

We claim that for every $x\in X$, $\core(\phi_x)$ is 
contained in the closure of the {\em same} gap of $\core(\phi_y)$. 
For, let $g$ be a gap of $\phi_y$, and let $T_g \subset X$ be
the set of points $t$ for which $\core(\phi_t) \subset \overline{g}$.
Since $\overline{g}$ is closed, by lemma~\ref{cores_semicontinuous} the set
$T_g$ is closed. Moreover, by lemma~\ref{perfect_core}, distinct gaps have
disjoint closures, and therefore if $g_1,g_2$ are distinct gaps,
$T_{g_1}$ and $T_{g_2}$ are disjoint.
Let $x_1,x_2 \in X$ be arbitrary, and let $\gamma$ be a path in $X$ from $x_1$ to
$x_2$. Then $\gamma$ is decomposed into closed subsets which are
the intersections $\gamma \cap T_g$ as $g$ varies over the gaps of $\phi_y$.
But there are only countably many gaps of $\phi_y$. On the other hand, any
decomposition of an interval into countably many closed subsets has only one
element, by a beautiful theorem of Sierpinski \cite{Sierpinski}. 
It follows that $T_g = X$, and every $\core(\phi_x)$
is contained in the same gap $g$ of $\phi_y$. We can therefore label $g$
unambiguously as $g_y$, and similarly construct $g_{y'}$ for every other $y' \in Y$.

Now, as $y$ varies in $Y$, the closures of gaps $\overline{g_y}$ do not
vary continuously, but merely upper semicontinuously. In particular, if $y_i \to y$
then
$$\lim_{i \to \infty} \overline{g_{y_i}} \subset \overline{g_y}$$
for every Hausdorff limit. Since each $\overline{g_{y_i}}$ is a closed arc, the
same is true of each Hausdorff limit. For each such discontinuous limit, i.e.
where $\lim_{i \to \infty} \overline{g_{y_i}} \neq \overline{g_y}$, we interpolate
a $1$--parameter family of closed arcs from $\lim_{i \to \infty} \overline{g_{y_i}}$
to $\overline{g_y}$ which are all contained in $\overline{g_y}$. Let $\G$ denote
the union of the set of
arcs $\overline{g_y}$ with $y \in Y$ and the arcs in the interpolating families. Then
$\G$ is a {\em connected} subset of the space of closed arcs in $S^1$. It follows
that the intersection
$$\bigcap_{\gamma \in \G} \gamma = \bigcap_{y \in Y} \overline{g_y}$$
is a {\em connected} arc, which contains $\core(X)$, and whose 
interior is in the complement
of $\core(Y)$. It follows that $\core(X)$ and $\core(Y)$ are unlinked, as claimed.
\end{proof}

\subsection{Pushforward of laminations}

Laminations of $S^1$ can be pushed forward by monotone maps.

\begin{defn}
Let $\Lambda$ be a lamination of $S^1_X$, and $\phi:S^1_X \to S^1_Y$ a
monotone map. Then $\phi$ induces a map from unordered pairs of points in $S^1_X$
to unordered pairs of points in $S^1_Y$. We let $\phi(\Lambda)$ denote the
image of $\Lambda$ in the complement of the diagonal.
\end{defn}

\begin{lem}\label{pushforward_is_lamination}
Let $\phi:S^1_X \to S^1_Y$ be monotone, and let $\Lambda$ be a lamination of $S^1_X$.
Then $\phi(\Lambda)$ is a lamination of $S^1_Y$.
\end{lem}
\begin{proof}
The map $\phi$ induces a continuous map from $S^1_X \times S^1_X \to S^1_Y \times S^1_Y$
which takes the diagonal to the diagonal. It follows that the image of $\Lambda$
is closed in $S^1_Y \times S^1_Y \backslash \text{diagonal}$. It remains to show
that it is unlinked. But monotone maps do not reverse the cyclic order of subsets;
the claim follows.
\end{proof}

Laminations can also be pulled back by monotone maps.

\begin{defn}
Let $\Lambda$ be a lamination of $S^1_Y$, and $\phi:S^1_X \to S^1_Y$ a monotone
map. Then $\Lambda$ determines a laminar equivalence relation $\sim_Y$ on $S^1_Y$, by 
construction~\ref{equivalence_circle}. Let $\sim_X$ be the equivalence relation
on $S^1_X$ whose equivalence classes are the preimages of equivalence classes in
$\sim_Y$. Then $\sim_X$ is a laminar relation, and induces a lamination of
$S^1_X$ by construction~\ref{equivalence_circle} which we denote $\phi^{-1}(\Lambda)$.
\end{defn}

The proof that $\sim_X$ is laminar follows immediately from the fact that
$\phi$ is monotone.

\subsection{Coarse geometry of the hyperbolic plane}

This subsection summarizes some basic facts about coarse geometry,
quasigeodesic and geodesics in the hyperbolic plane $\H^2$. We will use
the results in this subsection implicitly throughout the rest of the paper,
usually without comment. It is included here simply as a service to those
readers who might be unfamiliar with or hazy on this material.
 
The material in this subsection is completely standard; excellent
references are \cite{mG87} and \cite{Casson_Bleiler}.

\begin{defn}
Let $X$ be a metric space.
Let $k>1$ and $\epsilon >0$. A {\em $(k,\epsilon)$--quasigeodesic} is
a map $l:\R \to X$ such that for all $p,q \in \R$,
$$\frac 1 k d_X(l(p),l(q)) - \epsilon \le  d(p,q) \le k d_X(l(p),l(q)) + \epsilon$$
where $d_X(\cdot,\cdot)$ denotes distance in $X$, and $d(\cdot,\cdot)$ denotes
Euclidean distance in $\R$.
\end{defn}

\begin{lem}\label{quasi_near_geodesic}
Let $\gamma \subset \H^2$ be a $(k,\epsilon)$--quasigeodesic. Then there is a
constant $C(k,\epsilon)$ such that there is a complete geodesic $\gamma_s \subset \H^2$
which is $C(k,\epsilon)$ close to $\gamma$ in the Hausdorff metric.
\end{lem}

\begin{lem}\label{half_distance}
Let $\gamma_1,\gamma_2$ be two geodesics in $\H^2$ which are Hausdorff distance
$C$ apart on subsegments $\gamma_1',\gamma_2'$ of length $t$. Then there is a constant
$C_1$ which does not depend on $C$,
such that $\gamma_1'$ and $\gamma_2'$ are Hausdorff distance $C/2$ apart
on subsegments of length $t - C_1$.
\end{lem}

By applying lemma~\ref{half_distance} iteratively, one sees that for any $C$ and
any $\epsilon$ there is a $t(C)$ such that two geodesic segments of length $\ge t(C)$ 
which are Hausdorff distance $C$ apart are Hausdorff distance $\epsilon$ apart 
on their middle third subsegments. In particular,
two bi--infinite geodesics which are a finite Hausdorff distance apart are equal.

\begin{lem}\label{continuous_extension}
Let $\Lambda$ be a lamination of $S^1_X$, and let $h:S^1_X \to S^1_Y$ be a
homeomorphism. Suppose $S^1_X$ and $S^1_Y$ bound copies $\H^2_X,\H^2_Y$ of
the hyperbolic plane, and let $\Lambda_\geo,(h(\Lambda))_\geo$ be the geodesic
laminations determined by construction~\ref{circle_geodesic}. 
Then $h$ extends to a homeomorphism 
$$H:\overline{\H^2_X} \to \overline{\H^2_Y}$$
taking $\Lambda_\geo$ to $(h(\Lambda))_\geo$.
\end{lem}

See \cite{Casson_Bleiler} for the proofs of these facts.

\begin{note}
If $X,Y$ are Hausdorff distance $\le C$ apart in some third metric space $Z$, 
we sometimes abbreviate this by saying $X$ and $Y$ are {\em $C$ close in $Z$}, or
just {\em $C$ close} if $Z$ is understood.
\end{note}

\section{The theory of essential laminations}\label{laminations_theory}

In this section we define taut foliations and essential laminations, 
and present some of their fundamental theory and properties. 
None of the material in
this section is new, but it is presented here for the convenience of
the reader. A good reference for the general theory of foliations 
is \cite{CCI}. References for basic $3$--manifold topology are
\cite{Hempel}, \cite{Hatcher} and \cite{Jaco}.

\subsection{Taut foliations}

A $2$--dimensional foliation of a $3$--manifold is a partition into surfaces
called {\em leaves} whose local connected components have a product structure.

A foliation is {\em orientable} if the leaves can be oriented in a continuously
varying way. It is {\em co--orientable} if the transverse space can be
oriented. By passing to a cover of index at most $4$, we can assume that
our foliations are orientable and co--orientable. Note that this implies that
the ambient manifold is itself orientable. Throughout this paper we assume
that all our manifolds are orientable, and all our foliations are orientable
and co--orientable.

A basic atomic unit in the theory of foliations of $3$--manifolds is
the following:

\begin{defn}\label{Reeb_component}
Let $H$ be the closed upper half space in $\R^3$ minus the origin. $H$ is
foliated by its intersection with horizontal planes. Every leaf is either
a plane or a punctured plane. A nontrivial dilation $\phi$ centered at the origin
preserves this foliation of $H$, so it descends to a foliation of the
quotient manifold $S = H/\langle \phi \rangle$ which is a solid torus. This
foliation is called a {\em Reeb component}.
\end{defn}

\begin{figure}[ht]
\centerline{\relabelbox\small \epsfxsize 2.0truein
\epsfbox{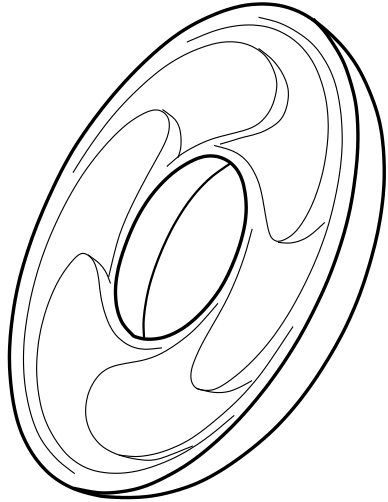}
\endrelabelbox}
\caption{A Reeb component is a foliation of a solid torus by planar
leaves. Each leaf is like a sock which is swallowed by the next sock.}
\label{reeb}
\end{figure}

See figure~\ref{reeb} for a cutaway of half a Reeb component.

Every closed $3$--manifold admits a $2$--dimensional foliation 
\cite{Lick},\cite{Thur_exist} but such foliations typically contain
Reeb components. Notice that the boundary torus of a Reeb component
is compressible, by a compression lying entirely within the component.
Conversely, foliations with topologically essential leaves are much
harder to construct.

\begin{defn}
A $2$--dimensional foliation $\F$ of a closed $3$--manifold $M$ is {\em taut} if 
no leaf is a sphere or projective plane, and there
is a map $\phi:S^1 \to M$ which is transverse to $\F$, and which intersects
every leaf.
\end{defn}

The condition that $\F$ has no sphere or projective plane leaf is a convention
to rule out some very special cases. By the Reeb stability theorem \cite{Reeb}
if $\F$ contains a sphere or projective plane leaf, then $\F$ is finitely
covered by a product foliation of $S^2 \times S^1$ by leaves $S^2 \times \text{point}$.

The relationship between taut foliations and Reeb components is complementary,
so that one has the following theorem, which is basically due to Novikov
\cite{sN65}:

\begin{thm}[Novikov]\label{Novikov_theorem}
A foliation $\F$ of an atoroidal $3$--manifold $M$ is taut iff it contains
no Reeb components.
\end{thm}

Taut foliations have several distinct lives: a topological life, a dynamical life,
and a geometric life. 

Firstly, they certify many useful topological 
properties of $M$. For instance, there is the following theorem 
which combines work of  Palmeira \cite{cP78} with earlier important work 
of Novikov \cite{sN65} and Rosenberg \cite{Rosenberg}. 

\begin{thm}[Palmeira, Novikov, Rosenberg]\label{Palmeira_theorem}
Let $M$ be a $3$--manifold which admits a taut foliation $\F$. Then
the universal cover $\til{M}$ is homeomorphic to $\R^3$, and the
leaves of $\til{\F}$ are planes. Moreover, there is a foliation $F$ of $\R^2$ by
lines such that the pair $(\til{M},\til{\F})$ is
topologically equivalent to a product
$$(\til{M},\til{\F}) = (\R^2,F)\times \R$$
\end{thm}

The {\em leaf space} of $\til{\F}$ is just the quotient space of $\til{M}$
by the equivalence relation whose equivalence classes are the leaves of
$\til{\F}$. We denote this leaf space by $L$.
It follows from theorem~\ref{Palmeira_theorem} that the leaf space of
$\til{\F}$ is simply connected; on the other hand, it is typically non--Hausdorff. 
It will become apparent, especially in \S~\ref{constructing_laminations} that
this non--Hausdorffness of $L$ is fundamentally the source of much of the
structure that we develop, so one should not be too upset to encounter it.
In any case, the action of $\pi_1(M)$ on $\til{M}$ descends to an action 
on $L$ by homeomorphisms. 

\begin{defn}\label{holonomy_homomorphism}
We denote the representation inducing the action of $\pi_1(M)$ on $L$ by
$$\rho_{\text{hol}}:\pi_1(M) \to \homeo(L)$$
and call this homomorphism the {\em holonomy homomorphism}.
\end{defn}

Since many readers might not have considered the subject, we 
should say a few words at this point about non--Hausdorff $1$--manifolds.
The kinds of non--Hausdorff $1$--manifolds we consider are obtained
from countably many copies of the open unit interval $I_i$ by identifying
open subsets in pairs. In this way, each $I_i$ embeds into the quotient
space, and each point in the quotient space is contained in the
interior of an embedded interval. One may pass from $L$ to its
{\em Hausdorffification}, which is just the maximal Hausdorff quotient,
and is obtained from $L$ by inductively identifying pairs of points which are not contained
in disjoint open subsets; i.e. it is obtained by quotienting out the {\em closed}
equivalence relation generated by the property of being nonseparated.
The Hausdorffification of $L$ is also simply--connected,
and is a more familar kind of object. It is homeomorphic to 
the underlying topological space of an $\R$--tree (although in fact,
the underlying space of an $\R$--tree is more general, and is not typically
assumed to be 2nd countable). In practice, one reasons about $L$ by thinking
about this Hausdorff quotient, and remembering that branch points of the
quotient correspond to unseparable sets of points in $L$. The study of
the action of $\pi_1(M)$ on $L$ by the holonomy homomorphism falls into
the domain of arboreal group theory. See \cite{RSS} for an important
example of the kinds of results that may be obtained by such methods.

A co--orientation on $\F$ pulls back to a co--orientation on $\til{\F}$, and
defines an orientation on $L$. So the Hausdorff quotient should be thought of
more as having the structure of an {\em oriented train track} (i.e. there is a combing
at the branch points into positive and negative directions) than a tree.

\begin{figure}[ht]
\centerline{\relabelbox\small \epsfxsize 2.0truein
\epsfbox{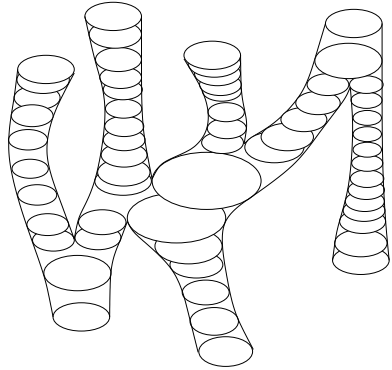}
\endrelabelbox}
\caption{This foliation of a topological ball by planes
exhibits two--sided branching.}
\label{branching}
\end{figure}

Taut foliations can be classified in terms of the kind of branching 
exhibited by the (Hausdorffified) leaf space $L$. If $L$ is Hausdorff, 
equivalently if its Hausdorffification does not branch at all, 
then of course it is homeomorphic to $\R$. In 
this case $\F$ is said to be {\em $\R$--covered}. If $L$ branches in
only one direction (e.g. the positive direction), 
we say $\F$ has {\em one sided branching}, and otherwise
we say $\F$ has {\em two sided branching}. If $\F$ has one sided branching,
then necessarily $\F$ is co--orientable. In this paper we will concentrate
on the case that $\F$ has two sided branching. Analogous results in the
case that $\F$ is $\R$--covered or has one sided branching are contained in
\cite{dC99b} and \cite{dC00} respectively. See figure~\ref{branching} for an
example of (part of) the universal cover of a foliation with two--sided branching.

The orientation on $L$ defines a {\em partial order} on $L$, as follows.

\begin{defn}
The canonical partial order on $L$ is defined as follows.
Let $\lambda,\mu$ be leaves of $\til{\F}$
If there is a positively oriented transversal to $\til{\F}$ from $\lambda$ to
$\mu$, then $$\mu>\lambda$$
Similarly, if there is a negatively oriented transversal, then $\lambda > \mu$.
Note that if $\lambda < \mu$ and $\mu < \lambda$ then $\mu = \lambda$, by
theorem~\ref{Palmeira_theorem}. If there is no transversal between $\mu$ and
$\lambda$, we say the leaves are {\em incomparable}.
\end{defn}

Note that the co--orientation of $\til{\F}$ lets us define unambiguously the
{\em positive} and {\em negative sides} of $\lambda$ in $M$. Every leaf
of $\til{\F} \backslash \lambda$ is either on the positive or negative side.
Moreover, if $\mu > \lambda$, then $\mu$ is on the positive side, and if
$\mu < \lambda$, then $\mu$ is on the negative side, but {\em not} conversely.
The reader should be careful to distinguish between the two notions.

Taut foliations equally well certify useful topological properties of
{\em surfaces}. The following theorem amalgamates a theorem of
Thurston \cite{Thurston_norm} and a (much harder) theorem of Gabai \cite{Gabai_fol}.
Here $\chi(\Sigma)$ denotes the Euler characteristic 

\begin{thm}[Gabai, Thurston]\label{finite_depth_existence_theorem}
Let $M$ be a compact connected irreducible orientable $3$-manifold whose boundary
$\partial M$ is a (possibly empty) union of tori. A properly embedded homologically
essential surface $\Sigma$ is a leaf of a taut foliation of $M$ if and only if it
minimizes $-\chi(\Sigma)$ amongst all proper embedded surfaces with no spherical
components in its homology class.
\end{thm}

\begin{rmk}
Given a surface $S$ in a manifold $M$ satisfying the hypotheses of 
theorem~\ref{finite_depth_existence_theorem}, the {\em Thurston norm} of $S$
is defined to be the sum of $-\chi(S_i)$ over all non-spherical components $S_i$
of $S$. So this theorem may be re-stated as saying that a properly embedded
surface without spherical components in a manifold $M$ (as above) is
a leaf of a taut foliation iff it is Thurston norm minimizing in its homology class.
The ``if" direction is part of Theorem 5.5, page 445 from \cite{Gabai_fol}; the ``only if"
is part of Corollary 1, page 118 from \cite{Thurston_norm}.
\end{rmk}

The second life of a taut foliation is dynamical. There is a basic
duality between minimal surfaces and volume preserving transverse
flows, which in its most rudimentary form is just the min cut --- max flow
theorem from graph theory. 
The next theorem of Sullivan \cite{Sullivan_cycles}, \cite{Sullivan}
makes this precise:

\begin{thm}[Sullivan]\label{minimal_surface_theorem}
Let $\F$ be a co--orientable $C^2$ foliation of $M$. The following are equivalent:
\begin{enumerate}
\item{$\F$ is taut.}
\item{$\F$ admits a volume preserving transverse flow for some volume form.}
\item{There is a closed $2$--form $\theta$ on $M$ which is positive on $T\F$.}
\item{There is a Riemannian metric on $M$ for which leaves of $\F$ are
(calibrated) minimal surfaces.}
\end{enumerate}
\end{thm}

If $\F$ is not $C^2$, there is a combinatorial version of this theorem due
to Hass, based on an idea of Thurston \cite{Hass}. We will not make use of
the theorem of Sullivan in this paper except as a suggestive analogy.

Finally, foliations reveal geometry of the underlying $3$--manifold. One of
the main applications of this paper will be to show that the existence of
a taut foliation on a $3$--manifold allows one to construct auxiliary objects
which certify geometric properties of $M$, for instance word hyperbolicity
of $\pi_1(M)$. This will be developed in detail in the sequel, and therefore
we postpone a discussion until the appropriate time.

\subsection{Essential laminations}

A $2$--dimensional {\em lamination} of a $3$--manifold is a decomposition
of a closed subset into surfaces, which come together locally in product
charts. The transverse structure to a lamination is not a manifold in
general, but rather an arbitrary (locally compact) space. The class of
laminations which most closely resemble taut foliations in their utility
are {\em essential laminations}.

Essential laminations were introduced into $3$--manifold topology in
\cite{Gabai_Oertel} as a simultaneous generalization of the concepts of
taut foliations and of incompressible surfaces.
Taut foliations can be turned into (nowhere dense)
essential laminations in a more or less trivial way by blowing up leaves
--- i.e. by replacing a leaf with a complementary $I$--bundle. {\em Genuine}
laminations are those essential laminations which do {\em not} arise from
this construction. However, it was not until
\cite{dGwK97} that the usefulness of the concept of a genuine
lamination was realized. Genuine laminations 
are characterized amongst essential laminations by
the property that some complementary region is not an $I$--bundle. 
Following \cite{Gabai_Oertel} we make this precise.

\begin{defn}
The complement of a $2$--dimensional lamination $\Lambda$ of a $3$--manifold $M$
falls into connected components called {\em complementary regions}. 
A lamination is {\em essential} if it contains no spherical leaf or torus
leaf bounding a solid torus, and furthermore, if $C$ is the metric completion
of a complementary region (with respect to the path metric on $M$), then
$C$ is irreducible, and $\partial C$ is both incompressible and {\em end
incompressible} in $C$. Here an end compressing disk is a properly embedded
$$D^2 - (\text{point in }\partial D^2) \subset C$$
which is not properly isotopic rel. $\partial$ in $C$ to an embedding in a leaf.
\end{defn}

\begin{figure}[ht]
\centerline{\relabelbox\small \epsfxsize 2.0truein
\epsfbox{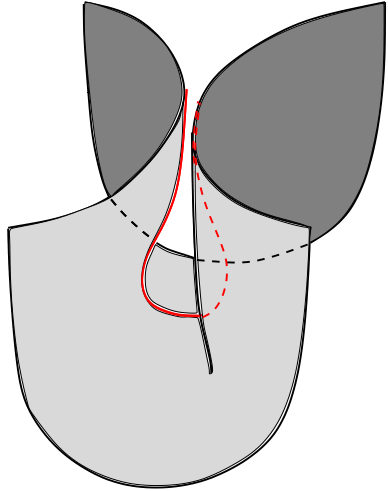}
\endrelabelbox}
\caption{An end compressing disk (boundary in red) is also
called an {\em essential monogon}.}
\label{monogon}
\end{figure}

Such an end compressing disk is also called an {\em essential ideal monogon}, or by
abuse of notation, an {\em essential monogon}. See 
figure~\ref{monogon} for the justification for this terminology. Note: some
authors prefer the term {\em essential disk-with-end}.

Another
way of phrasing the conditions above are that there should be no spherical leaf
or torus leaf bounding a solid torus, and complementary regions should contain
no essential surfaces (possibly with boundary and/or ideal points)
of positive Euler characteristic. The missing point in the boundary of an
end compressing disk should be thought of as an ideal point.
The Euler characteristic of a polygon with ideal points can be
calculated by doubling it: the double of an (ideal) monogon is a punctured sphere,
so a monogon has Euler characteristic $1/2$. An ideal bigon has Euler characteristic
$0$, an ideal triangle (a ``trigon'') has Euler characteristic $-1/2$ and so on.

A taut foliation is an example of an essential lamination. An incompressible
surface in an irreducible manifold is another example.

A complementary region to a lamination decomposes into a compact {\em gut}
piece and non--compact {\em interstitial regions} which are $I$--bundles
over non--compact surfaces. These interstitial regions are also referred to
in the literature as {\em interstitial $I$--bundles} and {\em interstices}.
These pieces meet along {\em interstitial
annuli}. Formally, the interstitial regions make up
the non--compact components of the {\em characteristic $I$--bundle}
of the complementary region (see \cite{Jaco}). For more details, see
\cite{dGwK97} or \cite{Gabai_Oertel}.

\begin{defn}
An essential lamination is {\em genuine} if some complementary region has nonempty
gut.
\end{defn}

Said another way, an essential lamination is genuine if some complementary
region contains some essential surface (possibly with boundary and/or ideal points)
of negative Euler characteristic.

The leaf space of $\til{\Lambda}$ is in general an {\em order tree}. Following
\cite{Gabai_Oertel}, an order tree can be defined as follows. 

\begin{defn}
An {\em order tree} is a set $T$ together with a collection $\mathcal{S}$ of
linearly ordered subsets called {\em segments}, each with distinct least
and greatest elements called the {\em initial} and {\em final} ends. If $\sigma$
is a segment, $-\sigma$ denotes the same subset with the reverse order, and
is called the {\em inverse} of $\sigma$. The following conditions should be
satisfied:
\begin{enumerate}
\item{If $\sigma \in \mathcal{S}$ then $-\sigma \in \mathcal{S}$}
\item{Any closed subinterval of a segment is a segment (if it has more than one element)}
\item{Any two elements of $T$ can be joined by a finite sequence of
segments $\sigma_i$ with the final end of $\sigma_i$ equal to the initial end of
$\sigma_{i+1}$}
\item{Given a cyclic word $\sigma_0\sigma_1 \cdots \sigma_{k-1}$ (subscripts
mod $k$) with the final end of $\sigma_i$ equal to the initial end of $\sigma_{i+1}$,
there is a subdivision of the $\sigma_i$ yielding a cyclic word 
$\rho_0\rho_1 \cdots \rho_{n-1}$ which becomes the trivial word when adjacent
inverse segments are cancelled}
\item{If $\sigma_1$ and $\sigma_2$ are segments whose intersection is a single element
which is the final element of $\sigma_1$ and the initial element of $\sigma_2$
then $\sigma_1 \cup \sigma_2$ is a segment}
\end{enumerate}
If all the segments are homeomorphic to subintervals of $\R$ with their order
topology, then $T$ is an {\em $\R$--order tree}.
\end{defn}

An order tree is topologized by the usual order topology on segments.
Order trees are not typically Hausdorff, but even if they are, 
there are many more possibilities than arise in the case of a foliation.

\begin{defn}
An essential lamination $\Lambda$ is {\em tight} if the leaf space of the universal
cover $\til{\Lambda}$ is Hausdorff.
\end{defn}

It follows that a taut foliation is tight iff it is $\R$--covered. Equivalently,
a lamination $\Lambda$ is tight if every arc $\alpha$ in $M$ is homotopic rel.
endpoints to an efficient arc which is either transverse or tangent to $\Lambda$.
Here an arc $\alpha$ is {\em efficient} if it does not contain a subarc
$\beta$ whose interior is disjoint from $\Lambda$, and which cobounds with
an arc $\beta'$ in a leaf of $\Lambda$ a disk
whose interior is disjoint from $\Lambda$.

If $\Lambda$ has no isolated leaves, then the associated order tree of $\til{\Lambda}$
is actually an {\em $\R$--order tree}. Any lamination can be transformed into
one without isolated leaves by blowing up isolated leaves to foliated interval bundles.
It follows that we can always consider $\R$--order trees for our applications.

Moreover, if $\Lambda$ is tight, a Hausdorff $\R$--order tree is just the underlying
topological space of an $\R$--tree. We refer to such a space as a {\em topological}
$\R$--tree to emphasize that the metric is not important. Finally, if
$\Lambda$ is a tight $1$--dimensional lamination of a surface, so that $\til{\Lambda}$
is a tight $1$--dimensional lamination of the plane, then the associated order tree
$T$ comes with a natural planar embedding, dual to $\til{\Lambda}$. See 
\cite{Gabai_Kazez_order} for more details.

Whether or not a lamination $\Lambda$ is tight, the following is true:
\begin{lem}\label{separation_constant}
Let $\Lambda$ be an essential lamination of a closed $3$--manifold $M$, and give
$M$ an arbitrary Riemannian metric. Then
there is an $\epsilon > 0$ such that every leaf $\lambda$ of $\til{\Lambda}$
is quasi--isometrically embedded in its $\epsilon$--neighborhood, and no
two incomparable leaves $\lambda,\mu$ of $\til{\Lambda}$ contain points which
are closer than $\epsilon$ in $\til{M}$.
\end{lem}
This lemma is an easy consequence of the compactness of $M$ and the defining
property of laminations, that they have local product charts. 
See \cite{dCnD02} for a proof. 

Such an $\epsilon$ is called a {\em separation constant} for $\Lambda$.

Genuine laminations certify important properties of the ambient manifold $M$.
The existence of the interstitial annuli gives a canonical collection
of knots in $M$ with important properties. Using these annuli, Gabai and Kazez
prove the following in \cite{dGwK98} and \cite{dGwK98b}. For the
definition of word--hyperbolicity of a group, see \cite{mG87}.

\begin{thm}[Gabai--Kazez \cite{dGwK98} 
Word hyperbolicity]\label{lamination_word_hyperbolic}
Let $M$ be an atoroidal $3$--manifold containing a genuine lamination $\Lambda$.
Then $\pi_1(M)$ is word--hyperbolic in the sense of Gromov.
\end{thm}

\begin{thm}[Gabai--Kazez \cite{dGwK98b} Finite MCG]\label{lamination_finite_MCG}
Let $M$ be an atoroidal $3$--manifold containing a genuine lamination $\Lambda$.
Then the mapping class group of $M$ is finite.
\end{thm}

Laminations come in all degrees of smoothness. Moreover, it is important to
distinguish between the smoothness of {\em leaves} and the smoothness of
the {\em transverse space}. For some applications in this paper, it will be
important for our laminations to be leafwise smooth. Fortunately, the situation
for $2$--dimensional laminations in $3$--manifolds is as simple as it could be. The
main theorem of \cite{dC01} is the following:

\begin{thm}\label{leafwise_smooth}
Let $\Lambda$ be a $2$--dimensional lamination in a smooth $3$--manifold $M$. Then
$\Lambda$ is isotopic to a lamination with smoothly immersed leaves.
\end{thm}

\begin{rmk}
$2$--dimensional laminations are also sometimes referred to informally
as {\em surface laminations}.
\end{rmk}

Amongst all genuine laminations, some are more useful than others.
If $M$ is not Haken, then \cite{Hatcher_Oertel} show that the gut regions
are all homeomorphic to handlebodies. They call such laminations {\em full}, where
the terminology is meant to imply that the complementary regions contain no
closed incompressible surface. Specializing further, we have the following.

\begin{defn}\label{very_full_genuine}
A genuine lamination is {\em very full} if all complementary regions are
finite sided ideal polygon bundles over $S^1$.
\end{defn}

The relationship between the topology of the guts and the topology of the
complementary regions is not straightforward in general. However, in the
case of a lamination with solid torus guts, the following lemma is proved in
\cite{dCnD02}:

\begin{lem}[Calegari--Dunfield \cite{dCnD02} Filling Lemma]\label{filling_lemma}
Let $\Lambda$ be a genuine lamination of a closed $3$--manifold $M$ with solid
torus guts. Then $\Lambda$ is a sublamination of a very full genuine lamination
$\overline{\Lambda}$.
Moreover, if $\Lambda$ is tight, so is $\overline{\Lambda}$.
\end{lem}

Very full genuine laminations are particularly nice. There is the following
theorem of Gabai and Kazez from \cite{dGwK97}:

\begin{thm}[Gabai--Kazez]\label{homotopy_implies_isotopy}
Let $M$ be a $3$--manifold with a very full genuine lamination $\Lambda$.
Then any self--homeomorphism of $M$ homotopic to the identity is isotopic
to the identity.
\end{thm}

{\em Tight} very full genuine laminations have another application, more central
to the theme of this paper. In \cite{dCnD02} Calegari and Dunfield show
that they give rise to a {\em universal circle}, as alluded to 
in \S~\ref{introduction}. In this paper, the very full genuine laminations we
produce, although not necessarily tight,
already come with the data of a universal circle, so the
construction in \cite{dCnD02} is superfluous for our purposes.

\subsection{Candel's uniformization theorem}\label{Candel_uniformization}

The classical uniformization theorem says that every Riemann surface
is conformally equivalent to a surface with a complete metric of constant
curvature $1,0$ or $-1$. Such a conformal equivalence is called a
{\em uniformizing map}.
If the curvature of the metric is negative, then the surface is said to
be of {\em hyperbolic type}, and the metric is unique.

Riemann surface laminations can be uniformized leafwise; an important
question is to determine when these choices of uniformizing maps are
all of the same type, and can be chosen in a continuously varying way.

For Riemann surfaces with all leaves hyperbolic, the complete answer is
given by a theorem of Candel \cite{aC93}:

\begin{thm}[Candel \cite{aC93} Uniformization for 
hyperbolic laminations]\label{Candel_uniformize}
Let $\Lambda$ be a $2$--dimensional lamination with a leafwise Riemannian metric such
that every leaf is of hyperbolic type.
Then the leafwise constant curvature hyperbolic metric determined uniquely by the conformal
structure of the leaves of $\Lambda$ varies {\em continuously} in the
transverse direction.
\end{thm}

Let $M$ be an atoroidal $3$--manifold, and $\Lambda$ an essential
lamination. Then the leaves of $\Lambda$ are all of hyperbolic type 
(see e.g. \cite{dC99b}) and therefore Candel's theorem applies. If
$\Lambda$ is an essential surface lamination with smooth leaves, 
then there is a $C^0$ Riemannian metric on $M$ (i.e. a continuous
section of the bundle of positive definite symmetric $2$--tensors) 
which restricts on each leaf
to a Riemannian metric of constant curvature $-1$. In general, if $\Lambda$ is 
transversely $C^n$, then the Riemannian metric on $M$ can be chosen to be $C^n$.

Using Candel's theorem, we can define the {\em circle bundle at infinity}
of an essential lamination.

\begin{defn}
Let $\Lambda$ be an essential lamination of $M$, and let $L$ be the leaf space
of $\til{\Lambda}$. Uniformize $\Lambda$ by theorem~\ref{Candel_uniformize}, so
that every leaf of $\til{\Lambda}$ is isometric to $\H^2$.

For each leaf $\lambda$ of $\til{\Lambda}$, let $S^1_\infty(\lambda)$ denote
the ideal boundary of $\lambda$ with respect to this metric. The {\em endpoint map}
$$e: UT_p\lambda \to S^1_\infty(\lambda)$$
takes a unit vector $v$ in $\lambda$ at $p$ to the endpoint at infinity of the
geodesic ray $\gamma_v \subset \lambda$ which emanates from $p$, and
satisfies $\gamma_v'(0) = v$.
\end{defn}

\begin{defn}
Let $\Lambda$ be an essential lamination of $M$, and let $L$ be the
leaf space of the pulled back lamination in the universal cover $\til{\Lambda}$.
Uniformize $\Lambda$ by theorem~\ref{Candel_uniformize}, so that every leaf of
$\til{\Lambda}$ is isometric to $\H^2$.
For each leaf $\lambda$ of $\til{\Lambda}$, let $S^1_\infty(\lambda)$ denote the
circle at infinity of $\lambda$, under its isometric identification with $\H^2$.
The {\em circle bundle at infinity} is the topological space whose underlying set
is the disjoint union
$$E_\infty = \bigcup_{\lambda \in L} S^1_\infty(\lambda)$$
and with the smallest topology so that the endpoint map
$$e:UT\til{\Lambda} \to E_\infty$$
is continuous.
\end{defn}

With this topology, $E_\infty$ is a {\em circle bundle} over $L$, whose fiber
over each $\lambda \in L$ is $S^1_\infty(\lambda)$.

Note that by Candel's theorem~\ref{Candel_uniformize}, 
for every efficient transversal $\tau$ to $\til{\Lambda}$, the restriction
$$e:UT\til{\Lambda}|_\tau \to E_\infty|_\tau$$
is a homeomorphism.

\subsection{Minimal sets}

Given a lamination $\Lambda$ a {\em minimal set} is a sublamination, defined
as follows:

\begin{defn}
Let $\Lambda$ be a lamination of a compact manifold.
A {\em minimal set} $\Lambda_m \subset \Lambda$ is a subset of $\Lambda$ which is
minimal with respect to inclusion, and satisfies the following properties:
\begin{enumerate}
\item{$\Lambda_m$ is nonempty.}
\item{$\Lambda_m$ is {\em saturated}. That is, it is a union of leaves of
$\Lambda$.}
\item{$\Lambda_m$ is closed.}
\end{enumerate}
\end{defn}

Minimal sets always exist. In fact, a nonempty sublamination $\Lambda_m$
is minimal iff every leaf is dense. A lamination which does not satisfy
this property contains some leaf $\lambda$ which is not dense; the closure
of $\lambda$ is a smaller sublamination. By transfinite induction, the closure
of every leaf contains a minimal set. Note that this construction uses the
axiom of choice.

Let $\F$ be a taut foliation of $M$. If $\F$ is not minimal, by the definition
of essentiality, every minimal set $\Lambda$ is an essential lamination. 
Such a lamination is either
{\em genuine}, or else all complementary regions are products. In this case,
either $\Lambda$ is a single fiber of a fibration of $M$ over $S^1$, or else
such complementary regions can be {\em collapsed} to give a new taut foliation
$\F'$ of $M$ which is minimal.

This collapsing procedure is the converse of the operation of {\em blowing--up}
or {\em Denjoying} a leaf. A thorough discussion is contained in \cite{CCI},
so we do not elaborate here. We summarize the discussion above in the following
lemma:

\begin{lem}\label{minimal_quotient}
Let $\F$ be a taut foliation of a $3$--manifold $M$. Then either
$\F$ contains a genuine sublamination, or a fiber of a fibration over $S^1$, or
else $M$ contains a taut foliation $\F'$ with every leaf dense.
\end{lem}

\section{Universal circles for taut foliations}\label{constructing_circle_section}

In this section we define a universal circle for a taut foliation, and
give an outline of the construction of a universal circle in \cite{dCnD02}.

\subsection{Definition of a universal circle}

\begin{defn}
Let $\F$ be a taut foliation of an atoroidal $3$--manifold $M$. A {\em universal
circle} for $\F$ is a circle $S^1_\u$ together with the following data:
\begin{enumerate}
\item{There is a faithful representation $$\rho_\u:\pi_1(M) \to \homeo^+(S^1_\u)$$}
\item{For every leaf $\lambda$ of $\til{\F}$ there is a monotone map
$$\phi_\lambda:S^1_\u \to S^1_\infty(\lambda)$$
Moreover, the map
$$\phi:S^1_\u \times L \to E_\infty$$
defined by $\phi(\cdot,\lambda) = \phi_\lambda(\cdot)$ is continuous. That is,
$(E_\infty,L,\phi)$ is a monotone family.}
\item{For every leaf $\lambda$ of $\til{\F}$ and every $\alpha \in \pi_1(M)$ the
following diagram commutes:
$$\begin{CD}
S^1_\u @>\rho_\u(\alpha)>> S^1_\u \\
@V\phi_{\lambda}VV @V\phi_{\alpha(\lambda)}VV \\
S^1_\infty(\lambda) @>\alpha>> S^1_\infty(\alpha(\lambda))
\end{CD}$$}
\item{If $\lambda$ and $\mu$ are incomparable leaves of $\til{\F}$ then the core
of $\lambda$ is contained in the closure of a single gap of $\mu$ and
{\it vice versa}.}
\end{enumerate}
\end{defn}

\begin{thm}[Thurston, Calegari--Dunfield \cite{dCnD02} Universal 
circles for foliations]\label{universal_circle_exists}
Let $\F$ be a co--oriented taut foliation of an atoroidal, oriented
$3$--manifold $M$. Then there is a universal circle for $\F$.
\end{thm}

In section \S~\ref{constructing_laminations}, 
we will see how the axiomatic definition of a universal
circle lets us construct transverse very full genuine laminations. However,
to analyze the properties of these laminations in more detail, we need to know
more about the construction of the universal circle. 

\subsection{Markers}\label{marker_subsection}

\begin{defn}
Let $\Lambda$ be an essential lamination of $M$ with hyperbolic leaves.
A {\em marker} for $\Lambda$ is a map
$$m:I \times \R^+ \to \til{M}$$ with the following properties:
\begin{enumerate}
\item{There is a closed set $K \subset I$ such that for each $k \in K$, the
image of $k \times \R^+$ in $\til{M}$ is a geodesic ray in a leaf of $\til{\Lambda}$.
Further, for $k \in I \backslash K$,
$$m(k \times \R^+) \subset \til{M} \backslash \til{\Lambda}$$
We call these rays the {\em horizontal rays} of the marker.}
\item{For each $t \in \R^+$, the interval $m(I \times t)$ is a tight transversal.
Further, there is a separation constant $\epsilon$ for $\Lambda$, such that
$$\text{length}(m(I \times t)) < \epsilon/3$$
We call these intervals the {\em vertical intervals} of the marker.}
\end{enumerate}
\end{defn}

For a marker $m$, a horizontal ray $m(k \times \R^+)$ in a leaf $\lambda$ of
$\til{\Lambda}$ is asymptotic to a unique point in $S^1_\infty(\lambda)$, which
we call the {\em endpoint} of $m(k \times \R^+)$. By abuse of notation, we
call the union of such endpoints, as $k$ varies over $K$, the {\em endpoints of
the marker $m$}.

If $\F$ is a foliation, then $K = I$ for each marker $m$, and the set of endpoints of
$m$ define an embedded interval in $E_\infty$ transverse to the foliation by circles.

Markers are related to, and arise in practice from {\em sawblades}, defined as follows:
\begin{defn}\label{sawblade_definition}
Let $\Lambda$ be an essential lamination of $M$ with hyperbolic leaves.
An {\em $\epsilon$--sawblade} for $\F$ is an embedded polygonal surface $P \subset M$
obtained from a square by gluing the right hand edge to a subset of the left hand
edge in such a way that the lowermost vertices are identified. In co--ordinates:
if we parameterize the square as $[0,1]\times [0,1]$ and denote the image of $P$ as
$P([0,1],[0,1])$, then under the identification, $P(1,[0,1])$ gets identified with a
subset $P(0,[0,t])$ with $0 < t \le 1$. Moreover, $P$ must satisfy the following properties:
\begin{enumerate}
\item{There is a closed subset $K \subset I$ including the endpoints of $I$, such that
for each $t \in K$, the subset $P([0,1],t) \subset \F$ is a geodesic arc in
a leaf $\lambda_t$ of $\Lambda$. For $t = 0$, the subset $P([0,1],0) \subset \F$ closed
up to a geodesic loop $\gamma \subset \lambda_0$.}
\item{For each $t \in [0,1]$, the subset $P(t,[0,1])$ is an embedded, tight transversal
to $\Lambda$ of length $\le \epsilon$. The transversal $P(1,[0,1])$ is contained in
the image of $P(0,[0,1])$, and the corresponding geodesic segments $P([0,1],t_1)$ and 
$P([0,1],t_2)$, where $P(1,t_1) = P(0,t_2)$ for $t_1 \in K$, join up to a geodesic segment in
the corresponding leaf of $\Lambda$; i.e. there is no corner along $P(0,[0,1])$.}
\end{enumerate}
\end{defn}
If $\epsilon$ is understood, we just say a {\em sawblade}.
Note that holonomy transport of the transversal $P(0,[0,1])$ around $\gamma$ induces
an embedding $K \to K$ taking one endpoint to itself. Here $\gamma$ is oriented
compatibly with the usual orientation on $I = [0,1]$. We call the positive direction
on $\gamma$ the {\em contracting direction} for the sawblade, and the
negative direction the {\em expanding direction}.

We show how to construct a marker from a sawblade.

\begin{construct}\label{sawblade_marker}
Let $P$ be a sawblade, and let $\til{P}$ be a component of the preimage
in $\til{M}$. $\til{P}$ is the universal cover of $P$, and the deck group of
the cover is $\pi_1(P) = \Z$, generated by the closed geodesic $\gamma$ as
in definition~\ref{sawblade_definition}.

Let $\tau$ be a lift of $P(0,I)$, and let
$K \subset I$ be as in definition~\ref{sawblade_definition}. 
Parameterize $\tau$ as $\tau(t)$ where $\tau(t)$ corresponds to the
lift of $P(0,t)$. Then for each $k \in K$, let $\lambda_k$ denote the
leaf of $\til{\Lambda}$ containing $\tau(k)$. By the second property of
a sawblade, the intersection $\lambda_k \cap \til{P}$ contains an
entire geodesic ray starting from $\tau(k)$. Together with complementary strips
of $\til{P}$, the union of these rays are a marker for $\Lambda$.
\end{construct}

Notice that the union of the markers constructed in construction~\ref{sawblade_marker},
over all lifts $\tau$ of $P(0,I)$, is exactly the preimage $\til{P}$.

By abuse of notation, we refer to the union of the
endpoints of the markers associated
to $\til{P}$ in construction~\ref{sawblade_marker}
as the {\em endpoints} of $\til{P}$.

Every closed geodesic $\gamma$ contained in a non--simply connected
leaf $\lambda$ of $\Lambda$ is the boundary geodesic of some $\epsilon$--sawblade,
for any positive $\epsilon$. If
$M$ is a closed $3$--manifold containing an essential lamination $\Lambda$ with
every leaf simply connected, then $M$ is $T^3$ (see \cite{tL02} for an elegant proof).
It follows that if $M$ is atoroidal, many sawblades can be constructed.
 
Let $P$ be a sawblade, and suppose $\Lambda$ is minimal. Then there is a uniform
constant $C$ such that for every leaf $\lambda$ of $\til{\Lambda}$, and for every
point $p \in \lambda$, there is a point $q \in \lambda$ within distance
$C$ in the path metric on $\lambda$, such that $q$ is contained in a lift of $P$.
It follows that $q$ is contained in a marker, and there is a geodesic ray $r$ through
$q$ such that holonomy transport of a suffiently short transversal $\tau(q)$ through
$q$ along $r$ keeps the length of the transversal smaller than $\epsilon/3$ for
all time. It is not hard from this to conclude the following lemma, proved in
\S 5.6 of \cite{dCnD02}:

\begin{lem}[Calegari--Dunfield]\label{sawblade_lifts_dense}
Let $\Lambda$ be a minimal essential lamination of an atoroidal $3$--manifold $M$, and
let $P$ be an $\epsilon$--sawblade for $\Lambda$. Then the set of
endpoints of lifts $\til{P}$ of $P$ is dense in $S^1_\infty(\lambda)$ for
every leaf $\lambda$ of $\til{\Lambda}$.
\end{lem}

From this lemma, it is not hard to conclude the following theorem, called
the {\em Leaf Pocket Theorem}:

\begin{thm}[Calegari--Dunfield \cite{dCnD02} Leaf Pocket theorem]\label{leaf_pocket_theorem}
Let $\Lambda$ be an essential lamination on an atoroidal $3$--manifold $M$.
Then for every leaf $\lambda$ of $\til{\Lambda}$, and every $\epsilon > 0$,
the set of endpoints of $\epsilon$--markers is dense in $S^1_\infty(\lambda)$.
\end{thm}

As explained above, each marker $M$ defines by the endpoint map, an embedded 
interval $e(M) \subset E_\infty$
transverse to the foliation by circles.

The following lemma is a restatement of lemma 6.11 in \cite{dCnD02}:

\begin{lem}[Calegari--Dunfield]
Let $e(m_1),e(m_2)$ be two endpoint intervals of markers $m_1,m_2$.
Then these intervals are either disjoint, or else their union is an embedded, ordered
interval transverse to the foliation of $E_\infty$ by circles.
\end{lem}

It follows that distinct markers whose endpoints intersect can be amalgamated,
and the unions give a $\pi_1(M)$--invariant
family of disjoint, embedded intervals in $E_\infty$ transverse to the foliation
by circles. We denote this family of intervals by $\M$, and denote a typical
element of $\M$ by $m$, or $m_i$ for some index $i$.

\subsection{Special sections}

In this section we indicate how to go from theorem~\ref{leaf_pocket_theorem} to
a proof of theorem~\ref{universal_circle_exists}.

\begin{defn}
Let $p \in S^1_\u$. The {\em special section} associated to $p$ is a section
$$\sigma_p:L \to E_\infty$$
defined by
$$\sigma_p(\lambda) = \phi_\lambda(p)$$
\end{defn}

The strategy in \cite{dCnD02} in constructing $S^1_\u$ is to construct
sufficiently many non--crossing sections $L \to E_\infty$, show that there
is a natural $\pi_1(M)$--invariant circular ordering on this set $\SS$ 
of sections, and define $S^1_\u$ to be the order completion of $\SS$.

For each leaf $\lambda$ of $\til{\F}$, and each point $p \in S^1_\infty(\lambda)$,
we construct a section $s_p:L \to E_\infty$ satisfying
$$s_p(\lambda) = p$$

The combinatorics of the construction are somewhat complicated, and involve
careful attention to orientations of transversals and circles in the
oriented circle bundle $E_\infty$.

Some simple cases, which nevertheless give an idea of what is going on, 
are illustrated in \S~\ref{example_circles_subsection}.
We suggest that the reader unfamiliar with \cite{dCnD02} go back and forth between
that subsection and this, in order to get a clear idea of our conventions.
Of course, for details, consult \cite{dCnD02}.

\begin{defn}
Let $I \subset L$ be an embedded interval. A section $\tau:I \to E_\infty|_I$
is {\em admissible} if its image does not cross any element $m \in \M$
transversely.
\end{defn}

Now suppose $I$ is an {\em oriented} interval in $L$, and let $\lambda_t$
with $t \in I$ denote the corresponding leaves of $\til{\F}$. Given
$p \in S^1_\infty(\lambda_t)$, a section $\tau:I \to E_\infty|_I$
with $\tau(0) = p$ is {\em leftmost} if it is never to the right of
any other admissible section $\tau'$ with $\tau'(0) = p$.

Here the orientation on the interval $I$ and on the cylinder $E_\infty|_I$
determine the meaning of the ``left'' and ``right'' sides of $\tau(I)$ and
$\tau'(I)$. Note that if the orientation on $I$ is reversed, then the left and
right sides are reversed too.

Since they agree in $S^1_\infty(\lambda_0)$, where they first
start to diverge, it makes sense to say that $\tau(I)$ stays to the left of
$\tau'(I)$, and it makes sense to say that thereafter $\tau(I)$ never crosses
$\tau'(I)$ from the left side. That is, if $t$ is a local maximum for the
subset of $I$ for which $\tau(t) = \tau'(t)$, then for all
$s>t$ with $s-t$ sufficiently small, $\tau(s)$ must be to the left of $\tau'(s)$.

Notice that if the orientation on $I$ agrees with the orientation on $L$ coming
from the co--orientation of $\F$, then leftmost admissible sections 
over $I$ are {\em clockwisemost}.
Conversely, if the orientation on $I$ disagrees with the orientation on $L$,
then leftmost admissible sections over $I$ are {\em anticlockwisemost}.

Now, let $l \subset L$ be a properly embedded copy of $\R$, intersecting
$\lambda$, and let $p \in S^1_\infty(\lambda)$. $l \backslash p$ consists of
two rays $l^\pm$. Give $l^+$ the usual orientation, agreeing with the order structure
on $L$, but reverse the orientation on $l^-$, so that it points in the {\em negative}
direction. Then define $s_p|_{l^\pm}$ to be leftmost admissible section with
$s_p(\lambda) = p$. Notice that $s_p|_{l^+}$ is {\em clockwisemost} and
$s_p|_{l^-}$ is {\em anticlockwisemost} amongst admissible sections, with
respect to the global order structure on $L$ and $E_\infty$.

This defines the section $s_p$ over the union of leaves which are comparable with
$\lambda$. Now, suppose $\mu_1,\mu_2$ are two leaves such that $\mu_1 < \lambda$,
and with the additional
property that there is a $1$--parameter family of leaves
$\nu_t$ with $t \in [0,1)$, such that
$\nu_t < \nu_s$ for $t < s$, and both $\mu_1$ and $\mu_2$ are positive
limits (in $L$) of $\nu_t$ as $t \to 1$.

Consider the union of the intervals $m_{1i},m_{2i}$ in $\M$ which intersect 
$S^1_\infty(\mu_1)$ and $S^1_\infty(\mu_2)$ respectively. As
$t \to 1$, more and more of the $m_{1i},m_{2i}$ intersect $S^1_\infty(\nu_t)$.
Since they are disjoint, they inherit a circular ordering as follows:
if $a,b,c$ are three such elements of $\M$, then there is a $t$ such that
all three intersect $S^1_\infty(\nu_t)$. Then the cyclic order on $a,b,c$
is just the cyclic order of the intersections $a,b,c \cap S^1_\infty(\nu_t)$.

Moreover, it is not hard to show that for any choice of
$m_{11},m_{12}$ and $m_{21},m_{22}$, the two unordered pairs cannot link in
any $S^1_\infty(\nu_t)$, and therefore all the $m_{1i}$ are contained in a
``gap'' of the circularly ordered set of $m_{2j}$. Since the $m_{2j}$ are
dense in $S^1_\infty(\mu_2)$, this gap defines a unique point $q \in S^1_\infty(\mu_2)$.
We can now define $s_p(\mu_2) = q$. We call this the method of {\em turning corners},
since it shows how to continue a leftmost section $s_p$ across nonseparated leaves
in $L$.

Since $L$ is simply connected, the Hausdorffification is a (topological) $\R$--tree $T$.
Given a point $\lambda \in L$, there is a unique embedded segment in $T$ from $\lambda$
to any other point. Back in $L$, this defines a sequence of leaves
$$\lambda = \lambda_0,\lambda_1,\lambda_2,\dots,\lambda_n = \mu$$
for every pair $\lambda,\mu$,
where $\lambda_{2i}$ and $\lambda_{2i+1}$ are comparable, and $\lambda_{2i-1}$
and $\lambda_{2i}$ are incomparable and nonseparated.

We define $s_p$ inductively, by leftmost sections over the oriented subintervals
$[\lambda_{2i}\lambda_{2i+1}]$ of $L$, and by the method of turning corners to go from
the section at $\lambda_{2i-1}$ to the section at $\lambda_{2i}$. Notice, in fact,
that the value of $s_p(\lambda_{2i})$ does not depend at all on the values of
$s_p$ on $\lambda_j$ with $j < 2i$, and depends only on the {\em leaf}
$\lambda_{2i-1}$.

We let $\SS$ denote the union of all the special sections $s_p$ as above, as $p$
ranges over $E_\infty$, and we identify $s_p$ and $s_q$ if they define the same
section $L \to E_\infty$.

For each leaf $\lambda$ of $\til{\F}$, let $\SS(\lambda)$ denote the set of
special sections $s_p$ where $p \in S^1_\infty(\lambda)$.

Finally, in lemma 6.25 of \cite{dCnD02} it is shown that the set $\SS$ of sections
$s_p$ as above is naturally circularly ordered. It follows that we can take
the order completion $\overline{\SS}$, which is homeomorphic with the order topology
to a closed subset of a circle. By collapsing the complementary gaps in this image,
we get a universal circle, which we call $S^1_\u$. That is, we have
$$\SS \to \overline{\SS}  \to S^1_\u$$
where the first map is an inclusion, and the second is a surjection.
Notice that the map $\overline{\SS} \to S^1_\u$ is at most $2$--$1$, and is $1$--$1$
away from countably many points. It follows that the natural inclusion
$\SS(\lambda) \to \SS$ extends to an {\em inclusion} $\SS(\lambda) \to S^1_\u$.

The map
$\phi_\lambda$ is defined on $s_p \in S^1_\u$ by
$$\phi_\lambda(s_p) = s_p(\lambda)$$
and these maps are collated, by varying over $\lambda \in L$, to $\phi$.
It is clear that $\phi:S^1_\u \times L \to E_\infty$ defined in this
way is continuous, and that $\phi_\lambda$ is monotone for each $\lambda$.

Moreover, it is clear that the natural action of $\pi_1(M)$ on $E_\infty$
induces an action on $\SS$ preserving the circular order, and therefore
induces a representation
$$\rho_\u:\pi_1(M) \to \homeo^+(S^1_\u)$$

\subsection{Examples of universal circles}\label{example_circles_subsection}

In this subsection we give some idea of the combinatorics of
universal circles.

\begin{exa}[Linear segment]\label{linear_segment_example}
Let $I \subset L$ be a closed interval, with lowest leaf $\lambda$ and highest
leaf $\lambda'$. Leftmost trajectories can run into each other, but not cross.
A leftmost ascending trajectory can coalesce with a leftmost descending trajectory. 
The set of special sections give the cylinder $E_\infty|_I$ the structure of
a ($1$--dimensional) {\em branched lamination}; see 
definition~\ref{branched_lamination_definition} for a general definition.

In the universal circle, the set of special sections which intersect
$\lambda$ at $x$ and $\lambda'$ at $x'$ is an interval, running positively
from $s_x$ to $s_{x'}$.

Here is another way to see the circular order on special sections in $I$.
Lift to the universal cover of the cylinder $E_\infty|_I$. Each
special section lifts to $\Z$ copies in the cover. In the cover, two
sections $s_y,s_z$ satisfy
$s_y < s_z$ iff there is a nontrivial positive transversal from
$s_y$ to $s_z$. This defines a total order upstairs, which is evidently
order isomorphic to $\R$. The action of the deck group on the cover of
the cylinder induces an action on the ordered set of lifts of special sections,
inducing a circular order on their quotient.

\begin{figure}[ht]
\centerline{\relabelbox\small \epsfxsize 3.0truein
\epsfbox{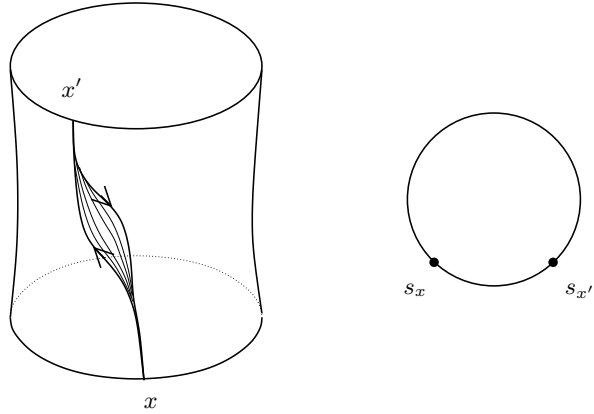}
\relabel {x}{$x$}
\relabel {x'}{$x'$}
\relabel {sx}{$s_x$}
\relabel {sx'}{$s_{x'}$}
\endrelabelbox}
\caption{The special sections might coalesce, but they don't cross.}
\label{linear_segment}
\end{figure}
\end{exa}

\begin{exa}[Nonseparated leaves]\label{nonseparated_leaves_example}
The next example incorporates positive branching.
Let $\lambda,\mu$ be two incomparable leaves which are nonseparated, and
such that there is a $1$--parameter family of leaves $\nu_t$ with 
$t \in [0,1)$, satisfying $\nu_t < \lambda,\mu$
for all $t$, and converging to both $\lambda$ and $\mu$ as $t \to 1$.

Every marker which intersects $\lambda$ or
$\mu$ will intersect $\nu_t$, for sufficiently large $t$.
As described in the previous subsection, this induces a circular order on
the union of a
dense subset of $S^1_\infty(\lambda)$ and $S^1_\infty(\mu)$, and by
comparing special sections in $S^1_\infty(\nu_t)$ for sufficiently large
$t$, these can be completed to
a circular order on the disjoint union
of {\em all} special sections $s_x$ where
$$x \in S^1_\infty(\lambda) \cup S^1_\infty(\mu)$$

In this circularly ordered set
the set of special sections $s_x$ with $x \in S^1_\infty(\lambda)$
is a half--open interval, containing a (locally) clockwisemost point, but
not a (locally) anticlockwisemost point, and similarly for the $s_y$ with
$y \in S^1_\infty(\mu)$.

Notice that if $\lambda,\mu$ were nonseparated, but the approximating sequence
$\nu_t$ satisfied $\nu_t > \lambda,\mu$ then the half--open intervals
of special sections would contain (locally) anticlockwisemost points instead.

\begin{figure}[ht]
\centerline{\relabelbox\small \epsfxsize 3.0truein
\epsfbox{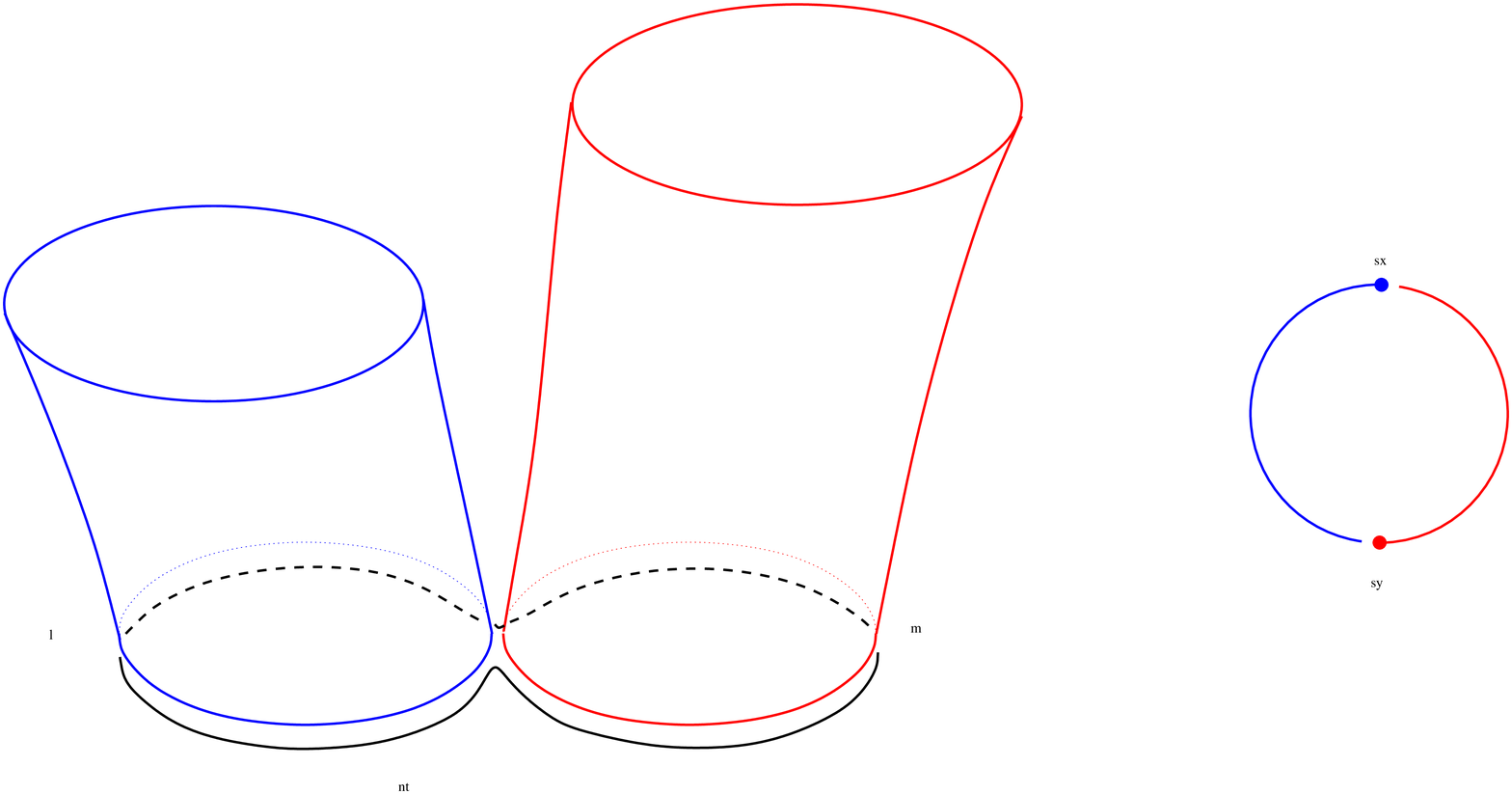}
\relabel {l}{$\mu$}
\relabel {m}{$\lambda$}
\relabel {nt}{$\nu_t$}
\relabel {sx}{$s_y$}
\relabel {sy}{$s_x$}
\endrelabelbox}
\caption{The special sections coming from each of the two nonseparated leaves
determine a half--open interval in the circular order on the union. Here, the
point $x$ is in $S^1_\infty(\lambda)$, and the point $y$ is in $S^1_\infty(\mu)$.}
\label{nonseparated_leaves}
\end{figure}
\end{exa}

\begin{exa}[More branching]\label{more_branching_example}
The next example includes both positive and negative branching.
In this case, we have nonseparating leaves $\mu,\lambda$ exhibiting positive
branching, nonseparating leaves $\nu,\lambda'$ exhibiting negative branching,
where $\lambda' > \lambda$. Let $x \in S^1_\infty(\lambda)$ be the point
determining the locally clockwisemost segment $s_x$ in the previous example,
and let $x'$ be the corresponding point (determining the locally 
anticlockwisemost segment)
in $S^1_\infty(\lambda')$.

There are two topologically distinct cases to consider: in the first, the special sections
$s_x$ and $s_{x'}$ do not agree on the entire interval $[\lambda,\lambda']$, although they
might agree on some closed subset of this interval, which might include either
or both of the endpoints.
In the second, the sections $s_x,s_{x'}$ {\em do} agree on the entire interval, and
therefore agree on all of $E_\infty$.
\end{exa}

These examples contain all the necessary information to show how to go from a
finite union $K$ of ordered subsegments in $L$, whose image in the Hausdorffification
of $L$ is connected, to a circle $S^1(K)$
which realizes the circular order on the set of
special sections associated to points in leaves $\lambda$ in $K$. By following
the model of example~\ref{nonseparated_leaves_example}, one can amalgamate the
circles associated to a pair of ordered segments whose endpoints are nonseparated.
Given $K_i,K_j$ disjoint, finite connected unions, we get circles $S^1(K_i)$ and
$S^1(K_j)$; if $K_i$ and $K_j$ contain a pair of nonseparated leaves, we can
follow example~\ref{nonseparated_leaves_example} to amalgamate $S^1(K_j)$ and
$S^1(K_i)$ into $S^1(K_i \cup K_j)$, completing the induction step. One must
verify that the result does not depend on the order in which one constructs
$K$ from ordered subsegments; implicitly, this is a statement about the
{\em commutativity} of the amalgamating operation in example~\ref{nonseparated_leaves_example}.
This commutativity is evident even in example~\ref{more_branching_example},
where one may choose to amalgamate the segment $[\lambda,\lambda']$ with
$\mu$ first and then $\nu$, or the other way around.

\begin{figure}[ht]
\centerline{\relabelbox\small \epsfxsize 3.5truein
\epsfbox{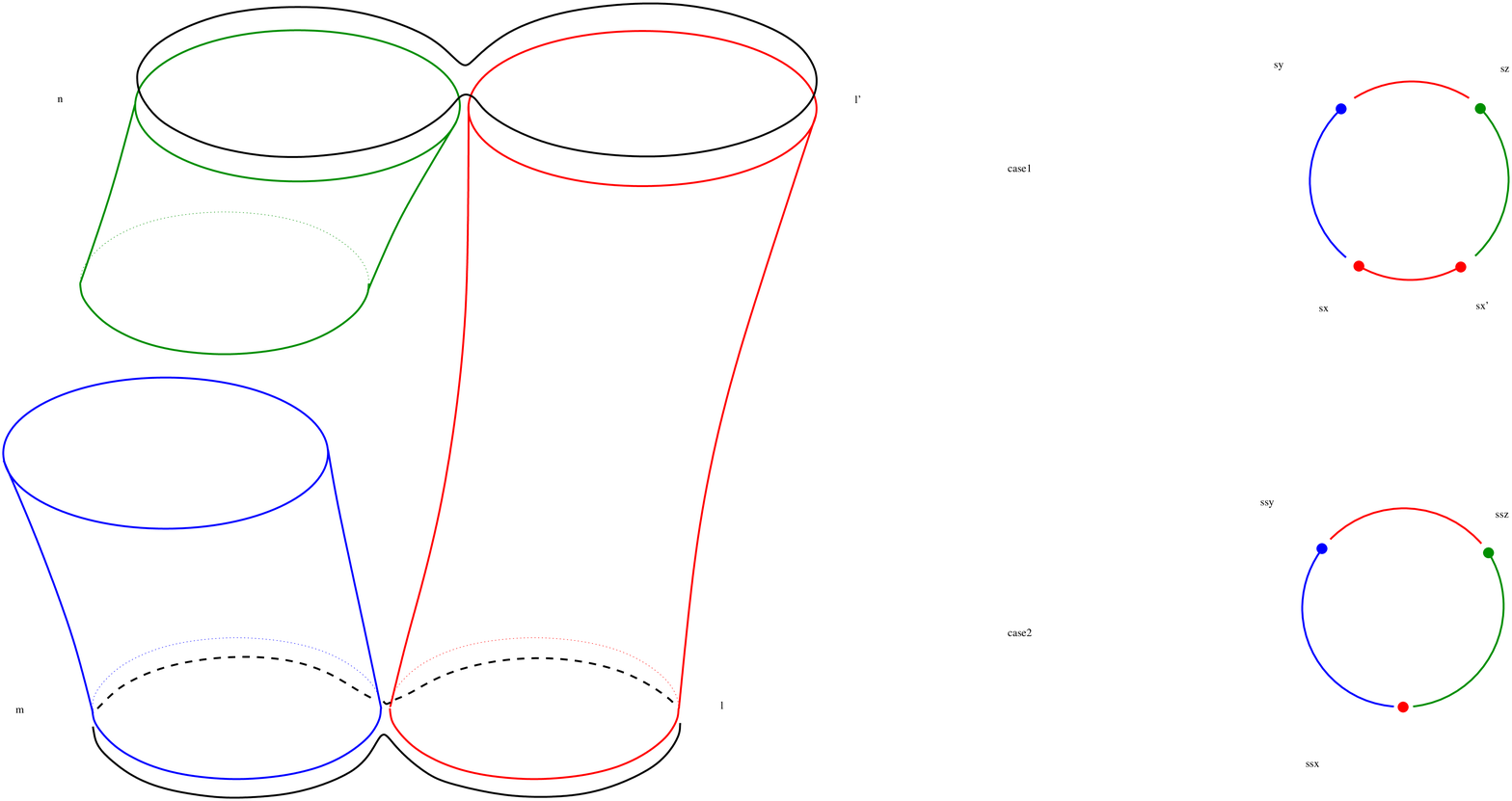}
\relabel {l}{$\lambda$}
\relabel {l'}{$\lambda'$}
\relabel {m}{$\mu$}
\relabel {n}{$\nu$}
\relabel {sx}{$s_x$}
\relabel {sx'}{$s_x'$}
\relabel {sy}{$s_y$}
\relabel {sz}{$s_z$}
\relabel {ssx}{$s_x = s_{x'}$}
\relabel {ssy}{$s_y$}
\relabel {ssz}{$s_z$}
\relabel {case1}{$\text{case 1:}$}
\relabel {case2}{$\text{case 2:}$}
\endrelabelbox}
\caption{In case $1$, $s_x$ and $s_{x'}$ differ somewhere on $E_\infty|_{[\lambda,\lambda]}$.
In case $2$, they are equal on all of $E_\infty$.}
\label{more_branching}
\end{figure}

\subsection{Special sections and cores}

Recall the notation $\SS(\lambda)$ to denote the set of special sections
associated to points $x \in S^1_\infty(\lambda)$. In this subsection we 
describe the relationship between $\SS(\lambda)$ and the core of $\phi_\lambda$.

\begin{lem}\label{sections_dense_in_core}
Let $\lambda$ be a leaf of $\til{\F}$. Then $\core(\phi_\lambda)$ is
contained in the closure $\overline{\SS(\lambda)} \subset S^1_\u$,
and the difference $\overline{\SS(\lambda)} \backslash \core(\phi_\lambda)$ consists
of at most countably many isolated points, at most one in each gap of $\phi_\lambda$.
\end{lem}
\begin{proof}
Given $p,q \in \core(\phi_\lambda)$,
either $p$ and $q$ are the boundary points of the closure of some gap, or
else $\phi_\lambda(p) \ne \phi_\lambda(q)$, and therefore there are $p',q' \in \SS(\lambda)$
which link $p,q$. It follows that every accumulation point of $\core(\phi_\lambda)$
is an accumulation point of $\SS(\lambda)$. Since $\core(\phi_\lambda)$ is
perfect, it follows that $\core(\phi_\lambda) \subset \overline{\SS(\lambda)}$.

Conversely, given $p,q \in S^1_\infty(\lambda)$ distinct points, we have
$\phi_\lambda(p) = p \ne q = \phi_\lambda(q)$, and therefore there are points
$p',q' \in \core(\lambda)$ which link $p,q$. In particular, $p$ and $q$ are not
both in the same gap region of $\phi_\lambda$, and therefore there is at most
one such point in each gap. Since $\phi_\lambda$ has only countably many gaps,
the lemma follows.
\end{proof}

An example where $\overline{\SS(\lambda)} \backslash \core(\phi_\lambda)$ might
contain isolated points is illustrated in figure~\ref{more_branching}.

Now, if $\lambda$ and $\mu$ are incomparable leaves, then $\phi_\mu(\SS(\lambda))$
is a single point of $S^1_\infty(\mu)$, and similarly for $\phi_\lambda(\SS(\mu))$.
Since $\phi_\lambda$ is $1$--$1$ on $\SS(\lambda)$, it follows that
$\SS(\lambda)$ and $\SS(\mu)$ are not linked as subsets of $S^1_\u$, and therefore
the same is true of $\core(\phi_\lambda)$ and $\core(\phi_\mu)$, by 
lemma~\ref{sections_dense_in_core}. This is the last defining property of a 
universal circle, and completes the sketch of the argument of 
theorem~\ref{universal_circle_exists}.

\section{Constructing invariant laminations}\label{constructing_laminations}

This section contains the first important new results in this paper. Given
a taut foliation $\F$ of an atoroidal $3$--manifold $M$, we construct a pair
of essential laminations $\Lambda^\pm_\split$ of $M$ transverse to $\F$ and
describe their properties.

\subsection{Minimal quotients}

New universal circles can be obtained from old in an uninteresting way:
given a point $p \in S^1_\u$, we can blow up the orbit of $p$ to obtain
a new universal circle $\overline{S^1_\u}$ and a monotone map to $S^1_\u$ whose
gaps are the interiors of the preimages of the points in the orbit of $p$.

These blown up universal circles have the property that there are distinct
points $p, q \in S^1_\u$ whose images are identified under {\em every} map
$\phi_\lambda$. We make the following definition:

\begin{defn}
A universal circle is {\em minimal} if for any distinct $p,q \in S^1_\u$ there
is some $\lambda$ such that $\phi_\lambda(p) \ne \phi_\lambda(q)$.
\end{defn}

In the next lemma, we show that any universal circle which is not minimal
is obtained from a minimal universal circle by blow up.

\begin{lem}\label{minimal_universal_circle}
Let $S^1_\u$ be a universal circle for $\F$. Then there is a minimal universal
circle $S^1_m$ for $\F$ with monotone maps $\phi_\lambda^m:S^1_m \to S^1_\infty(\lambda)$
and a monotone map $m:S^1_\u \to S^1_m$ such that for all $\lambda \in L$
$$\phi_\lambda^m \circ m = \phi_\lambda$$
\end{lem}
\begin{proof}
If $S^1_\u$ is not minimal, define an equivalence relation on $S^1_\u$ by
$p \sim q$ if $\phi_\lambda(p) = \phi_\lambda(q)$ for all $\lambda \in L$. Let
$\gamma^p \subset S^1_\u$ be the interiors of the two closed arcs from two
such distinct $p,q$ with $p \sim q$. Then for each $\lambda \in L$, either
$\gamma^+$ is contained in a single gap of $\phi_\lambda$, or $\gamma^-$ is. Moreover, if 
{\em both} $\gamma^-$ and $\gamma^+$ were contained in gaps of $\phi_\lambda$,
the map $\phi_\lambda$ would be constant, which is absurd.

Now, by lemma~\ref{cores_semicontinuous}, closures of 
gaps of $\phi_\lambda$ vary upper
semicontinuously as a function of $\lambda \in L$. It follows that the subset of
$\lambda \in L$ for which $\gamma^+$ is contained in a gap of $\phi_\lambda$ is
{\em closed}, and similarly for $\gamma^-$. But $L$ is path connected, so either
$\gamma^+$ is contained in a gap of $\phi_\lambda$ for {\em every} $\lambda$, 
or $\gamma^-$ is.

It follows that the equivalence classes of $\sim$ are a $\rho_\u(\pi_1(M))$--equivariant
collection of closed disjoint intervals of $S^1_\u$, and single points, and
therefore the quotient space of $S^1_\u$ by this decomposition defines a new circle
with a $\pi_1(M)$ action induced by the quotient map
$$m:S^1_\u \to S^1_m$$
By construction, for each $\lambda \in L$ the equivalence relation on $S^1_\u$ defined
by $\phi_\lambda$ is coarser than the equivalence relation defined by $m$, and
therefore $\phi_\lambda$ factors through $m$ to a unique map
$\phi_\lambda^m:S^1_m \to S^1_\infty(\lambda)$ satisfying
$$\phi_\lambda = \phi_\lambda^m \circ m$$
\end{proof}

Note that the construction of a universal circle in \S~\ref{constructing_circle_section}
produces a minimal circle.

\subsection{Laminations of $S^1_\u$}

The main purpose of this section is to prove that a minimal universal circle
for a taut foliation with $2$--sided branching
admits a pair of nonempty laminations $\Lambda^\pm_\u$ which are preserved by
the action of $\pi_1(M)$, acting via the representation $\rho_\u$.

\begin{construct}\label{circle_lams_construct}
Let $\lambda \in L$.  Let $L^+(\lambda),L^-(\lambda)$ 
denote the two connected components of $L \backslash \lambda$, where the labelling
is such that $L^+(\lambda)$ consists of the leaves on the positive side of $\lambda$,
and $L^-(\lambda)$ consists of the leaves on the negative side.

Recall that for $X \subset L$, the set $\core(X)$ denotes the union, over
$\lambda \in X$, of the sets $\core(\phi_\lambda)$. As in construction~\ref{boundary_hull}
we can associate to the subset $\core(X)$ the lamination of $\H^2$ which is the
boundary of the convex hull of the closure of $\core(X)$, and thereby construct
the corresponding lamination $\Lambda(\core(X))$ of $S^1$.

Then define
$$\Lambda^+(\lambda) = \Lambda(\core(L^+(\lambda)))$$
and
$$\Lambda^+_\u = \overline{\bigcup_{\lambda \in L} \Lambda^+(\lambda)}$$
and similarly for $\Lambda^-(\lambda)$ and $\Lambda^-_\u$, where the
closure is taken in the space of unordered pairs of distinct points in $S^1_\u$.
\end{construct}

Observe the following property of $\Lambda^+(\lambda)$.

\begin{lem}\label{trivial_map}
Let $\lambda,\mu$ be leaves of $\til{\F}$. Then $\phi_\mu(\Lambda^+(\lambda))$
is trivial unless $\mu < \lambda$.
\end{lem}
\begin{proof}
If $\mu \in L^+(\lambda)$ then by definition, $\core(\mu) \subset \core(L^+(\lambda))$
and therefore every leaf of $\Lambda^+(\lambda)$ is contained in the closure of
a gap of $\mu$. If $\mu \in L^-(\lambda)$ but $\mu$ is incomparable with $\lambda$,
then $\mu$ is incomparable with every element of $L^+(\lambda)$, and therefore
by theorem~\ref{unlinked_extension}, $\core(L^+(\lambda))$ is contained in the closure
of a single gap of $\mu$, and therefore $\phi_\mu(\Lambda^+(\lambda))$ is trivial in
this case too.
\end{proof}

We are now ready to establish the key property of $\Lambda^\pm_\u$: that they
are {\em laminations} of $S^1_\u$.

\begin{thm}\label{union_is_lamination}
Let $\F$ be a taut foliation of an atoroidal $3$--manifold $M$, 
and let $S^1_\u$ be a minimal universal
circle for $\F$. Then $\Lambda^\pm_\u$ are laminations of $S^1_\u$ which are
preserved by the natural action of $\pi_1(M)$. Furthermore, if $L$ branches
in the positive direction, then $\Lambda^+_\u$ is nonempty, and if $L$ branches
in the negative direction, then $\Lambda^-_\u$ is.
\end{thm}
\begin{proof}
We first show that no leaf of $\Lambda^+(\lambda)$ links any leaf of $\Lambda^+(\mu)$,
for $\mu,\lambda \in L$. There are three cases to consider

{\noindent \bf Case (i): $\lambda \in L^-(\mu)$ and $\mu \in L^-(\lambda)$}

In this case, $L^+(\lambda)$ and $L^+(\mu)$ are disjoint, and moreover they
are {\em incomparable}. That is, for every $\nu_1 \in L^+(\lambda)$ and
$\nu_2 \in L^+(\mu)$ the leaves $\nu_1$ and $\nu_2$ are incomparable. It follows
from the definition of a universal circle that for all such pairs, the core
of $\phi_{\nu_1}$ is contained in the closure of a single gap of $\phi_{\nu_2}$,
and {\it vice versa}. Since $L^+(\lambda)$ and $L^+(\mu)$ are path connected,
theorem~\ref{unlinked_extension} implies that $\core(L^+(\mu))$ and
$\core(L^+(\lambda))$ are unlinked. It follows that no leaf of
$\Lambda^+(\lambda)$ links any leaf of $\Lambda^+(\mu)$, as claimed.

{\noindent \bf Case (ii): $\lambda \in L^-(\mu)$ and $\mu \in L^+(\lambda)$}

In this case, we have $L^+(\lambda) \subset L^+(\mu)$ and therefore
$$\core(L^+(\lambda)) \subset \core(L^+(\mu))$$
so the claim is proved in this case too.

{\noindent \bf Case (iii): $\lambda \in L^+(\mu)$ and $\mu \in L^+(\lambda)$}

In this case, observe that $L^-(\lambda) \subset L^+(\mu)$ and
$L^-(\mu) \subset L^+(\lambda)$, and therefore
$$L = L^+(\mu) \cup L^+(\lambda)$$
Since $S^1_\u$ is minimal, every point in $S^1_\u$ is a limit of
a sequence of points in $\core(\phi_{\lambda_i})$ for some sequence
$\lambda_i$. It follows that $\core(L)$ is all of $S^1_\u$, and
therefore $\core(L^+(\lambda)) \cup \core(L^+(\mu)) = S^1_\u$.

Now, if two subsets $X, Y \subset S^1$ satisfy 
$\overline{X} \cup \overline{Y} = S^1$, then
the boundaries of the convex hulls of $\overline{X}$ and $\overline{Y}$ do
not cross in $\H^2$. For, if $l,m$ are boundary geodesics of $H(\overline{X})$
and $H(\overline{Y})$ respectively which cross in $\H^2$, 
then $l,m$ both bound open half spaces 
$l^+,m^+$ which are disjoint from $H(\overline{X})$ and $H(\overline{Y})$
respectively. Moreover, since $l,m$ are transverse, the intersection
$l^+ \cap m^+$ contains an open sector in $\H^2$, 
which limits to some nonempty interval in
$S^1$ which by construction is disjoint from both $X$ and $Y$. But this
contradicts the defining property of the pair $X,Y$. 
This contradiction proves the claim in this case too.

It remains to show that $\Lambda^+_\u$ is nonempty when $L$ branches in the
positive direction. Now, for any $\lambda \in L$, $\core(\phi_\lambda)$
is perfect by lemma~\ref{perfect_core}. It suffices
to show $\core(L^+(\lambda))$ is not equal to $S^1_\u$.

If we can find another leaf $\mu$ with $\lambda \in L^-(\mu)$ and $\mu \in
L^-(\lambda)$, then as above, $\core(L^+(\lambda))$ and $\core(L^+(\mu))$ are unlinked
as subsets. It follows that the subset $\core(L^+(\lambda))$ is contained in the closure of
a single interval in the complement of $\core(L^+(\mu))$ 
and conversely, and therefore neither core is dense. 
To see that such a $\mu$ exists, note that if there is $\nu$ with
$\nu < \mu$ and $\nu < \lambda$ but $\mu,\lambda$ incomparable, then $\mu$ will
have the desired properties.

Since $L$ branches in the positive direction, there is $\nu$ and some leaves
$\lambda',\mu$ with $\nu < \mu,\lambda'$ and $\lambda',\mu$ incomparable.
Since $\F$ is taut, if $\pi(\lambda')$ and $\pi(\lambda)$ denote the projections
of $\lambda,\lambda'$ to $M$, there is some transverse positively oriented 
arc $\gamma$ from $\pi(\lambda')$ to $\pi(\lambda)$. Lifting to $\til{M}$, 
we see there is some $\alpha \in \pi_1(M)$ such that 
$\alpha(\lambda') < \lambda$. Then $\alpha(\mu)$ is the desired leaf.

The corresponding properties for $\Lambda^-_\u$ are proved by reversing the
orientation on $L$.
\end{proof}

\subsection{Branched surfaces and branched laminations}

\begin{construct}\label{construct_essential_lamination}
Let $\Lambda^\pm_\u$ be the invariant laminations of $S^1_\u$ provided
by theorem~\ref{union_is_lamination}. For each $\lambda \in L$,
there are lamination of $S^1_\infty(\lambda)$ given by the pushforward
$\phi_\lambda(\Lambda^\pm_\u)$. By construction~\ref{circle_geodesic}, these
laminations of $S^1_\infty(\lambda)$ span geodesic laminations of $\lambda$, which
we denote by $\Lambda^\pm_\geo(\lambda)$. Then define
$$\til{\Lambda}^\pm_\geo = \bigcup_{\lambda \in L} \Lambda^\pm_\geo(\lambda)$$
\end{construct}

Note the tilde notation to be consistent with the convention that
$\til{\Lambda}^\pm_\geo$ covers an object in $M$. The objects $\til{\Lambda}^\pm_\geo$
are not yet necessarily $2$--dimensional laminations; rather they are
{\em branched laminations}, to be defined shortly.
On the other hand, they have the important property that
the branch locus of each leaf
is a {\em $1$--manifold} (that is, there are no double points of the branch locus)
and moreover, the sheets come with a parameterization by leaves of $\Lambda^\pm_\u$
that lets us split them open in a canonical way to a lamination.

The definition we give here of a branched lamination is not the most general possible,
since for us, every branched lamination comes together with an ordinary lamination which it
fully carries. Branched laminations are a generalization of {\em branched surfaces};
see \cite{uO86} for a definition and basic properties of branched surfaces.

\begin{defn}\label{branched_lamination_definition}
A {\em branched lamination fully carrying a lamination} 
$K \subset M$ is given by the following data:
\begin{enumerate}
\item{An open submanifold $N \subset M$}
\item{A $1$--dimensional foliation $X_V$ of $N$}
\item{A lamination $\Lambda$ of $N$ transverse to $X_V$, intersecting every leaf of $X_V$}
\item{A surjective map $\psi:N \to N$ from $N$ to itself which is {\em monotone}
on each leaf $x$ of $X_V$}
\end{enumerate}
The underlying space of the
branched lamination itself is the image $K = \psi(\Lambda)$, 
thought of as a subset of $M$. 
We say that the lamination $\Lambda$ is {\em fully carried by} 
$K$, and is obtained by {\em splitting $K$ open}.
\end{defn}

Notice that with this definition, we allow the possibility that $K = N = M$, which
would happen for instance if $\Lambda$ is a foliation.

Let us describe our strategy to realize $\til{\Lambda}^\pm_\geo$ as branched
laminations, which fully carry split open laminations $\til{\Lambda}^\pm_\split$.
 
Firstly, observe that we can define in generality a branched lamination as a structure on
$M$ which is {\em locally} modelled on the structure in 
definition~\ref{branched_lamination_definition}, and for which the
$1$--dimensional foliations $X_V$ in local charts are required to piece together to give
a global transverse $1$--dimensional foliation, but for which the laminations
$\Lambda$ and the map $\psi$ are only defined locally, with no conditions on
how they might piece together globally. General branched laminations do not
always fully carry laminations.

Another way of thinking of a branched lamination is as the total space of a
distribution defined on a closed subset of $M$ which is integrable, but
not uniquely. That is, through every point, there is a complete integral submanifold
tangent to the distribution, but such submanifolds might not be disjoint.
The {\em branch locus} of the branched lamination consists of the union of
the boundaries of the subsets where such distinct integral submanifolds agree.
In particular, the branch locus has the structure of a union of $1$--manifolds.
In the case of a branched surface, this branch locus is a finite union of circles,
and one typically requires this union of circles to be in general position with
respect to each other.

Given a branched lamination $K$, one can always find an {\em abstract} 
lamination ``carried'' by $K$ which consists of the disjoint union
of some collection of such integral submanifolds, topologized leafwise with the
path topology, and as a lamination by the compact open topology. The difficulty
is in embedding this abstract lamination in $N$ transverse to the foliation $X_V$.
This amounts to finding a local {\em order structure} on the leaf space of this
abstract lamination. Once this order structure is obtained, the process of
recovering $\psi$ from $K$ is more or less the same as the usual process of
blowing up some collection of leaves of a foliation or lamination, as described
in \cite{CCI}. In our case, the leaves of the abstract laminations
carried by $\til{\Lambda}^\pm_\geo$ are the unions $\bigcup_{\lambda \in L} \phi_\lambda(l)$
where $l$ is a leaf of $\Lambda^\pm_\u$, and we are implicitly thinking of
$\phi_\lambda(l)$ as a geodesic in $\lambda$ by construction~\ref{circle_geodesic}. 
The desired local order structure on the leaves of these abstract laminations
comes from the local order structure on the order trees which are the leaf spaces
of the geodesic laminations of $\H^2$ constructed from $\Lambda^\pm_\u$. 
In this way, the abstract laminations may be realized as laminations in $\til{M}$
fully carried by $\til{\Lambda}^\pm_\geo$. This is the summary of our strategy. Now
we go into detail.

To establish the desired properties of $\til{\Lambda}^\pm_\geo$, we must first understand
how the laminations $\Lambda^\pm_\geo(\lambda)$ vary as a function of $\lambda$.

Let $\tau$ be a transversal to $\til{\F}$. The cylinder $UT\til{\F}|_\tau$,
thought of as a circle bundle over $\tau$, carries two natural families of
sections. The first family of sections comes from the structure maps
$e$ and $\phi$.

\begin{construct}
Let $\tau$ be a transversal to $\til{\F}$. The endpoint map defines an
embedding $e:UT\til{\F}|_\tau \to E_\infty$. The structure map of the
universal circle $\phi:S^1_\u \times L \to E_\infty$, composed with $e^{-1}$,
defines a canonical collection of sections of the circle bundle 
$$UT\til{\F}|_\tau \to \tau$$
as follows. If we let $\iota:\tau \to L$ denote the embedding induced by
the quotient map $\til{M} \to L$, then the arcs $p \times \iota(\tau)$
with $p \in S^1_\u$ map to a family of arcs in $E_\infty|_{\iota(\tau)}$.
In the case of the universal circles constructed in \S~\ref{constructing_circle_section}
these are the restriction of the special sections to $\iota(\tau)$.
Then $e^{-1}$ pulls these back to define a family of sections 
of $UT\til{\F}|_\tau$, which by abuse of notation we call the
{\em special sections} over $\tau$. If $p \in S^1_\u$, we denote by
$\sigma(p)|_\tau$ the special section corresponding to $p$ over $\tau$.
\end{construct}

The second family of sections comes from the geometry of $\til{M}$.

\begin{construct}
A Riemannian metric on $M$ pulls back to a Riemannian metric on $\til{M}$.
Parallel transport with respect to the Levi--Civita connection does not preserve the
$2$--dimensional distribution $T\til{\F}$, but the combination of the
Levi--Civita connection of the metric on $\til{M}$ together with orthogonal
projection to $T\til{\F}$ defines an orthogonal (i.e. metric preserving)
connection on $T\til{\F}$.

If $\tau$ is a transversal, this connection defines a trivialization of
$UT\til{\F}|_\tau$ by parallel transport along $\tau$. We call the fibers 
of this trivialization the {\em geometric sections} over $\tau$.
\end{construct}

Let $\nu$ denote the unit normal vector field to $\F$, and $\til{\nu}$ the
unit normal vector field to $\til{\F}$.

\begin{lem}\label{uniform_comparison_of_connections}
There is a uniform modulus $f:\R^+ \to \R^+$ with
$$\lim_{t \to 0} f(t) = 0$$
such that for any $p \in S^1_\u$, any $q \in \til{M}$ and $\tau(t)$ any integral
curve of $\til{\nu}$ through $q$ parameterized by arclength, then if 
$$r = \sigma(p)|_{\tau(0)} \in UT_q\til{\F}$$
and $\sigma'(\cdot)$ denotes the geometric section over $\tau$ obtained by
parallel transporting $r$, we have
$$\text{arcwise distance from }\sigma(p)|_{\tau(t)} \text{ to } \sigma'(t) \text{ in }
UT_{\tau(t)}\til{\F} \le f(t)$$
\end{lem}
\begin{proof}
This just follows from the compactness of $UT\F$ and the continuity of $e$ and
$\phi$.
\end{proof}

Said another way, lemma~\ref{uniform_comparison_of_connections} says that a
geometric section and a special section which agree at some point cannot
move apart from each other too quickly. Since the geometric sections are defined by
an orthogonal connection, it follows that if $\sigma(p)$ and
$\sigma(q)$ are two special sections, the angle between them cannot vary too quickly.
This lets us prove the following.

\begin{lem}\label{continuity_of_image_lamination}
The laminations $\Lambda^\pm_\geo(\lambda)$ vary continuously on compact
subsets of $\til{M}$, as a function of $\lambda \in L$. Moreover, the sets
$\til{\Lambda}^\pm_\geo$ are closed as subsets of $\til{M}$.
\end{lem}
\begin{proof}
The continuity of $\Lambda^\pm_\geo(\lambda)$ on compact subsets of $\til{M}$
follows from the fact that the leaves $\lambda$ themselves vary
continuously on compact subsets, together with the continuity of $e$ and $\phi$.

Now we show that the unions $\til{\Lambda}^\pm_\geo$ are closed.
Let $\lambda_i \to \lambda$ and $p_i \in \lambda_i \to p \in \lambda$ be a
sequence of leaves of $\til{\F}$ and points in those leaves. Since
$p_i \to p$, it follows that for sufficiently large $i$, the leaves
$\lambda_i$ are all comparable, and contained in an interval $I \subset L$,
so without loss of generality, we can assume that all $\lambda_i$ are contained
in $I$. Let $\tau$ be an orthogonal trajectory to $\til{\F}$ through $p$,
parameterized by arclength, and let $q_i \in \lambda_i$ be equal to 
$\tau \cap \lambda_i$. Suppose that $p_i \in \Lambda^+_\geo(\lambda_i)$ for
each $i$. We must show that $p \in \Lambda^+_\geo(\lambda)$.
Now, since $p_i \in \Lambda^+_\geo(\lambda_i)$, there is a leaf $l_i$ of
$\Lambda^+_\geo(\lambda_i)$ with $p_i \in l_i$. Geometrically, $l_i$ is
just a geodesic in $\lambda$ with respect to its hyperbolic metric.
Let $l_i^\pm \in S^1_\u$ be a pair of points which span a leaf $l_i'$
of $\Lambda^+_\u$
which maps to $l_i$ under $\phi_\lambda$ by the pushforward construction and
construction~\ref{circle_geodesic}. Then by lemma~\ref{uniform_comparison_of_connections},
the angle between the special sections over $\tau$ defined by $l_i^\pm$ cannot
vary too quickly. But in $UT_{p_i}\lambda_i$, the angle between the endpoints of
$l_i$ is $\pi$, since $p_i$ lies on the geodesic $l_i$.
It follows that as $i \to \infty$, the angle between the special
sections over $\tau$ defined by $l_i^\pm$ converges to $\pi$, and therefore the
pushforward leaves $\phi(l_i')$ span geodesics in $\Lambda^+_\geo(\lambda)$ which contain
points converging to $p$. Since $\Lambda^+_\geo(\lambda)$ is closed in $\lambda$,
the point $p \in \Lambda^+_\geo(\lambda)$, as claimed.
\end{proof}

The next lemma shows, as promised, that $\til{\Lambda}_\geo^\pm$ are
branched laminations which can be split open. The following lemma is somewhat 
{\it ad hoc}. However, the basic idea is very simple, and is
precisely as described in the paragraphs
following definition~\ref{branched_lamination_definition}. Namely, the branched
laminations $\til{\Lambda}^\pm_\geo$ can be split open because they are
parameterized by abstract laminations whose
leaf spaces already have well--defined local order structures.

\begin{lem}\label{image_is_branched_lamination}
$\til{\Lambda}_\geo^\pm$ are branched laminations of $\til{M}$, fully carrying
laminations $\til{\Lambda}^\pm_\split$ which are preserved by the action of $\pi_1(M)$.
\end{lem}
\begin{proof}
For the sake of notation, we restrict to $\til{\Lambda}_\geo^+$.

Fix some small $\epsilon$, and for each leaf $\lambda$ of $L$ let $N(\lambda)$ be
the subset of points in $\lambda$ which are distance $< \epsilon$ from 
$\Lambda^+_\geo(\lambda)$, and let
$$\til{N} = \bigcup_\lambda N(\lambda)$$
The nearest point map (in the path metric on $\lambda$) defines a retraction from
$N(\lambda)$ to $\Lambda^+_\geo(\lambda)$, away from the set of points which are
equally close to two leaves; call these {\em ambivalent} points. 
The preimages of this retraction, together with the
points equally close to two leaves, give a 
$1$--dimensional foliation of $N(\lambda) \backslash \Lambda^+_\geo(\lambda)$
by open intervals, with at most one ambivalent point on each open interval, as the
midpoint. 

If $\Lambda^+_\geo$ foliates some region,
then the integral curves of the orthogonal distribution define a foliation of
the foliated region of $\Lambda^+_\geo$. Together, this defines a $1$--dimensional
foliation of $N(\lambda)$. By lemma~\ref{continuity_of_image_lamination},
these foliations vary continuously from leaf to leaf of $\lambda$, and
define a $1$--dimensional foliation $X_V$ of $\til{N}$, which is an open neighborhood
of $\til{\Lambda}_\geo^+$.

If $l,m$ are leaves of $\Lambda^+_\geo(\lambda)$ which are both in the closure of
the same complementary region, and which contain points which are $<2\epsilon$ apart,
then there are at most two points $p,q$ in this complementary region which are
distance exactly $\epsilon$ in $\lambda$ to both $l$ and $m$. Call such points
{\em cusps}. The set of cusps in each leaf $\lambda$ of $\til{\F}$ are isolated;
futhermore, by lemma~\ref{continuity_of_image_lamination} the set of cusps in
$\lambda$ varies continuously as a function of $\lambda$, thereby justifying the
notation $p(\lambda)$ for a family of leafwise cusps, with the possibility of
birth--death pairs in the sense of Morse theory when two distinct cusp points
$p(\lambda),q(\lambda)$ coalesce at some leaf $\lambda_0$ and disappear for
nearby leaves on one side. It follows that the
union of all cusps defines a locally finite collection $\til{\fc}$ of properly
embedded lines in $\til{M}$ which covers a link $\fc \subset M$.
By abuse of notation, we call $\til{\fc}$ the {\em cusps} of $N$.
Observe that the cusps parameterize the branching of the leaf space of $X_V$, as follows.
For each point $p \in \til{\fc}$ there is a $1$--parameter family
$\gamma_t$ of leaves of $X_V$, with
$t \in [0,1)$, such that the limit of the
as $t \to 1$ is a union of two leaves $\gamma_1^\pm$
together with the point $p$, which is in the closure of both $\gamma_1^+$ and $\gamma_1^-$.
We refer to such a family of leaves of $X_V$ as a {\em bifurcating family}.

To show that $\til{\Lambda}^+_\geo$ is a branched lamination fully carrying
a lamination, we must first define a map $\psi:\til{N} \to \til{N}$ 
which is monotone on each leaf
of $X_V$. For convenience, we use construction~\ref{circle_geodesic} to think
of $\Lambda^+_\u$ as a geodesic lamination of a copy $\H^2_\u$ of the
hyperbolic plane bounded by $S^1_\u$. Notice that each leaf $\gamma$ 
of $X_V$ is contained in a leaf $\lambda$ of $\til{\F}$. The leaf $\gamma$ might be
bounded or unbounded in $\lambda$, the latter case occurring for instance if
$\Lambda^+_\geo(\lambda)$ is a foliation. A bounded endpoint of
$\gamma$ determines a complementary region to $\Lambda^+_\u$ in $\H^2_\u$.
Pick a point in such a complementary region.
An unbounded end determines an endpoint in $S^1_\infty(\lambda)$,
which determines its preimage under $\phi_\lambda^{-1}$ in $S^1_\u$. This
preimage might be a point or an interval; for concreteness, if it is an interval,
pick its anticlockwisemost point. Span the two points constructed in this way
by a geodesic $\gamma_\u$. This geodesic $\gamma_\u$ can be thought of as a
``preimage'' to $\gamma$. Note that by our choice of ideal 
endpoints for $\gamma_\u$ that $\gamma_\u$
does not cross any leaves of $\Lambda^+_\u$ whose endpoints are identified by
$\phi_\lambda$. It follows that $\gamma_\u$ crosses exactly those leaves of
$\Lambda^+_\u$ which correspond to leaves of $\Lambda^+(\lambda)$ crossed by
$\gamma$. We define a monotone map $\psi: \gamma_\u \to \gamma$ which takes
each intersection $\gamma_\u \cap \Lambda^+_\u$ to the corresponding
intersection $\gamma \cap \Lambda^+_\geo(\lambda)$, and takes complementary intervals
either to the corresponding intervals, or collapses them to points if the
corresponding leaves in $\Lambda^+_\u$ are identified in $\Lambda^+_\geo(\lambda)$.

We want to make the assignment $\gamma \to \gamma_\u$ continuously as a function of
$\gamma$, at least away from the cusps $\til{\fc}$. 
This amounts to choosing the endpoints of $\gamma_\u$ in complementary
regions to $\Lambda^+_\u$ in $\H^2_\u$ continuously as a function of $\gamma$. Since
the complementary regions are all homeomorphic to disks, and are therefore
contractible, there is no obstruction to making such a choice.
It is clear that this construction can be done in a $\pi_1(M)$ equivariant manner, where
we think of $\pi_1(M)$ acting on the leaves of $\Lambda^+_\u$ and
permuting the complementary regions as {\em sets}.
Along the cusps $\til{\fc}$, one must be slightly more careful.
If $\gamma_t$ with $t \in [0,1)$ limiting to $\gamma_1^\pm$
is a bifurcating family, we must choose $(\gamma_t)_\u$ and
$(\gamma_1^\pm)_\u$ so that there is an equality
$$(\gamma_1^-)_\u \cup p \cup (\gamma_1^+)_\u = \lim_{t \to 1} (\gamma_t)_\u$$
for some $p$ in a complementary region to $\H^2_\u$. Again, the contractibility of
complementary regions implies that this can be done, even equivariantly.

For each $\gamma$, the graph of $\psi:\gamma_\u \to \gamma$ defines an interval
$\psi(\gamma_\u)$ in the product $\H^2_\u \times \til{N}$.
The disjoint union of intervals $\psi(\gamma_\u)$ as $\gamma$ varies over leaves of $X_V$
is itself an open $3$--manifold $\til{N'}$ homeomorphic to $\til{N}$ as a 
subspace of $\H^2_\u \times \til{N}$. 
Moreover, the intersections of the geodesics $\gamma_\u$ 
with leaves of $\Lambda^+_\u$ defines
a lamination $\til{\Lambda}^+_\split$
of $\til{N'}$ that maps by $\psi$ to $\til{\Lambda}^+_\geo$.

The action of $\pi_1(M)$ on the base $3$--manifold $N$ induces an action on
$\til{N'}$ as follows. Since we want the actions on $\til{N'}$ and $\til{N}$ to
be semiconjugate under the monotone map $\psi$, we must just decide how an
element $\alpha \in \pi_1(M)$ should act on point preimages of $p \in \til{N}$.
Now, for each $p \in \til{N}$, either $\psi^{-1}(p)$ is a point, or an
interval in a complementary region of
$\Lambda^+_\u$ with endpoints on distinct leaves $l,m$ of $\Lambda^+_\u$ which map to
the same leaf of some $\Lambda^+_\geo(\lambda)$.
But then for $\alpha \in \pi_1(M)$,
the preimage of $\psi^{-1}(\alpha(p))$ is also a complementary interval,
with endpoints on leaves $\alpha(l),\alpha(m)$ of $\Lambda^+_\u$.
The interval we
define $\alpha:\psi^{-1}(p) \to \psi^{-1}(\alpha(p))$ to be the
unique affine homeomorphism which takes the endpoint on $l$ to the endpoint
on $\alpha(l)$, and the endpoint on $m$ to the endpoint on $m$. Here we mean
affine with respect to the length induced as a geodesic segment in $\H^2$.
This is the desired action.
\end{proof}

Notice that we can choose $\psi:N \to N$ to have point
preimages which are as small as desired. It follows that the laminations
$\til{\Lambda}^\pm_\split$ can be chosen to intersect leaves of $\til{\F}$ in lines
which are uniformly $(k,\epsilon)$ quasigeodesic, 
for any choice of $k>1,\epsilon > 0$.

Define $\Lambda^\pm_\split$ to be the laminations of $M$ covered by
$\til{\Lambda}^\pm_\split$.

Notice too that if $\Lambda^+_\split$ for instance is a genuine lamination, there
is a choice of partition into guts and interstitial regions for which the cores
of the interstitial annuli are exactly the cusps $\fc$.

\begin{thm}\label{laminations_essential}
Let $\F$ be a taut foliation of an atoroidal $3$--manifold $M$. Suppose $\F$
has two--sided branching. Then $M$
admits laminations $\Lambda^\pm_\split$ which are
are essential laminations of $M$, which are transverse to $\F$, and which
intersect the leaves of $\F$ in curves which are uniformly $(k,\epsilon)$
quasigeodesic, for any $k>1,\epsilon >0$.
\end{thm}
\begin{proof}
We construct $\til{\Lambda}^\pm_\split$ as in lemma~\ref{image_is_branched_lamination},
covering laminations $\Lambda^\pm_\split$ in $M$.

By construction, the leaves of $\til{\Lambda}^\pm_\split$ are all planes,
so $\Lambda^\pm_\split$ do not contain any spherical leaves or torus
leaves bounding a solid torus, and complementary regions admit no compressing disks. 
Moreover, since $M$ admits a taut foliation $\F$
by hypothesis, $\til{M}$ is homeomorphic to $\R^3$, so complementary regions
admit no essential spheres. It remains to show that there are no compressing 
monogons.

If $D$ is a compressing monogon for $\Lambda^+_\split$, 
there are points $p,q$ in $\partial D$ contained in
a leaf $\lambda$ of $\Lambda^+_\split$ which are
arbitrarily close together in $D$ but arbitrarily far apart in $\lambda$.
Lift $D, p, q, \lambda$ to $\til{M}$, where by abuse of notation we refer to
them by the same names. Since $p,q$ are arbitrarily close in $\til{M}$,
they are contained in comparable leaves $\mu_1,\mu_2$ of $\til{\F}$.
Suppose $p \in \mu_1 \cap \lambda$. Let $\tau$ be a short orthogonal
trajectory from $\mu_1$ to $\mu_2$. The endpoints of the quasigeodesic
$\lambda \cap \mu_1$ determine a leaf of $\Lambda^+_\u$, which determines
a pair of special sections of $UT_\tau\til{\F}$.
By lemma~\ref{uniform_comparison_of_connections} and the uniformity of $k,\epsilon$, 
the angle between these special
sections stays close to $\pi$ along $\tau$, for $\tau$ sufficiently short. It follows
that there is a short path in $\lambda$, starting from $p$,
from $\mu_1$ to some $p' \in \mu_2$. But $\til{\Lambda}^+_\split \cap \mu_2$
is a $(k,\epsilon)$--quasigeodesic lamination, so $p'$ and $q$ can be joined by
a short path in $\lambda \cap \mu_2$, and therefore $p$ and $q$ are close in
$\lambda$, contrary to the definition of $D$.

It follows that no such compressing monogon $D$ exists, and the laminations
$\Lambda^\pm_\split$ are essential, as claimed.
\end{proof}

\subsection{Straightening interstitial annuli}

In this subsection we show that each complementary region to $\Lambda^\pm_\split$
can be exhausted by a sequence of guts, for some partition 
into guts and $I$--bundles, such that the
interstitial annuli are {\em transverse} to $\F$. This implies that complementary
regions are {\em solid tori}. Note that this does not address the question, left
implicit in the last subsection, of whether or not the 
laminations $\Lambda^\pm_\split$ are {\em genuine}; but it does show that
{\em if} they are genuine, then they are very full.

Each leaf of $\til{\Lambda}^\pm_\split$ is transverse to the foliation
$\til{\F}$, and therefore it inherits a codimension one foliation, whose leaves
are the intersection with leaves of $\til{\F}$. We show that this foliation
branches in at most one direction.

\begin{lem}\label{one_sided_branching}
Let $\pi$ be a leaf of $\til{\Lambda}^+_\split$. The induced foliation
$\pi \cap \F$ of $\pi$ does not branch in the positive direction, and similarly
for leaves of $\til{\Lambda}^-_\split$.
\end{lem}
\begin{proof}
Let $l$ be a leaf of $\Lambda^+(\lambda)$. That is, a leaf of $\Lambda(\core(L^+(\lambda)))$,
thought of as an unordered pair of distinct points
in $S^1_\u$. By lemma~\ref{trivial_map}, the image
$\phi_\mu(l)$ is trivial unless $\mu < \lambda$. 

The subset of $L$ consisting of leaves $\mu$ with $\mu < \lambda$ does not branch
in the positive direction. Consider the union
$$\Pi(l) = \bigcup_{\lambda \in L} \phi_\lambda(l) \subset \til{M}$$
where for each $\lambda$, we think of $\phi_\lambda(l)$ as a leaf of 
$\Lambda^+_\geo(\lambda)$. Let $H \subset L$ be the subset of leaves $\lambda$ with
$\phi_\lambda(l)$ nontrivial. Then $H$ does not branch in the positive direction.
Moreover, $\Pi(l)$ is carried by the branched lamination
$\til{\Lambda}^+_\geo$, and naturally embeds into the split open lamination
$\til{\Lambda}^+_\split$ as a union of leaves $\pi_1,\pi_2,\cdots$,
corresponding to the connected components $H_1,H_2, \cdots$ of $H$.
Moreover, for each leaf $\lambda \in H_i$, the intersection 
$\pi_i \cap \lambda = \phi_\lambda(l)$ is a single line. It follows that the induced
foliation of each $\pi_i$ does not branch in the positive direction.

Now, the leaves $l$ of laminations $\Lambda^+(\lambda)$ with $\lambda$ in $\til{\F}$ are
dense in $\Lambda^+_\u$. If $\pi$ is a limit of leaves $\pi_i$ where the induced
foliation of $\pi_i$ does not branch in the positive direction, the same is true for
$\pi$. To see why this is true, let $J$ be the subset of $L$ which $\pi$ intersects.
Lemma~\ref{continuity_of_image_lamination} implies that the set of leaves of
$\til{\F}$ which $\pi$ intersects in a single component is both open and closed
in $J$, and is therefore equal to $J$. It follows that if $\pi$ branches in the
positive direction, then $J$ branches in the positive direction. In this case,
$\pi$ intersects leaves $\mu_1,\mu_2$ of $\til{\F}$ which are incomparable
but satisfy $\mu_1> \lambda,\mu_2 > \lambda$ for some third leaf $\lambda$ of $\til{\F}$.
But this means that $\pi_i$ intersects both 
$\mu_1$ and $\mu_2$ for sufficiently big $i$, contrary
to the fact that $H_i$ does not branch in the positive direction.

This contradiction proves the claim, and the lemma follows.
\end{proof}

\begin{lem}\label{horizontally_foliated}
Let $\Lambda^\pm_\split$ be the laminations constructed in 
theorem~\ref{laminations_essential}. Then there is a system
of interstitial annuli $A_i^\pm$ for $\Lambda^\pm_\split$, such that,
(supressing the superscript $\pm$ for the moment) each
$A_i$ satisfies the following properties:
\begin{enumerate}
\item{The intersection of $A_i$ with the foliation $\F$ induces a nonsingular
product foliation of $A_i = S^1 \times I$ by intervals $\text{point} \times I$.}
\item{There is a uniform $\epsilon$, which may be chosen as small as desired, such
that each leaf of the induced foliation of each $A_i$ has length $\le \epsilon$.
Moreover, every point $p$ in an interstitial region can be connected to
a point in the lamination by an arc contained in a leaf of 
$\F$ of length $\le \epsilon/2$.}
\end{enumerate}
We say that such an interstitial system is {\em horizontally foliated}.
\end{lem}
\begin{proof}
We do the construction upstairs in $\til{M}$. For convenience, we concentrate
on $\til{\Lambda}^+_\split$. By abuse of notation, we denote
$\til{\Lambda}^+_\split \cap \lambda$ by $\Lambda^+_\split(\lambda)$,
for $\lambda$ a leaf of $\til{\F}$. We suppose that we have performed the splitting
in such a way that the geodesic curvatures of the leaves of $\Lambda^+_\split(\lambda)$
are uniformly pinched as close to $0$ as we like.

Recall that we can split open $\til{\Lambda}^+_\geo$ so that 
the laminations $\Lambda^+_\split(\lambda)$ for $\lambda$ a leaf of
$\til{\F}$ are as close as desired to geodesic laminations.
We define the interstitial regions to be
precisely the set of points $p$ in each leaf $\lambda$ of $\til{\F}$, not in
$\Lambda^+_\split(\lambda)$, and which are contained in an arc in $\lambda$ of
length $\le \epsilon$ between two distinct boundary leaves. This obviously
satisfies the desired properties.
\end{proof}

\begin{lem}\label{straighten_interstitial_region}
Let $A_i$ be a horizontally foliated system of interstitial annuli for 
$\Lambda^+_\split$. Then, after possibly throwing away
annuli bounding compact interstitial regions, the system $A_i$ can be isotoped so that at
the end of the isotopy, each annulus is transverse or tangent to $\F$.
\end{lem}
\begin{proof}
We assume before we start that we have thrown away annuli bounding compact interstitial
regions. The key to the proof of this lemma is the fact that the foliation of leaves
of $\til{\Lambda}^+_\split$ by the intersection with $\til{\F}$ does
not branch in the positive direction. This lets us inductively push local minima
in the positive direction, until they cancel local maxima. The first step is
to describe a homotopy from each $A_i$ to some new $A_i$ which is either transverse or
tangent to $\F$. At each stage of this homotopy, we require that the image of $A_i$
be foliated by arcs of its intersection with leaves of $\F$. It is clear that
there is no obstruction to doing this. We fix notation: let $C$ be a complementary
region, and $C_i$ the interstitial $I$--bundle bounded by the $A_i$.

Let $\fc_i$ be the core of an interstitial annulus $A_i$. Suppose $\fc_i$ is
not transverse or tangent to $\F$. Then $\fc_i$ must have at least one local
maximum and one local minimum, with respect to the foliation $\F$. Either
$A_i$ bounds a compact $I$ bundle over a disk, or else the universal
cover $\til{A_i}$ is noncompact, and $\til{\fc_i}$ has infinitely many local
maxima and minima. By hypothesis, we have already thrown away compact $I$--bundles,
so we may assume $\til{A_i}$ is noncompact.
Let $p$ be a local minimum on $\til{\fc_i} \cap \lambda$, 
and $p^\pm$ neighboring local maxima on $\til{\fc_i} \cap \lambda^\pm$ for 
leaves $\lambda,\lambda^\pm$ of $\til{\F}$. Then by construction, 
$\lambda < \lambda^+,\lambda^-$ in the partial order on $L$. (It should be
remarked that the $\pm$ notation
reflects the order of the points $p^-,p,p^+$ in the arc $\til{\fc_i}$, and not the
order structure of the related leaves $\lambda^-,\lambda,\lambda^+$ in $L$.)

By lemma~\ref{one_sided_branching}, $\lambda^+$ and $\lambda^-$
are comparable; without loss of generality, we can assume $\lambda^- \le \lambda^+$.
Then there is $q$ on $\til{\fc_i}$ between $p$ and $p^+$ with
$q \in \til{\fc_i} \cap \lambda^-$. The points $q$ and $p^-$ are contained in
arcs $I_q,I_{p^-}$ of $\til{A_i}$ which bound a rectangle $R \subset \til{A_i}$. 
The arcs $I_q,I_{p^-}$ also bound a rectangle $R' \subset \lambda^-$ of a
complementary region to $\Lambda^+(\lambda^-)_\split$. The union $R \cup R'$ is
a cylinder which bounds an interval bundle over a disk $D \times I$ in a
complementary region. We push $R$ across this $D \times I$ to $R'$, and then slightly
in the positive direction, cancelling the local minimum at $p$ with the local maximum at
$p'$. Do this equivariantly with respect to the action of
$\pi_1(A_i)$ on the lift $\til{A_i}$. 
After finitely many moves of this kind, all maxima and minima are cancelled, and
we have produced new immersed annuli $A_i'$ either transverse or tangential to $\F$,
and homotopic to the original $A_i$. If $A_i'$ is tangent to $\F$, it finitely covers
some annular complementary region to $\Lambda^+_\split \cap \lambda$ for some leaf
$\lambda$ of $\F$. Since it is homotopic to an embedded annulus, by
elementary $3$--manifold topology the degree of this
covering map must be one, and $A_i'$ must be embedded. See for example
chapter 13 of \cite{Hempel}.

If $A_i'$ is transverse to $\F$, it is either embedded, or cuts off finitely many
$\text{bigon} \times I$ where the edges of the bigons are transverse to $\F$. 
By inductively pushing arcs across innermost embedded bigons, we can reduce the number of
bigons by two at a time. We can do this unless there is a single arc of self--intersection
of $A_i'$ which corresponds to both cusps of a (non--embedded) bigon. But, since
$A_i'$ is homotopic to $A_i$ which is embedded, the number of arcs of intersection must
be even, for homological reasons. It follows that all bigons can be cancelled, and $A_i'$
is homotopic to $A_i''$ which is embedded. By further cancelling bigons of intersection
of $A_i''$ with $A_j''$ for distinct indices $i,j$ we can assume the union of the
$A_i''$ are disjointly embedded. Let $C_i''$ be the $I$--bundle bounded by the $A_i''$.
By construction, $C_i$ and $C_i''$ are homotopy equivalent in $C$. Since
$C_i,C_i''$ and $C$ are all Haken, again, by standard $3$--manifold topology, 
$C_i$ and $C_i''$ are isotopic in $C$, and therefore the system $A_i$ is isotopic to the
system $A_i''$.

Compare lemma 2.2.2 of \cite{dC00}.
\end{proof}

\begin{thm}\label{complementary_regions_polygon_bundles}
Every complementary region to $\Lambda^\pm_\split$ is a finite--sided ideal
polygon bundle over $S^1$. Moreover, after possibly removing finitely many
isolated leaves and collapsing bigon bundles over $S^1$, 
the laminations $\Lambda^\pm_\split$ are {\em minimal}.
\end{thm}
\begin{proof}
By lemma~\ref{horizontally_foliated} and lemma~\ref{straighten_interstitial_region},
we can exhaust each complementary region by a sequence of guts $\fG_i$ bounded
by interstitial annuli transverse to $\F$. It follows that the boundary of each
$\fG_i$ is foliated by the intersection with $\F$, and therefore each boundary
component is a torus. Since $M$ is irreducible and
atoroidal, these tori either bound solid
tori, or are contained in balls. But by construction, the core of each essential
annulus is transverse to $\F$, and is therefore essential in $\pi_1(M)$ by
theorem~\ref{Palmeira_theorem}. So every torus bounds a solid torus, which is
necessarily on the $\fG_i$ side, and therefore each $\fG_i$ is a solid torus.

For distinct $i,j$, the core of a complementary annulus to the interstitial annulus
is a longitude contained in both $\partial \fG_i$ and $\partial \fG_j$. It follows
that the inclusion of each $\fG_i$ into $\fG_{i+1}$ is a homotopy equivalence, and therefore
the union is an open solid torus. Moreover, this inclusion takes interstitial regions
to interstitial regions, and therefore each interstitial region is of the form
$S^1 \times I \times \R^+$. It follows that each complementary region is a finite--sided
ideal polygon bundle over $S^1$, as claimed.

If some complementary region is actually a bigon bundle over $S^1$, after lifting
to $\til{M}$, the boundary leaves
intersect each leaf of $\til{\F}$ in quasigeodesics which are asymptotic at infinity.
It follows that such regions arose by unnecessarily splitting open a leaf of
$\til{\Lambda}^\pm_\geo$; we blow such regions down, identifying their boundary
leaves.

To see that $\Lambda^\pm_\split$ are minimal, after possibly removing finitely many
isolated leaves, observe that if $\Lambda$ is a minimal sublamination of $\Lambda^+_\split$,
then the construction in lemma~\ref{horizontally_foliated} still applies, and
therefore by lemma~\ref{straighten_interstitial_region}, the complementary regions
of this minimal sublamination are also finite sided ideal polygon bundles over $S^1$.
It follows that every leaf of $\Lambda^+_\split \backslash \Lambda$ must be a
suspension of one of the finitely many diagonals of one of the finitely many 
ideal polygons. The theorem follows.
\end{proof}

\subsection{Anosov flows}

In this subsection we address the question of when the laminations $\Lambda^\pm_\split$
are genuine, and not merely essential, at least when $\F$ has $2$--sided branching. 

\begin{lem}\label{endpoints_dense}
Let $S^1_\u$ be a minimal universal circle.
The endpoints of leaves of $\Lambda^+_\u$ are dense in $S^1_\u$, and similarly for
$\Lambda^-_\u$.
\end{lem}
\begin{proof}
Suppose not. Then there is some interval $I \subset S^1_\u$ which does not intersect
any leaf of $\Lambda^+_\u$. Since $S^1_\u$ is minimal, there is some leaf $\lambda$ of
$\til{\F}$ such that $\core(\phi_\lambda)$ intersects the interior of $I$. It follows
that $\phi_\lambda(I)$ is an interval in $S^1_\infty(\lambda)$ which does not
intersect a leaf of $\phi_\lambda(\Lambda^+_\u)$. But this is contrary to the fact
from theorem~\ref{complementary_regions_polygon_bundles} that
complementary regions to $\Lambda^+_\split$ are finite sided ideal polygon bundles over
$S^1$.
\end{proof}

It is clear from theorem~\ref{complementary_regions_polygon_bundles} that
$\Lambda^\pm_\split$ are genuine iff for some leaf $\lambda$ of $\til{\F}$, the
geodesic laminations $\Lambda^\pm_\geo(\lambda)$ are not foliations.

\begin{exa}\label{universal_fan}
Suppose $p \in S^1_\u$ is invariant under the action of $\pi_1(M)$ on $S^1_\u$.
We let $\Lambda_p$ be the lamination of $S^1_\u$ consisting of all unordered pairs
$p,q$ where $q \in S^1_\u \backslash p$. By construction~\ref{circle_geodesic} this
corresponds to the geodesic lamination of $\H^2$ by all geodesics with one endpoint
at $p$.

Note that for each leaf $\lambda$ of $\til{\F}$, the pushforward lamination
$\phi_\lambda(\Lambda_p)_\geo$ consists of the geodesic lamination of $\lambda$ 
by all geodesics
with one endpoint at $\phi_\lambda(p)$.

We call the foliation of $\H^2$ by all geodesics with one 
endpoint at some $p \in S^1_\infty$
the {\em geodesic fan} centered at $p$. By abuse of notation, we also refer to the
corresponding lamination of $S^1_\infty$ as the geodesic fan centered at $p$.
\end{exa}

At first glance, it appears as though example~\ref{universal_fan} is 
the typical example of a taut foliation for which $\Lambda^+_\split$ is
essential but not genuine.
However, this is somewhat misleading, as we will shortly see.

\begin{lem}\label{no_invariant_point}
Suppose for some $\F$ that there is a point
$p \in S^1_\u$ which is invariant under the action of $\pi_1(M)$ on $S^1_\u$.
Then $\F$ does not have two--sided branching.
\end{lem}
\begin{proof}
Let $\lambda$ be some leaf of $\til{\F}$, and let $\mu_1 > \lambda$ and
$\nu_2,\nu_3 > \lambda$ leaves such that $\mu_1$ is incomparable with $\nu_2,\nu_3$
and $\nu_2,\nu_3$ are incomparable with each other. Such leaves can certainly
be found if $\til{\F}$ branches in the positive direction. Since $\F$ is
taut, there is a positive transversal from the projection to $M$ of each $\nu_i$ to
the projection of $\mu_1$. Lifting to $\til{M}$, there exist elements 
$\alpha_2,\alpha_3 \in \pi_1(M)$ such that $\alpha_i(\mu_1) = \mu_i > \nu_i$
for $i = 2,3$. Then the $\mu_i$ are all translates of each other, are mutually
incomparable, and are all $> \lambda$ for some $\lambda$.

It follows that $L^+(\mu_i)$ are disjoint and incomparable
for $i = 1,2,3$ and therefore $\core(L^+(\mu_i))$ is contained in the closure of
a single gap of $\core(L^+(\mu_j))$ for $i \ne j$. But this implies that
the sets $\core(L^+(\mu_i))$ do not contain a common point of intersection.
Since $p$ is preserved by the action of $\pi_1(M)$, if it were contained
in $\core(L^+(\mu_i))$ for some $i$, it would be contained in 
$\core(L^+(\mu_i))$
for all three, contrary to what we have just shown. It follows that $p$
is not contained in $\core(L^+(\mu_1))$, and therefore is not contained in
$\core(L^+(\alpha(\mu_1)))$ for {\em any} $\alpha \in \pi_1(M)$.

But if $\til{\F}$ branches in the negative direction, by the tautness of
$\F$ we can find an element $\beta \in \pi_1(M)$ such that
$\beta(\mu_1)$ and $\mu_1$ are incomparable, and both satisfy
$\mu_1,\beta(\mu_1) < \lambda'$ for some $\lambda'$. But then
the union of $L^+(\mu_1)$ and $L^+(\beta(\mu_1))$ is all of $L$, and therefore
$\core(L^+(\mu_1)) \cup \core(L^+(\beta(\mu_1))) = S^1_\u$, 
so $p$ is contained in one of them, which is a contradiction.

It follows that $\F$ does not have $2$--sided branching, as claimed.
\end{proof}

\begin{construct}\label{pair_of_points}
Let $\G$ be a foliation of $\H^2$ by geodesics, and suppose $\G$ is not a geodesic
fan. Then $\G$ does not branch, and the leaf space of $\G$ is homeomorphic to $\R$.
Corresponding to the two ends of $\R$ there are exactly two points in $S^1_\infty$ 
which are not the endpoints of any leaf of $\G$. Call these the {\em ideal leaves} of $\G$.
\end{construct}

\begin{lem}\label{invariant_alternative}
Let $\F$ be a taut foliation of an atoroidal $3$--manifold $M$, 
and suppose for every leaf $\lambda$ of
$\til{\F}$ that $\Lambda^+_\geo(\lambda)$ is a foliation. Then every foliation
$\Lambda^+_\geo(\lambda)$ is a geodesic fan centered at some unique $s(\lambda) \in
S^1_\infty(\lambda)$.
\end{lem}
\begin{proof}
Let $J \subset L$ be the leaves of $\til{\F}$ for which $\Lambda^+_\geo(\lambda)$ is not
a geodesic fan. Then by construction~\ref{pair_of_points}, to each 
$\lambda \in J$ we can associate two points $p^\pm(\lambda) \in S^1_\infty(\lambda)$ 
which are the ideal leaves of the foliation $\Lambda^+_\geo(\lambda)$. 
Let $\gamma_\lambda$ 
be the geodesic spanned by $p^\pm$.

By lemma~\ref{continuity_of_image_lamination}, $J$ is open as a subset of $L$, and
the union 
$$\til{G} = \bigcup_{\lambda \in J} \gamma_\lambda$$
is a locally finite union of complete planes transverse to $\til{\F}$. This
union covers some compact surface $G \subset M$ transverse to $\F$, and the
intersection with leaves of $\F$ defines a foliation of $G$. It follows that $G$
consists of a union of incompressible tori and Klein bottles. But $M$ is
atoroidal, so $J$ is empty. The lemma follows.
\end{proof}

To characterize those taut foliations for which $\Lambda^\pm_\split$ are essential
but not genuine, we must introduce the notion of an {\em Anosov flow}.

\begin{defn}
An {\em Anosov flow} $\phi_t$ on a $3$--manifold $M$, with orbit space the $1$--dimensional
foliation $X$, is a flow which preserves a
decomposition of the tangent bundle
$$TM = E^s \oplus E^u \oplus TX$$
and such that the time $t$ flow uniformly expands $E^u$ and contracts $E^s$.
That is, there are constants $\mu_0 \ge 1,\mu_1 > 0$ so that
$$\|d\phi_t(v)\| \le \mu_0 e^{-\mu_1 t}\|v\| \text{ for any } v \in E^s,t \ge 0$$
and
$$\|d\phi_{-t(v)}\| \le \mu_0 e^{-\mu_1 t}\|v\| \text{ for any } v \in E^u, t \ge 0$$
The $1$--dimensional foliations obtained by integrating $E^s$ and $E^u$ are 
called the {\em strong stable} and {\em strong unstable} foliations, and
we denote them $X^{ss},X^{su}$ respectively. Furthermore, the bundles $TX \oplus E^s$
and $TX \oplus E^u$ are integrable, by a theorem of Anosov \cite{Anosov}, and
are tangent to $2$--dimensional foliations called the {\em weak stable} and
{\em weak unstable} foliations, denoted $\F^{ws},\F^{wu}$ respectively.
\end{defn}

See \cite{sF94} and \cite{sF97} for more details, and some important results in the
theory of Anosov flows on $3$--manifolds.

\begin{exa}
Let $\F$ be $\F^{ws}$ for some Anosov flow $X$ on $M$. Then every leaf of
$\til{\F}$ is foliated by flowlines of $X$. Suppose flowlines of $X$ are quasigeodesic
in leaves of $\til{\F}$. Then after straightening flowlines leafwise, the foliation
of each leaf $\lambda$ by $X$ is a geodesic fan, asymptotic to some 
$p \in S^1_\infty(\lambda)$.
\end{exa}

The main result of this subsection is that such examples are the only possibility,
when $\Lambda^\pm_\split$ are essential but not genuine.

\begin{thm}\label{foliation_is_stable_Anosov}
Let $\F$ be a taut foliation of $M$, and suppose that $\Lambda^+_\split$ is
essential but not genuine. Then there is an Anosov flow $\phi_t$ of $M$
such that $\F$ is the weak stable foliation of $\phi_t$,
and $\Lambda^+_\split$ is the weak unstable foliation.
\end{thm}
\begin{proof}
Constructing the flow is easy; most of the proof will be concerned with verifying
that it has the requisite properties. 

By lemma~\ref{invariant_alternative}, for every leaf $\lambda$ of $\til{\F}$, the
lamination $\Lambda^+_\geo(\lambda)$ is a geodesic fan, asymptotic to some
unique $p(\lambda) \in S^1_\infty(\lambda)$.

Let $\til{Y}$ be the unit length
vector field on $\til{M}$ contained in $T\til{\F}$ which on a leaf 
$\lambda$ points towards $p(\lambda) \in S^1_\infty(\lambda)$. Here we are identifying
$UT_p\lambda$ with $S^1_\infty(\lambda)$ for each $p \in \lambda$ by the
endpoint map $e$. Then $\til{Y}$ descends to a nowhere vanishing leafwise geodesic
vector field $Y$ on $\F$. We will show that if $\phi_t$ denotes the time $t$
flow generated by $Y$, then $\phi_t$ is an Anosov flow, and that
$\F$ is the weak stable foliation for $\phi_t$.

We define $E^u$ as follows. Each point $q \in S^1_\u$ determines a
geodesic $\gamma_q(\lambda)$ in each leaf $\lambda$ of $\til{\F}$ by setting
$\gamma_q(\lambda)$ equal to the unique geodesic from $\phi_\lambda(q)$ to
$p(\lambda)$. As we let $\lambda$ vary but fix $q$, the $\gamma_q(\lambda)$
sweep out a (possibly disconnected) union of planes transverse to $\til{\F}$,
whose leaves intersect leaf of $\til{\F}$ exactly in the flowlines of $\til{Y}$.
We define $E^u$ to be the orthogonal distribution to $\til{Y}$ in the
tangent space to these leaves. 

For simplicity, we first treat the case that $\F$ is minimal, and then show how
to modify the argument for general $\F$.

Recall definition~\ref{sawblade_definition} of a {\em sawblade} 
from \S~\ref{marker_subsection},
and the definition of the {\em contracting} and {\em expanding} directions.

Let $\gamma$ be a closed embedded geodesic contained in a leaf $\lambda$ of $\F$.
Let $\til{\lambda}$ be a covering leaf of $\lambda$ in $\til{M}$, stabilized
by the corresponding element $[\gamma] \in \pi_1(M)$. Since
$p(\til{\lambda})$ is defined intrinsically by the foliation 
$\Lambda^+_\geo(\til{\lambda})$,
it follows that $[\gamma]$ fixes $p(\til{\lambda})$, and $\gamma$ is a closed
orbit of $Y$. Let $\til{\gamma}$ be the corresponding axis of $[\gamma]$ on
$\til{\lambda}$. Then one endpoint of $\til{\gamma}$ is $p(\til{\lambda})$;
let $r \in S^1_\infty(\til{\lambda})$ be
the other endpoint of $\til{\gamma}$. By hypothesis, $\til{\gamma}$ is
equal to $\phi_{\til{\lambda}}(l)_\geo$ for some leaf $l$ of $\Lambda^+_\u$.

Let $\tau$ is an embedded interval in $L$ containining $\til{\lambda}$
as an endpoint. For sufficiently short $\tau$, the intervals
$\tau$ and $[\gamma](\tau)$ are completely comparable; moreover, for some
choice of orientation on $\gamma$, we can assume $[\gamma](\tau) \subset \tau$.

Then the set of leaves $\phi_{\nu}(l)_\geo$ for $\nu \in \tau$ is an
embedded rectangle $R$ in $\til{M}$, such that $\gamma(R) \subset R$, and
we can find an embedded
$\epsilon$--sawblade $P$ for $\F$ in $M$ with $\gamma$ as a boundary circle. 
Notice that $R$ is tangent to $E^u\oplus TY$.
Since $P$ is embedded, there
is a lower bound on the length of an arc in $M$ from $P$ to itself which is
not homotopic into $P$. By minimality of $\F$, there is a uniform $R$ such
that for any leaf $\lambda$ of $\F$, and every point $p \in \lambda$,
the ball of radius $R$ about $p$ in $\lambda$ (in the path metric) intersects
$P$.

Return to the universal cover. Then preimages of $P$
intersect every leaf $\lambda$ of $\til{\F}$ in a union of bi--infinite geodesics
and geodesic rays, contained in flowlines of $\til{Y}$,
which intersect the ball of radius $R$ about any point in
$\lambda$, as measured in $\lambda$. If $\til{P}$ is one component of the
preimage, then $\til{P} \cap \lambda$ contains a geodesic ray $\delta$ in
the contracting direction. Let $q \in \lambda$ be a point far
from the geodesic containing $\delta$. Then there is a translate
$\alpha(\til{P})$ with $\alpha \in \pi_1(M)$ which intersects the ball in
$\lambda$ of radius $R$ about $q$, and whose intersection with
$\lambda$ contains a geodesic ray $\delta'$. By the choice of $q$, the rays
$\delta$ and $\delta'$ are not contained in the same geodesic. Moreover, since
$P$ is compact and embedded in $M$, it does not accumulate on itself, so the
rays $\delta$ and $\delta'$ are not asymptotic to the same point in $S^1_\infty(\lambda)$.
On the other hand, $\delta$ and $\delta'$ are both contained in flowlines of
$\til{Y}$, which {\em are} asymptotic in the positive direction; it follows that
the contracting direction of $\gamma$ is the negative direction.
{\it A priori}, holonomy around the sawblade $P$ is merely non--increasing for
some nearby leaf. But in fact, this argument shows that the holonomy is actually
{\em strictly decreasing} for all leaves sufficiently close to $\gamma$.
The same argument shows that there is another sawblade $P'$ on the other side of
$\gamma$, for which $\gamma$ is also the contracting direction.

We show now that flow along $\til{Y}$ eventually strictly increases the length of
any integral curve of $E^u$. Let $\tau$ be a short integral curve, and
let $\til{\tau}$ be a lift to $\til{M}$. Let
$\lambda$ be a leaf of $\til{\F}$ which intersects $\til{\tau}$ and also some lift
$\til{\gamma}$ of $\gamma$. Then the flowline of $\til{Y}$ through $\tau \cap \lambda$
is eventually asymptotic to the (negatively oriented) $\til{\gamma}$. But we have
just seen that holonomy around $\gamma$ is strictly decreasing, so flow along $\til{Y}$
eventually blows up any arbitrarily short transversal to $\gamma$ to a transversal
of definite size. If follows that flow along $\til{Y}$ eventually blows up the
length of any arbitrarily short $\til{\tau}$, as claimed. By covering an integral
curve with such short curves, and using the compactness of $M$,
we can find uniform estimates for the rate of this blow up,
as required.

If $\F$ is not minimal, the argument is basically the same, except that we must
use the fact that almost every geodesic ray in a leaf of $\F$ 
is asymptotic to some minimal set to extend the arguments to all of $\F$.
\end{proof}

It is now easy to deduce our main theorem.

\begin{lam_exist_thm}
Let $\F$ be a co--orientable taut foliation of a closed, orientable
algebraically atoroidal $3$--manifold $M$. Then either $\F$ has $2$--sided branching
and is the weak stable foliation of an Anosov flow, 
or else there are a pair of very full
genuine laminations $\Lambda^\pm_\split$ transverse to $\F$.
\end{lam_exist_thm}
\begin{proof}
If $\F$ has two--sided branching, then this follows from 
theorem~\ref{laminations_essential}, theorem~\ref{complementary_regions_polygon_bundles} 
and theorem~\ref{foliation_is_stable_Anosov}.

If $\F$ has one--sided branching, this follows from
theorem 4.1.1 of \cite{dC00}. If $\F$ is $\R$--covered, this follows from
theorem 5.3.13 of \cite{dC99b}.
\end{proof}

\begin{cor}
Let $M$ be a closed $3$--manifold which admits a taut foliation.
Then either $M$ is toroidal, or admits an Anosov flow whose weak stable and
unstable foliations have two--sided branching,
or else $\pi_1(M)$ is word hyperbolic,
the mapping class group of $M$ is finite, and every self--homeomorphism
of $M$ homotopic to the identity is isotopic to the identity.
\end{cor}
\begin{proof}
This follows from theorem A together with theorem~\ref{lamination_word_hyperbolic},
theorem~\ref{lamination_finite_MCG} and theorem~\ref{homotopy_implies_isotopy},
all due to Gabai--Kazez.
\end{proof}

\section{The dynamics of $\Lambda^\pm_\split$}\label{lamination_dynamics}

This section is basically descriptive. No important lemmas or theorems are proved,
but we try to describe in some geometric detail the interaction of the
foliation $\F$ with the laminations $\Lambda^\pm_\split$. Finally, in 
\S~\ref{optimistic_picture} we pose some questions, and paint an optimistic
and conjectural picture of the interaction of $\Lambda^\pm_\split$ and $\F$,
at least in the best case.

\subsection{The structure of gut regions}\label{gut_region_structure}

In this subsection we describe the structure of a gut region of $\Lambda^\pm_\split$.
This subsection is basically descriptive; no theorems are proved here, and only one
straightforward lemma.
For convenience of notation, we concentrate on $\Lambda^+_\split$; of course,
the case of $\Lambda^-_\split$ is completely analogous. 

By theorem~\ref{complementary_regions_polygon_bundles} we know that complementary
regions are finite sided ideal polygon bundles over $S^1$, and therefore
gut regions are finite sided {\em neutered} ideal polygon bundles over $S^1$.
This terminology is standard: an ideal polygon is {\em neutered} by removing
neighborhoods of its ideal points. The usual neighborhoods one takes are the
intersections with small horoballs centered at the points.

Recall that the collection of cusps constructed in
lemma~\ref{image_is_branched_lamination} are the cores of a system of interstitial
annuli $A_i$ for a decomposition of complementary regions 
of $\Lambda^\pm_\split$ into guts and interstitial $I$--bundles, and that
there is a corresponding decomposition of $\Lambda^\pm_\geo$ into
guts and interstitial $I$--bundles, where some of the $I$--bundles
might limit on branch circles or lines. Throw away annuli
bounding compact $I$--bundles. Now observe that each
$\til{A_i}$ intersects each leaf $\lambda$ of $\til{\F}$ in at most two intervals,
and therefore by lemma~\ref{one_sided_branching}, we can conclude that the
$A_i$ are actually {\em transverse} to $\F$. Note {\it a posteriori} that we can
conclude from theorem~\ref{complementary_regions_polygon_bundles} that every
interstitial annular system can be isotoped to be transverse to $\F$; that is,
we need to allow components which are tangential to $\F$, 
as is proved in lemma~\ref{straighten_interstitial_region}.

Recall further that the branch locus of $\Lambda^\pm_\geo$ consists 
of a collection of geodesic lines and circles in leaves of $\F$. Some of these
might be in the closure of an interstitial region of $\Lambda^\pm_\geo$.

Let $\fG$ be a gut region of $\Lambda^+_\split$. As remarked above, it 
is bounded by a system of interstitial annuli 
which are transverse to $\F$. Let $\til{\fG}$ be a cover of $\fG$ 
in $\til{M}$. Then $\til{\fG}$ is also the universal cover of $\fG$.
Topologically, $\fG$ is a solid torus, and $\til{\fG}$ is a solid cylinder.
If $\gamma$ denotes the core circle of $\fG$, we can also think of $\gamma$ by
abuse of notation as the generator of $\pi_1(\fG) = \Z$ which acts on $\til{\fG}$
by deck transformations. Let $\lambda_t$, with $t \in (-\infty,\infty)$ 
parameterize the leaves of $\til{\F}$ which intersect $\til{\fG}$, and 
suppose this parameterization is chosen so that the action of $\gamma$ 
on $L$ satisfies $$\gamma(\lambda_t) = \lambda_{t+1}$$

The boundary $\partial \fG$ consists of two parts: the annular components $A_i \subset
\partial \fG$, and the {\em laminar boundary} $\partial \fG \cap \Lambda^+_\split$.
Note that if $\Lambda^\pm$ are co-orientable, this decomposition defines the structure of a 
{\em sutured manifold} on $\fG$, in
the sense of Gabai; see \cite{Gabai_fol} for a definition and basic properties of
sutured manifolds.
We denote these subsets of $\partial \fG$ by $\partial_v\fG$ and $\partial_h\fG$,
consistent with the usual notation from \cite{Gabai_fol}. These lift
to $\partial_v \til{\fG}$ and $\partial_h \til{\fG}$ in the obvious way.
The boundary $\partial \til{\fG}$ is foliated by circles of intersection with
leaves of $\til{\F}$.

There are three distinct classes of interstitial regions. Recall the map
$\psi$ from lemma~\ref{image_is_branched_lamination}.

\begin{figure}[ht]
\centerline{\relabelbox\small \epsfxsize 5.0truein
\epsfbox{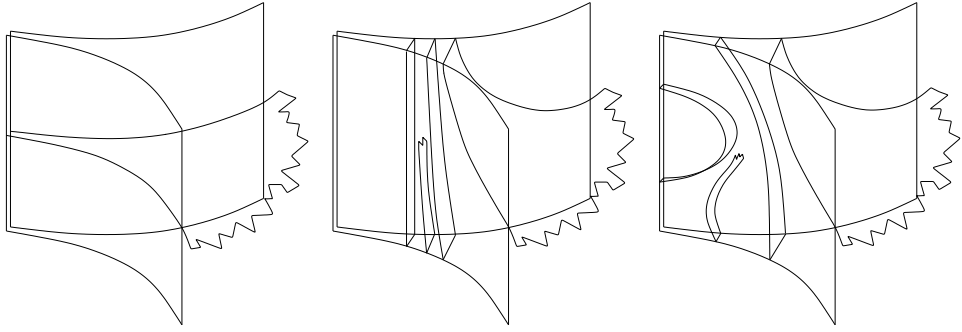}
\endrelabelbox}
\caption{Three different kinds of interstitial region.}
\label{interstices}
\end{figure}

\begin{enumerate}
\item{If an interstitial region $R$ of $\Lambda^+_\geo$ contains no branch locus, that
is, if the corresponding interstitial region $R'$ of $\Lambda^+_\split$ maps 
homeomorphically to $R$ by $\psi$, 
then the foliation of $\partial \fG$ by circles transverse to $\F$ 
can be extended to the entire interstitial region. The lift of the
interstitial region intersects exactly the leaves $\lambda_t$ of $\til{\F}$; i.e. the
same set of leaves that $\til{\fG}$ intersects. We call this kind of interstitial
region a {\em cusp}.}
\item{If an interstitial region $R$ of $\Lambda^+_\geo$ contains a circle branch component
$\fc$, this circle gets split open to a tangential interstitial annulus in the
corresponding interstitial region $R'$ of $\Lambda^+_\split$. The leaves of
$\F \cap R'$ spiral around this annulus and limit on to it. In $\til{M}$, the annulus
is covered by a rectangle contained in a leaf $\nu$ of $\til{\F}$ which is a limit
of $\lambda_t$ as either $t \to -\infty$ or $t \to \infty$. In the figure, the
spiralling is from the positive side, and $\nu$ is a limit of $\lambda_t$ as
$t \to -\infty$. We call this kind of interstitial region a {\em (positive or negative)
annular spiral}.}
\item{If an interstitial region $R$ of $\Lambda^+_\geo$ contains a line branch component,
it might conceivably contain infinitely many such components, which we denote $\gamma_i$. Each of these gets split
open to a rectangle contained in the interior of an interstitial region $R'$
of $\Lambda^+_\split$ bounded by a transverse interstitial annulus. Moreover, the
leaves $\F \cap R$ spiral out to fill all of the preimage of $R$ under
the collapsing map $\psi$, and limit on the union of split open rectangles in $R'$
contained in distinct leaves $\nu_i$ of $\til{\F}$ which are all limits of
$\lambda_t$ as $t \to -\infty$.
Note that there is no claim that the $\nu_i$ fall into finitely many orbit classes
under the action of $\gamma$. Note that since there are at least infinitely many
$\nu_i$ which are limits of $\lambda_t$, the spiralling in this case {\em must} be
from the positive side, by lemma~\ref{one_sided_branching}. We call this kind of
interstitial region a {\em positive rectangular spiral}.}
\end{enumerate}

These three classes of interstitial region are 
illustrated in figure~\ref{interstices}. 

We say that two interstitial regions bounding a gut component of $\Lambda^+_\split$
are {\em adjacent} if they have a common boundary leaf.

\begin{lem}\label{annular_spirals}
Adjacent interstitial regions bounding a gut component of $\Lambda^+_\split$ do
not both contain negative annular spirals.
\end{lem}
\begin{proof}
If so, then there are at least two distinct leaves $\nu,\nu'$ of $\til{\F}$
which are positive limits of $\lambda_t$ as $t \to \infty$, and which intersect the same 
leaf of $\til{\Lambda^+_\split}$,
contrary to lemma~\ref{one_sided_branching}. 
\end{proof}

\subsection{Dynamics of $\Lambda^\pm_\split$}

We continue to use the setup and notation from subsection~\ref{gut_region_structure}.
Throughout this subsection, for convenience, we work with $\til{\Lambda^+_\geo}$ instead
of $\til{\Lambda^+_\split}$, since the relationship with the geometry of
leaves of $\til{\F}$ is clearer. By abuse of notation, we refer to the
{\em leaves} of $\til{\Lambda^+_\geo}$, by which we mean the images of leaves
of $\til{\Lambda^+_\split}$ under the monotone map $\psi$ from 
lemma~\ref{image_is_branched_lamination}. The relationship between leaves of
$\til{\Lambda^+_\geo}$ and leaves of $\Lambda^+_\u$ is more straightforward:
each leaf $l_\u$ of $\Lambda^+_\u$ determines a geodesic $l_\u(\lambda)$ in each leaf
$\lambda$ of $\til{\F}$, as a leaf of 
$\Lambda^+_\geo(\lambda) = (\phi_\lambda(\Lambda^+_\u))_\geo$. The union of
the geodesics $l_\u(\lambda)$ as $\lambda$ varies over $L$ is a union of leaves
of $\til{\Lambda^+_\geo}$. Moreover, every leaf of $\til{\Lambda^+_\geo}$ arises this
way, although possibly not uniquely. For a leaf $l$ of $\til{\Lambda^+_\geo}$,
we let $l_\u$ denote a leaf of $S^1_\u$ associated to $l$ by this construction.
In particular, $$l \cap \kappa = (\phi_\kappa(l_\u))_\geo$$
for every leaf $\kappa$ of $\til{\F}$ which $l$ intersects.

\begin{figure}[ht]
\centerline{\relabelbox\small \epsfxsize 3.0truein
\epsfbox{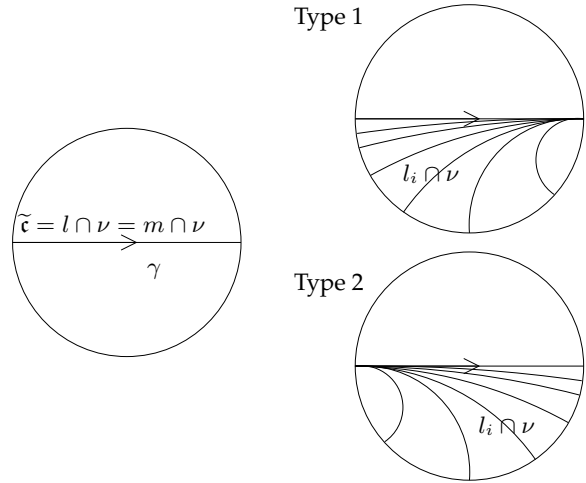}
\adjustrelabel <-12pt, 0pt> {l=m}{$\til{\fc} = l \cap \nu =m \cap \nu$}
\adjustrelabel <0pt, -15pt> {g}{$\gamma$}
\adjustrelabel <-5pt, 0pt> {case1}{Type 1}
\adjustrelabel <-5pt, 0pt> {case2}{Type 2}
\relabel {li}{$l_i \cap \nu$}
\adjustrelabel <-20pt, 0pt> {li2}{$l_i \cap \nu$}
\endrelabelbox}
\caption{Either the attracting or the repelling fixed point of $\gamma$ on $\nu$ 
is asymptotic to nearby leaves of $\Lambda^+(\nu)$ on either side.}
\label{branch_circle}
\end{figure}

Suppose $l,m$ are boundary leaves of $\til{\Lambda^+_\geo}$ which both bound a common
annular interstitial region after splitting, and 
let $\fc$ denote the corresponding branch circle of
$\Lambda^+_\geo$. Let $l_\u,m_\u$ denote associated leaves of $\Lambda^+_\u$. 
Let $l_i,m_i$ be leaves of $\til{\Lambda^+_\geo}$ which
accumulate on $l,m$ respectively, with associated leaves $(l_i)_\u,(m_i)_\u$ of
$\Lambda^+_\u$, and let $\nu$ be the leaf of $\til{\F}$ containing
the lift $\til{\fc}$, which is a geodesic in $\nu$, and a leaf of
$\Lambda^+_\geo(\nu)$. Note that we have an equality
$$\til{\fc} = l \cap \nu = m \cap \nu$$
The element $\gamma \in \pi_1(M)$ stabilizes $\nu$, and
therefore acts as a translation on $\til{\fc}$. Since $l_i \to l$, we must have
$l_i \cap \nu \to l \cap \nu$. On the other hand, $\gamma$ stabilizes
$\nu$, and permutes $l_i \cap \nu$. Since this sequence of geodesics accumulates
on $l \cap \nu$, they must all share a common endpoint with $\til{\fc}$,
or else some $l_i \cap \nu$ would cross $\gamma(l_i \cap \nu)$ transversely,
contrary to the definition of a lamination. There are two possibilities:
either the $l_i \cap \nu$ share the attracting fixed point of $\gamma$ with
$l \cap \nu$, or else they share the repelling fixed point. We call these
type 1 and type 2 respectively; these types are illustrated in figure~\ref{branch_circle}

The projection $\pi(l)$ of the leaf $l$ to $M$ is an annulus which is a boundary leaf of
$\Lambda^+_\geo$. The projection of the leaf $l_i$ spirals around this annulus
under holonomy transport around the core of the annulus. In type 1, the
projection of $l_i$ accumulates on the projection of $l$ in the positive direction,
and in type 2, the accumulation is in the negative direction. That is, holonomy
transport around the core of the annulus $\pi(l)$ is contracting in type 1,
and expanding in type 2. Here we are implicitly orienting
the core of the annulus so that the {\em positive} direction agrees with the
co--orientation of $\F$.

\begin{figure}[ht]
\centerline{\relabelbox\small \epsfxsize 4.0truein
\epsfbox{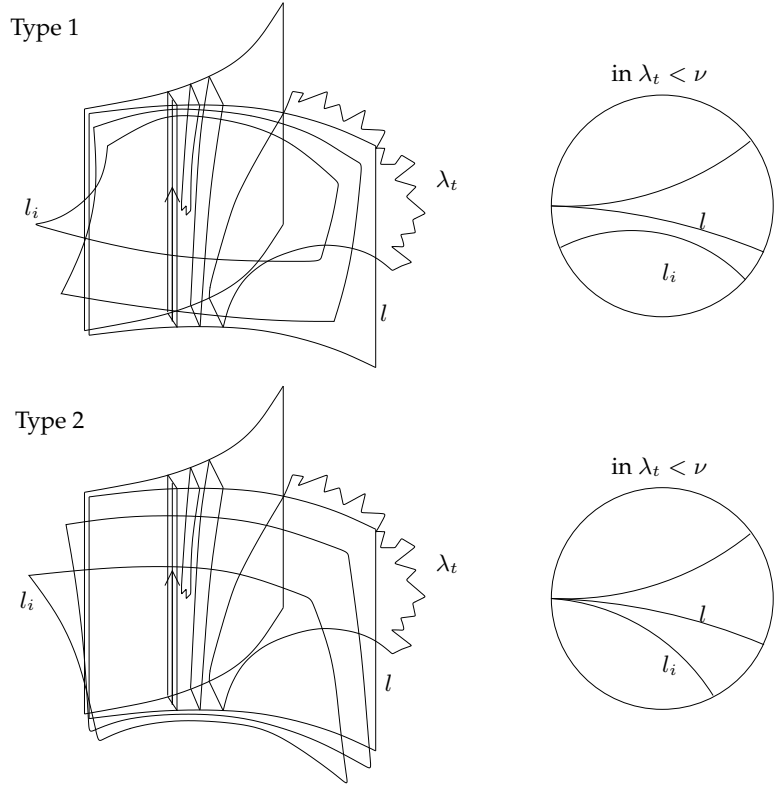}
\relabel {case1}{Type 1}
\relabel {case2}{Type 2}
\relabel {l}{$l$}
\relabel {l2}{$l$}
\relabel {l3}{$l$}
\relabel {l4}{$l$}
\relabel {li}{$l_i$}
\relabel {li2}{$l_i$}
\relabel {li3}{$l_i$}
\relabel {li4}{$l_i$}
\relabel {lt}{$\lambda_t$}
\relabel {lt1}{$\lambda_t$}
\relabel {inlt}{in $\lambda_t < \nu$}
\relabel {inlt2}{in $\lambda_t < \nu$}
\endrelabelbox}
\caption{In type 1, nearby leaves $l_i$ to $l$ spiral around $l$ in the positive direction;
in type 2, they spiral around in the negative direction. This determines the
geometry of $l_i \cap \lambda_t$ for large $\lambda_t$. In particular, $l_i \cap \lambda_t$
is asymptotic to $l \cap \lambda_t$ in type 2, but {\em not necessarily} in type 1.}
\label{branch_spiral}
\end{figure}

As described in \S~\ref{gut_region_structure},
the intersections $\lambda_t \cap \til{\fG}$ spiral around 
$\nu$, accumulating on the geodesic $\til{\fc}$. It follows that if our
annular spiral is {\em positive}, then in type 2, the geodesics
$l_i \cap \lambda_t$ and $l \cap \lambda_t$ are asymptotic. Similarly, if
the annular spiral is {\em negative}, then in type 1, the geodesics
$l_i \cap \lambda_t$ and $l \cap \lambda_t$ are asymptotic. This is illustrated in
figure~\ref{branch_spiral}.

Figure~\ref{branch_spiral} illustrates types 1 and 2 for a positive annular spiral.
Notice that although $l \cap \lambda_t$ and $l_i \cap \lambda_t$ must be asymptotic
in type 2, they are not necessarily asymptotic in type 1.

\subsection{An optimistic picture}\label{optimistic_picture}

The constructions and results above lead to many other natural questions and potential
connections with other areas of foliation theory and $3$--manifold topology. We discuss
some of these connections now.

{\noindent \bf Universal circles for $\Lambda^\pm_\split$:}
If $\Lambda^\pm_\split$ are genuine and tight, then theorem 3.8 of \cite{dCnD02} shows
that there are a pair of circles $(S^1_\u)^\pm$ and laminations $\mathscr{L}^\pm$
of $(S^1_\u)^\pm$ respectively whose leaves are in bijective correspondence with the
leaves of $\til{\Lambda^\pm_\split}$. Moreover, there is a natural action of
$\pi_1(M)$ on either circle, preserving the laminations $\mathscr{L}^\pm$.
In this case, what is the relationship between the three circles 
$S^1_\u,(S^1_\u)^+,(S^1_\u)^-$?

{\noindent \bf Good dynamical pairs for $\Lambda^\pm_\split$:}
If $\Lambda^\pm_\split$ are carried by a pair of branched surfaces which are a
{\em good dynamic pair} in the sense of Mosher (see \cite{lM00} for a definition), 
then by theorem 4.10.4 of \cite{lM00}, $M$ admits a pseudo--Anosov flow $\psi_t$.
Can $\psi_t$ be made transverse or almost transverse to $\F$? Corollary 3.9 of
\cite{dCnD02} constructs a universal circle associated to a pseudo--Anosov flow; what
is the relationship of this universal circle with the other universal circles above?

{\noindent \bf Finite depth foliations:}
If $\F$ is finite depth, Mosher and Gabai independently construct a pseudo--Anosov flow
on $M$ almost transverse to $\F$. The leaves of the singular stable and unstable foliations
of the flow can be split open to a pair of genuine laminations. What is the relationship
of these laminations to $\Lambda^\pm_\split$?

{\noindent \bf Anosov flows:}
If $\F$ is the weak stable foliation of an Anosov flow, is there nevertheless a universal
circle for $\F$ for which $\Lambda^\pm_\split$ are genuine? If $\Lambda^+_\split$
is (monotone equivalent to ) a taut foliation, what is the relationship between its
universal circle, and the universal circle of $\F$?

In the most optimistic picture for the structure of $\F,\Lambda^\pm_\split$, all circles
and structures constructed by various means are compatible to the extent that this makes sense.
This motivates the definition of the following object, called a {\em pseudo--Anosov package} 
for $M$.

\begin{defn}
Let $M$ be an algebraically atotoroidal $3$--manifold, 
and $\F$ a taut foliation. A {\em pseudo--Anosov package}
$\Psi$ for $M$ consists of the following structure:
\begin{enumerate}
\item{The universal cover $\til{M}$ has the structure of a product $D^2_\u \times \R$ in
such a way that the action of $\pi_1(M)$ descends to an action on the $D^2_\u$ factor.
The action of $\pi_1(M)$ on $\til{M}$ extends to a continuous action on the
product of the {\em closed} disk with $\R$. Moreover, there should be a natural
identification $\partial D^2_\u = S^1_\u$, compatible with the representation $\rho_\u$.}
\item{The laminations $\Lambda^\pm_\u$ determine geodesic laminations $(\Lambda^\pm_\u)_\geo$
of $D^2_\u$ with finite sided complementary regions which are transverse to each other and
{\em bind} $D^2_\u$; that is, complementary regions to the union of the two geodesic laminations
are finite sided compact polygons.}
\item{For each leaf $\lambda$ of $\til{\F}$, the projection of $\lambda$ to $D^2_\u$ is
{\em convex}; i.e. it is the interior of a region bounded by an embedded union of complete
geodesics of $(\Lambda^\pm_\u)_\geo$}
\item{The laminations $\Lambda^\pm_\split$ can be blown down to the stable/unstable foliations
of a pseudo--Anosov flow $\psi_t$ transverse, or almost transverse to $\F$}
\item{The quotient of the circle $S^1_\u$ by the laminar relations determined by
$\Lambda^\pm_\u$ is topologically a sphere $S^2_\u$, which is naturally homeomorphic to
the ideal boundary $S^2_\infty$ of $\til{M}$, and determines a Peano map $P:S^1_\u \to S^2_\u$
which can be approximated by embeddings.}
\end{enumerate}
\end{defn}

\begin{figure}[ht]
\centerline{\relabelbox\small \epsfxsize 4.0truein
\epsfbox{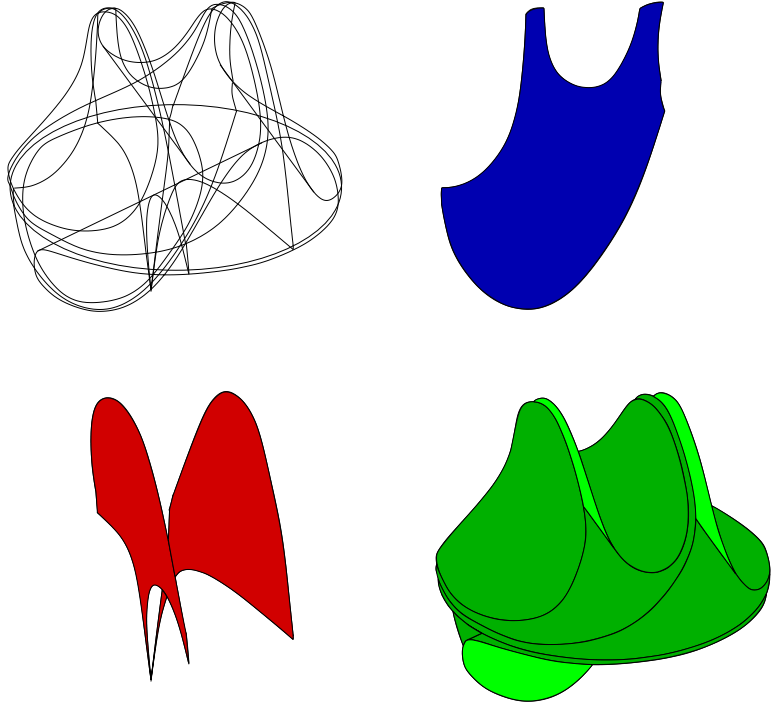}
\endrelabelbox}
\caption{Leaves of $\til{\Lambda^\pm_\split}$ (red and blue)
run up seams in either direction, and bind
the leaves of $\til{\F}$ (green)}
\label{splits}
\end{figure}

Properties 3 and 5 above imply that for each leaf $\lambda$ of $\til{\F}$, the
Peano map $P:S^1_\u \to S^2_\u$ factors through the monotone map 
$\phi_\lambda:S^1_\u \to S^1_\infty(\lambda)$ and induces a continuous map from $S^1_\infty(\lambda)$
to $S^2_\infty(\pi_1(M))$, continuously extending the inclusion $\lambda \to \til{M}$.
This so--called {\em continuous extension property} is actually established by Fenley for many
classes of taut foliations, including all those with {\em quasigeodesic} transverse
(or almost transverse) pseudo--Anosov flows. See \cite{Fenley_JAMS} for a detailed discussion.

{\noindent \bf Pseudo--Anosov packages:}
Given a pseudo--Anosov package $\Psi$, all the data of the package can be recovered from
the representation $\rho_\u:\pi_1(M) \to \homeo^+(S^1_\u)$ except for the foliation $\F$.
When are distinct taut foliations compatible with the same pseudo--Anosov package?
Are there only finitely many pseudo--Anosov packages (up to isotopy and
the ambiguity of $\F$) on a fixed $3$--manifold?

Figure~\ref{splits} shows how leaves of $\til{\Lambda^\pm_\split}$ and
of $\til{\F}$ should interlock. Notice how the blue leaf branches in the negative
direction, and the red leaves branch in the positive direction. Notice too how the
leaves of the laminations are asymptotic into the ``seams'' of $\til{\F}$.

\end{document}